\advance\hsize-3mm

\voffset=-11mm
\advance\vsize28mm

\def\pisac{Nenad Antoni\'{c} \& Marko Erceg \& Martin Lazar}
\def\naslova{Localisation principle for one-scale H-measures}

\def\draft{{\tt arXiv} version}
\def\dnum{dvostruko}


\catcode`\{=1 
\catcode`\}=2 
\catcode`\$=3 
\catcode`\&=4 
\catcode`\#=6 
\catcode`\^=7 \catcode`\^^K=7 
\catcode`\_=8 \catcode`\^^A=8 
\catcode`\^^I=10 
\chardef\active=13 \catcode`\~=\active 
\catcode`\^^L=\active \outer\def^^L{\par} 

\message{Preloading the plain format: codes,}




\catcode`@=11

\mathcode`\^^@="2201 
\mathcode`\^^A="3223 
\mathcode`\^^B="010B 
\mathcode`\^^C="010C 
\mathcode`\^^D="225E 
\mathcode`\^^E="023A 
\mathcode`\^^F="3232 
\mathcode`\^^G="0119 
\mathcode`\^^H="0115 
\mathcode`\^^I="010D 
\mathcode`\^^J="010E 
\mathcode`\^^K="3222 
\mathcode`\^^L="2206 
\mathcode`\^^M="2208 
\mathcode`\^^N="0231 
\mathcode`\^^O="0140 
\mathcode`\^^P="321A 
\mathcode`\^^Q="321B 
\mathcode`\^^R="225C 
\mathcode`\^^S="225B 
\mathcode`\^^T="0238 
\mathcode`\^^U="0239 
\mathcode`\^^V="220A 
\mathcode`\^^W="3224 
\mathcode`\^^X="3220 
\mathcode`\^^Y="3221 
\mathcode`\^^Z="8000 
\mathcode`\^^[="2205 
\mathcode`\^^\="3214 
\mathcode`\^^]="3215 
\mathcode`\^^^="3211 
\mathcode`\^^_="225F 
\mathcode`\ ="8000 
\mathcode`\!="5021
\mathcode`\'="8000 
\mathcode`\(="4028
\mathcode`\)="5029
\mathcode`\*="2203 
\mathcode`\+="202B
\mathcode`\,="613B
\mathcode`\-="2200
\mathcode`\.="013A
\mathcode`\/="013D
\mathcode`\:="303A
\mathcode`\;="603B
\mathcode`\<="313C
\mathcode`\=="303D
\mathcode`\>="313E
\mathcode`\?="503F
\mathcode`\[="405B
\mathcode`\\="026E 
\mathcode`\]="505D
\mathcode`\_="8000 
\mathcode`\{="4266
\mathcode`\|="026A
\mathcode`\}="5267
\mathcode`\^^?="1273 


\sfcode`\)=0 \sfcode`\'=0 \sfcode`\]=0

\delcode`\(="028300
\delcode`\)="029301
\delcode`\[="05B302
\delcode`\]="05D303
\delcode`\<="26830A
\delcode`\>="26930B
\delcode`\/="02F30E
\delcode`\|="26A30C
\delcode`\\="26E30F

\chardef\@ne=1
\chardef\tw@=2
\chardef\thr@@=3
\chardef\sixt@@n=16
\chardef\@cclv=255
\mathchardef\@cclvi=256
\mathchardef\@m=1000
\mathchardef\@M=10000
\mathchardef\@MM=20000


\message{registers,}

\count10=21 
\count11=9 
\count12=9 
\count13=9 
\count14=9 
\count15=9 
\count16=-1 
\count17=-1 
\count18=3 
\count19=255 
\countdef\insc@unt=19 
\countdef\allocationnumber=20 
\countdef\m@ne=21 \m@ne=-1 
\def\wlog{\immediate\write\m@ne} 


\countdef\count@=255
\dimendef\dimen@=0
\dimendef\dimen@i=1 
\dimendef\dimen@ii=2
\skipdef\skip@=0
\toksdef\toks@=0


\outer\def\newcount{\alloc@0\count\countdef\insc@unt}
\outer\def\newdimen{\alloc@1\dimen\dimendef\insc@unt}
\outer\def\newskip{\alloc@2\skip\skipdef\insc@unt}
\outer\def\newmuskip{\alloc@3\muskip\muskipdef\@cclvi}
\outer\def\newbox{\alloc@4\box\chardef\insc@unt}
\let\newtoks=\relax 
\outer\def\newhelp#1#2{\newtoks#1#1\expandafter{\csname#2\endcsname}}
\outer\def\newtoks{\alloc@5\toks\toksdef\@cclvi}
\outer\def\newread{\alloc@6\read\chardef\sixt@@n}
\outer\def\newwrite{\alloc@7\write\chardef\sixt@@n}
\outer\def\newfam{\alloc@8\fam\chardef\sixt@@n}
\def\alloc@#1#2#3#4#5{\global\advance\count1#1by\@ne
  \ch@ck#1#4#2
  \allocationnumber=\count1#1%
  \global#3#5=\allocationnumber
  \wlog{\string#5=\string#2\the\allocationnumber}}
\outer\def\newinsert#1{\global\advance\insc@unt by\m@ne
  \ch@ck0\insc@unt\count
  \ch@ck1\insc@unt\dimen
  \ch@ck2\insc@unt\skip
  \ch@ck4\insc@unt\box
  \allocationnumber=\insc@unt
  \global\chardef#1=\allocationnumber
  \wlog{\string#1=\string\insert\the\allocationnumber}}
\def\ch@ck#1#2#3{\ifnum\count1#1<#2%
  \else\errmessage{No room for a new #3}\fi}

\newdimen\maxdimen \maxdimen=16383.99999pt 
\newskip\hideskip \hideskip=-1000pt plus 1fill 
\newskip\centering \centering=0pt plus 1000pt minus 1000pt
\newdimen\p@ \p@=1pt 
\newdimen\z@ \z@=0pt 
\newskip\z@skip \z@skip=0pt plus0pt minus0pt
\newbox\voidb@x 

\outer\def\newif#1{\count@\escapechar \escapechar\m@ne
  \expandafter\expandafter\expandafter
   \edef\@if#1{true}{\let\noexpand#1=\noexpand\iftrue}%
  \expandafter\expandafter\expandafter
   \edef\@if#1{false}{\let\noexpand#1=\noexpand\iffalse}%
  \@if#1{false}\escapechar\count@} 
\def\@if#1#2{\csname\expandafter\if@\string#1#2\endcsname}
{\uccode`1=`i \uccode`2=`f \uppercase{\gdef\if@12{}}} 

\message{parameters,}


\pretolerance=100
\tolerance=200 
\hbadness=1000
\vbadness=1000
\linepenalty=10
\hyphenpenalty=50
\exhyphenpenalty=50
\binoppenalty=700
\relpenalty=500
\clubpenalty=150
\widowpenalty=150
\displaywidowpenalty=50
\brokenpenalty=100
\predisplaypenalty=10000
\doublehyphendemerits=10000
\finalhyphendemerits=5000
\adjdemerits=10000
\tracinglostchars=1
\uchyph=1
\defaulthyphenchar=`\-
\defaultskewchar=-1
\newlinechar=-1
\delimiterfactor=901
\showboxbreadth=5
\showboxdepth=3

\hfuzz=0.1pt
\vfuzz=0.1pt
\overfullrule=5pt
%
%
\hsize=162mm	
\vsize= 246mm	
\maxdepth=4pt
\splitmaxdepth=\maxdimen
\boxmaxdepth=\maxdimen
\delimitershortfall=5pt
\nulldelimiterspace=1.2pt
\scriptspace=0.5pt
\parindent=20pt

\parskip=0pt plus 1pt
\abovedisplayskip=12pt plus 3pt minus 9pt
\abovedisplayshortskip=0pt plus 3pt
\belowdisplayskip=12pt plus 3pt minus 9pt
\belowdisplayshortskip=7pt plus 3pt minus 4pt
\topskip=10pt
\splittopskip=10pt
\parfillskip=0pt plus 1fil

\thinmuskip=3mu
\medmuskip=4mu plus 2mu minus 4mu
\thickmuskip=5mu plus 5mu

\newskip\smallskipamount \smallskipamount=3pt plus 1pt minus 1pt
\newskip\medskipamount \medskipamount=6pt plus 2pt minus 2pt
\newskip\bigskipamount \bigskipamount=12pt plus 4pt minus 4pt
\newskip\normalbaselineskip \normalbaselineskip=12pt
\newskip\normallineskip \normallineskip=1pt
\newdimen\normallineskiplimit \normallineskiplimit=0pt
\newdimen\jot \jot=3pt
\newcount\interdisplaylinepenalty \interdisplaylinepenalty=100
\newcount\interfootnotelinepenalty \interfootnotelinepenalty=100

\def\magstephalf{1095 }
\def\magstep#1{\ifcase#1 \@m\or 1200\or 1440\or 1728\or 2074\or 2488\fi\relax}


\message{fonts,}


\def\@magscale#1{ scaled \magstep #1}
\def\@halfmag{ scaled \magstephalf}
\def\@ptscale#1{ scaled #100}

%
%
 \font\vrm  = cmr5    			
 \font\vmi  = cmmi5   			
    \skewchar\vmi ='177  		
 \font\vsy  = cmsy5   			
    \skewchar\vsy ='60 		
 \font\vmib = cmmib5   		
 \font\vbsy = cmbsy5  		
 \font\vbf  = cmbx5   			
 \font\vsf  = cmss8	\@ptscale5    			

\font\veufb = eufb5			
\font\veufm = eufm5			

\font\vmsam = msam5		
\font\vmsbm = msbm5		


%
%
%
 \font\virm  = cmr6    			
 \font\vimi  = cmmi6   			
    \skewchar\vimi ='177  		
 \font\visy  = cmsy6   			
    \skewchar\visy ='60 		
 \font\vimib = cmmib6   		
 \font\vibsy = cmbsy6  		
 \font\vibf  = cmbx6   			
 \font\visf  = cmss8	\@ptscale6    			


\font\vieufb = eufb6			
\font\vieufm = eufm6			

\font\vimsam = msam6		
\font\vimsbm = msbm6		


%
%
%
 \font\viirm  = cmr7    			
 \font\viimi  = cmmi7   			
    \skewchar\viimi ='177  		
 \font\viisy  = cmsy7   			
    \skewchar\viisy ='60 		
 \font\viiex  = cmex7   			
 \font\viimib = cmmib7   		
 \font\viibsy = cmbsy7  		
 \font\viibf  = cmbx7   			
 \font\viisf  = cmss8	\@ptscale7    			
 \font\viisfb = cmssbx10	\@ptscale7 		

\font\viieufb = eufb7			
\font\viieufm = eufm7			

\font\viimsam = msam7		
\font\viimsbm = msbm7		


%
%
%
 \font\viiirm  = cmr8    			
 \font\viiimi  = cmmi8   			
    \skewchar\viiimi ='177  		
 \font\viiisy  = cmsy8   			
    \skewchar\viiisy ='60 		
\font\viiiex  = cmex8   			
 \font\viiimib = cmmib8   		
 \font\viiibsy = cmbsy8  		
 \font\viiiit  = cmti8   			
 \font\viiisl  = cmsl8   			
 \font\viiibf  = cmbx8   			
 \font\viiitt  = cmtt8   			
    \hyphenchar\viiitt = -1         		
 \font\viiisf  = cmss8   			
 \font\viiisfb = cmssbx10	\@ptscale8 		


\font\viiieufb = eufb8			
\font\viiieufm = eufm8			

\font\viiimsam = msam8		
\font\viiimsbm = msbm8		


%
%
%
 \font\ixrm  = cmr9    			
 \font\ixmi  = cmmi9   		
    \skewchar\ixmi ='177  		
 \font\ixsy  = cmsy9   			
    \skewchar\ixsy ='60 		
\font\ixex  = cmex9   			
 \font\ixmib = cmmib9   		
 \font\ixbsy = cmbsy9  		
 \font\ixit  = cmti9   			
 \font\ixsl  = cmsl9   			
 \font\ixbf  = cmbx9   			
 \font\ixtt  = cmtt9   			
    \hyphenchar\ixtt = -1         		
 \font\ixsf  = cmss9   			
 \font\ixsfb = cmssbx10	\@ptscale9 		

\font\ixeufb = eufb9			
\font\ixeufm = eufm9			

\font\ixmsam = msam9		
\font\ixmsbm = msbm9		


%
%
%
 \font\xrm  = cmr10    			
 \font\xmi  = cmmi10   		
    \skewchar\xmi ='177  		
 \font\xsy  = cmsy10   			
    \skewchar\xsy ='60 		
 \font\xex  = cmex10   		
 \font\xmib = cmmib10   		
 \font\xbsy = cmbsy10  		
 \font\xit  = cmti10   			
 \font\xsl  = cmsl10   			
 \font\xbf  = cmbx10   			
 \font\xtt  = cmtt10   			
    \hyphenchar\xtt = -1         		
 \font\xsf  = cmss10   			
 \font\xsfb = cmssbx10 		

\font\xeufb = eufb10			
\font\xeufm = eufm10			

\font\xmsam = msam10		
\font\xmsbm = msbm10		


%
%
%
 \font\xirm  = cmr10        	\@halfmag		
 \font\ximi  = cmmi10       	\@halfmag		
    \skewchar\ximi ='177     				
 \font\xisy  = cmsy10      	\@halfmag		
    \skewchar\xisy ='60     				
 \font\xiex  = cmex10       	\@halfmag		
 \font\ximib = cmmib10       	\@halfmag		
 \font\xibsy = cmbsy10     	\@halfmag		
 \font\xiit  = cmti10      		\@halfmag		
 \font\xisl  = cmsl10       	\@halfmag		
 \font\xibf  = cmbx10       	\@halfmag		
 \font\xitt  = cmtt10      	\@halfmag		
    \hyphenchar\xitt = -1         				
 \font\xisf  = cmss10       	\@halfmag		
 \font\xisfb = cmssbx10     	\@halfmag		

\font\xieufb = eufb10    	\@halfmag		
\font\xieufm = eufm10    	\@halfmag		

\font\ximsam = msam10    	\@halfmag		
\font\ximsbm= msbm10    	\@halfmag		


%
%
%
 \font\xiirm  = cmr12    			
 \font\xiimi  = cmmi12   		
    \skewchar\xiimi ='177  		
 \font\xiisy  = cmsy10	\@magscale1   	
    \skewchar\xiisy ='60 		
 \font\xiiex  = cmex10	\@magscale1   	
 \font\xiimib = cmmib10\@magscale1   	
 \font\xiibsy = cmbsy10\@magscale1  	
 \font\xiiit  = cmti12	   		
 \font\xiisl  = cmsl12		   	
 \font\xiibf  = cmbx12		  	
 \font\xiitt  = cmtt12   			
    \hyphenchar\xiitt = -1         		
 \font\xiisf  = cmss12   			
 \font\xiisfb = cmssbx10\@magscale1 	

\font\xiieufb = eufb10	\@magscale1		
\font\xiieufm = eufm10	\@magscale1		

\font\xiimsam = msam10	\@magscale1		
\font\xiimsbm = msbm10	\@magscale1		


%
%
%
 \font\xivrm  = cmr12	\@magscale1    		
 \font\xivmi  = cmmi12	\@magscale1   	
    \skewchar\xivmi ='177  		
 \font\xivsy  = cmsy10	\@magscale2   	
    \skewchar\xivsy ='60 		
 \font\xivex  = cmex10		\@magscale2   	
 \font\xivmib = cmmib10\@magscale2   	
 \font\xivbsy = cmbsy10\@magscale2  	
 \font\xivit  = cmti12	\@magscale1   	
 \font\xivsl  = cmsl12	\@magscale1   	
 \font\xivbf  = cmbx12	\@magscale1   	
 \font\xivtt  = cmtt12	\@magscale1   		
    \hyphenchar\xivtt = -1         		
 \font\xivsf  = cmss12	\@magscale1   		
 \font\xivsfb = cmssbx10    \@magscale2 	

\font\xiveufb = eufb10	\@magscale2		
\font\xiveufm = eufm10	\@magscale2		

\font\xivmsam = msam10	\@magscale2		
\font\xivmsbm = msbm10	\@magscale2		


%
%
%
 \font\xviirm  = cmr17		    		
 \font\xviimi  = cmmi12	\@magscale2   	
    \skewchar\xviimi ='177  		
 \font\xviisy  = cmsy10	\@magscale3   	
    \skewchar\xviisy ='60 		
 \font\xviiex  = cmex10		\@magscale3   	
 \font\xviimib = cmmib10\@magscale3   	
 \font\xviibsy = cmbsy10\@magscale3  	
 \font\xviiit  = cmti12	\@magscale2   	
 \font\xviisl  = cmsl12	\@magscale2   	
 \font\xviibf  = cmbx12	\@magscale2   	
 \font\xviitt  = cmtt12	\@magscale2   		
    \hyphenchar\xviitt = -1         		
 \font\xviisf  = cmss17				
 \font\xviisfb = cmssbx10    \@magscale3 	

\font\xviieufb = eufb10	\@magscale3		
\font\xviieufm = eufm10	\@magscale3		

\font\xviimsam = msam10	\@magscale3		
\font\xviimsbm = msbm10	\@magscale3		


%
%
%
 \font\xxrm  = cmr17		\@magscale1      		
 \font\xxmi  = cmmi12	\@magscale3   	
    \skewchar\xxmi ='177  		
 \font\xxsy  = cmsy10	\@magscale4   	
    \skewchar\xxsy ='60 		
 \font\xxmib = cmmib10\@magscale4   	
 \font\xxbsy = cmbsy10\@magscale4  	
 \font\xxit  = cmti12	\@magscale3   	
 \font\xxsl  = cmsl12	\@magscale3   	
 \font\xxbf  = cmbx12	\@magscale3   	
 \font\xxtt  = cmtt12	\@magscale3   		
    \hyphenchar\xxtt = -1         		
 \font\xxsf  = cmss17		\@magscale1  		
 \font\xxsfb = cmssbx10    \@magscale4 	
 \font\xxvrm  = cmr17		\@magscale2      		
 \font\xxvmi  = cmmi12	\@magscale4   	
    \skewchar\xxvmi ='177  		
 \font\xxvsy  = cmsy10	\@magscale5   	
    \skewchar\xxvsy ='60 		
 \font\xxvit  = cmti12	\@magscale4   	
 \font\xxvsl  = cmsl12	\@magscale4   	
 \font\xxvbf  = cmbx12 \@magscale4   	
 \font\xxvtt  = cmtt12	\@magscale4   		
    \hyphenchar\xxvtt = -1         		
 \font\xxvsf  = cmss17		\@magscale2  		
\newskip\ttglue

\def\petit{\def\rm{\fam\z@\viiirm}%
  \textfont0=\viiirm \scriptfont0=\virm \scriptscriptfont0=\vrm 
  \textfont1=\viiimi \scriptfont1=\vimi \scriptscriptfont1=\vmi 
	\def\mit{\fam\@ne}
  \textfont2=\viiisy \scriptfont2=\visy \scriptscriptfont2=\vsy 
	\def\cal{\fam\tw@}
  \textfont3=\viiiex \scriptfont3=\viiex \scriptscriptfont3=\viiex 
    \def\it{\fam\itfam\viiiit}%
    \textfont\itfam=\viiiit 
    \def\sl{\fam\slfam\viiisl}%
    \textfont\slfam=\viiisl 
    \def\bf{\fam\bffam\viiibf}%
    \textfont\bffam=\viiibf \scriptfont\bffam=\vibf 
       \scriptscriptfont\bffam=\vbf 
    \def\tt{\fam\ttfam\viiitt}%
    \textfont\ttfam=\viiitt 
    \def\sf{\fam\sffam\viiisf}%
    \textfont\sffam=\viiisf \scriptfont\sffam=\visf 
       \scriptscriptfont\sffam=\vsf 
    \def\mib{\fam\mibfam\viiimib}%
    \textfont\mibfam=\viiimib \scriptfont\mibfam=\vimib 
       \scriptscriptfont\mibfam=\vmib 
    \def\bsy{\fam\bsyfam\viiibsy}%
    \textfont\bsyfam=\viiibsy \scriptfont\bsyfam=\vibsy 
       \scriptscriptfont\bsyfam=\vbsy 
    \def\fm{\fam\fmfam\viiieufm}%
    \textfont\fmfam=\viiieufm \scriptfont\fmfam=\vieufm 
       \scriptscriptfont\fmfam=\veufm 
    \def\fb{\fam\fbfam\viiieufb}%
    \textfont\fbfam=\viiieufb \scriptfont\fbfam=\vieufb 
       \scriptscriptfont\fbfam=\veufb 
    \def\sam{\fam\samfam\viiimsam}%
    \textfont\samfam=\viiimsam \scriptfont\samfam=\vimsam 
       \scriptscriptfont\samfam=\vmsam 
    \def\sbm{\fam\sbmfam\viiimsbm}%
    \textfont\sbmfam=\viiimsbm \scriptfont\sbmfam=\vimsbm 
       \scriptscriptfont\sbmfam=\vmsbm 
    \def\sfb{\fam\sfbfam\viiisfb}%
    \textfont\sfbfam=\viiisfb \scriptfont\sfbfam=\viisfb 
       \scriptscriptfont\sfbfam=\viisfb 
    \tt \ttglue=.5em plus.25em minus.15em 
  \normalbaselineskip=9pt 
  \setbox\strutbox=\hbox{\vrule height8.5pt depth3.5pt width\z@}%
  \normalbaselines\rm}

%
%
%
\def\borgis{\def\rm{\fam\z@\ixrm}%
  \textfont0=\ixrm \scriptfont0=\virm \scriptscriptfont0=\vrm 
  \textfont1=\ixmi \scriptfont1=\vimi \scriptscriptfont1=\vmi 
	\def\mit{\fam\@ne}
  \textfont2=\ixsy \scriptfont2=\visy \scriptscriptfont2=\vsy 
	\def\cal{\fam\tw@}
  \textfont3=\ixex \scriptfont3=\viiex \scriptscriptfont3=\viiex 
    \def\it{\fam\itfam\ixit}%
    \textfont\itfam=\ixit 
    \def\sl{\fam\slfam\ixsl}%
    \textfont\slfam=\ixsl 
    \def\bf{\fam\bffam\ixbf}%
    \textfont\bffam=\ixbf \scriptfont\bffam=\vibf 
       \scriptscriptfont\bffam=\vbf 
    \def\tt{\fam\ttfam\ixtt}%
    \textfont\ttfam=\ixtt 
    \def\sf{\fam\sffam\ixsf}%
    \textfont\sffam=\ixsf \scriptfont\sffam=\visf 
       \scriptscriptfont\sffam=\vsf 
    \def\mib{\fam\mibfam\ixmib}%
    \textfont\mibfam=\ixmib \scriptfont\mibfam=\vimib 
       \scriptscriptfont\mibfam=\vmib 
    \def\bsy{\fam\bsyfam\ixbsy}%
    \textfont\bsyfam=\ixbsy \scriptfont\bsyfam=\vibsy 
       \scriptscriptfont\bsyfam=\vbsy 
    \def\fm{\fam\fmfam\ixeufm}%
    \textfont\fmfam=\ixeufm \scriptfont\fmfam=\vieufm 
       \scriptscriptfont\fmfam=\veufm 
    \def\fb{\fam\fbfam\ixeufb}%
    \textfont\fbfam=\ixeufb \scriptfont\fbfam=\vieufb 
       \scriptscriptfont\fbfam=\veufb 
    \def\sam{\fam\samfam\ixmsam}%
    \textfont\samfam=\ixmsam \scriptfont\samfam=\vimsam 
       \scriptscriptfont\samfam=\vmsam 
    \def\sbm{\fam\sbmfam\ixmsbm}%
    \textfont\sbmfam=\ixmsbm \scriptfont\sbmfam=\vimsbm 
       \scriptscriptfont\sbmfam=\vmsbm 
    \def\sfb{\fam\sfbfam\ixsfb}%
    \textfont\sfbfam=\ixsfb \scriptfont\sfbfam=\viisfb 
       \scriptscriptfont\sfbfam=\viisfb 
    \tt \ttglue=.5em plus.25em minus.15em 
  \normalbaselineskip=11pt 
  \setbox\strutbox=\hbox{\vrule height8.5pt depth3.5pt width\z@}%
  \normalbaselines\rm}

%
%
%
\def\garmond{\def\rm{\fam\z@\xrm}%
  \textfont0=\xrm \scriptfont0=\viirm \scriptscriptfont0=\vrm 
  \textfont1=\xmi \scriptfont1=\viimi \scriptscriptfont1=\vmi 
	\def\mit{\fam\@ne}
  \textfont2=\xsy \scriptfont2=\viisy \scriptscriptfont2=\vsy 
	\def\cal{\fam\tw@}
  \textfont3=\xex \scriptfont3=\viiex \scriptscriptfont3=\viiex 
    \def\it{\fam\itfam\xit}%
    \textfont\itfam=\xit 
    \def\sl{\fam\slfam\xsl}%
    \textfont\slfam=\xsl 
    \def\bf{\fam\bffam\xbf}%
    \textfont\bffam=\xbf \scriptfont\bffam=\viibf 
       \scriptscriptfont\bffam=\vbf 
    \def\tt{\fam\ttfam\xtt}%
    \textfont\ttfam=\xtt 
    \def\sf{\fam\sffam\xsf}%
    \textfont\sffam=\xsf \scriptfont\sffam=\viisf 
       \scriptscriptfont\sffam=\vsf 
    \def\mib{\fam\mibfam\xmib}%
    \textfont\mibfam=\xmib \scriptfont\mibfam=\viimib 
       \scriptscriptfont\mibfam=\vmib 
    \def\bsy{\fam\bsyfam\xbsy}%
    \textfont\bsyfam=\xbsy \scriptfont\bsyfam=\viibsy 
       \scriptscriptfont\bsyfam=\vbsy 
    \def\fm{\fam\fmfam\xeufm}%
    \textfont\fmfam=\xeufm \scriptfont\fmfam=\viieufm 
       \scriptscriptfont\fmfam=\veufm 
    \def\fb{\fam\fbfam\xeufb}%
    \textfont\fbfam=\xeufb \scriptfont\fbfam=\viieufb 
       \scriptscriptfont\fbfam=\veufb 
    \def\sam{\fam\samfam\xmsam}%
    \textfont\samfam=\xmsam \scriptfont\samfam=\viimsam 
       \scriptscriptfont\samfam=\vmsam 
    \def\sbm{\fam\sbmfam\xmsbm}%
    \textfont\sbmfam=\xmsbm \scriptfont\sbmfam=\viimsbm 
       \scriptscriptfont\sbmfam=\vmsbm 
    \def\sfb{\fam\sfbfam\xsfb}%
    \textfont\sfbfam=\xsfb \scriptfont\sfbfam=\viisfb 
       \scriptscriptfont\sfbfam=\viisfb 
    \tt \ttglue=.5em plus.25em minus.15em 
 \normalbaselineskip=12pt 
  \setbox\strutbox=\hbox{\vrule height8.5pt depth3.5pt width\z@}%
  \normalbaselines\rm}

%
%
%
\def\hypocicero{\def\rm{\fam\z@\xirm}%
  \textfont0=\xirm \scriptfont0=\viiirm \scriptscriptfont0=\virm 
  \textfont1=\ximi \scriptfont1=\viiimi \scriptscriptfont1=\vimi 
	\def\mit{\fam\@ne}
  \textfont2=\xisy \scriptfont2=\viiisy \scriptscriptfont2=\visy 
	\def\cal{\fam\tw@}
  \textfont3=\xiex \scriptfont3=\viiiex \scriptscriptfont3=\viiex 
    \def\it{\fam\itfam\xiit}%
    \textfont\itfam=\xiit 
    \def\sl{\fam\slfam\xisl}%
    \textfont\slfam=\xisl 
    \def\bf{\fam\bffam\xibf}%
    \textfont\bffam=\xibf \scriptfont\bffam=\viiibf 
       \scriptscriptfont\bffam=\vibf 
    \def\tt{\fam\ttfam\xitt}%
    \textfont\ttfam=\xitt 
    \def\sf{\fam\sffam\xisf}%
    \textfont\sffam=\xisf \scriptfont\sffam=\viiisf 
       \scriptscriptfont\sffam=\visf 
    \def\mib{\fam\mibfam\ximib}%
    \textfont\mibfam=\ximib \scriptfont\mibfam=\viiimib 
       \scriptscriptfont\mibfam=\vimib 
    \def\bsy{\fam\bsyfam\xibsy}%
    \textfont\bsyfam=\xibsy \scriptfont\bsyfam=\viiibsy 
       \scriptscriptfont\bsyfam=\vibsy 
    \def\fm{\fam\fmfam\xieufm}%
    \textfont\fmfam=\xieufm \scriptfont\fmfam=\viiieufm 
       \scriptscriptfont\fmfam=\vieufm 
    \def\fb{\fam\fbfam\xieufb}%
    \textfont\fbfam=\xieufb \scriptfont\fbfam=\viiieufb 
       \scriptscriptfont\fbfam=\vieufb 
    \def\sam{\fam\samfam\ximsam}%
    \textfont\samfam=\ximsam \scriptfont\samfam=\viiimsam 
       \scriptscriptfont\samfam=\vimsam 
    \def\sbm{\fam\sbmfam\ximsbm}%
    \textfont\sbmfam=\ximsbm \scriptfont\sbmfam=\viiimsbm 
       \scriptscriptfont\sbmfam=\vimsbm 
    \def\sfb{\fam\sfbfam\xisfb}%
    \textfont\sfbfam=\xisfb \scriptfont\sfbfam=\viiisfb 
       \scriptscriptfont\sfbfam=\viisfb 
    \tt \ttglue=.5em plus.25em minus.15em 
   \normalbaselineskip=13,2pt 
   \setbox\strutbox=\hbox{\vrule height8.5pt depth3.5pt width\z@}%
   \normalbaselines\rm}

%
%
%
\def\cicero{\def\rm{\fam\z@\xiirm}%
  \textfont0=\xiirm \scriptfont0=\ixrm \scriptscriptfont0=\viirm 
  \textfont1=\xiimi \scriptfont1=\ixmi \scriptscriptfont1=\viimi 
	\def\mit{\fam\@ne}
  \textfont2=\xiisy \scriptfont2=\ixsy \scriptscriptfont2=\viisy 
	\def\cal{\fam\tw@}
  \textfont3=\xiiex \scriptfont3=\ixex \scriptscriptfont3=\viiex 
    \def\it{\fam\itfam\xiiit}%
    \textfont\itfam=\xiiit 
    \def\sl{\fam\slfam\xiisl}%
    \textfont\slfam=\xiisl 
    \def\bf{\fam\bffam\xiibf}%
    \textfont\bffam=\xiibf \scriptfont\bffam=\ixbf 
       \scriptscriptfont\bffam=\viibf 
    \def\tt{\fam\ttfam\xiitt}%
    \textfont\ttfam=\xiitt 
    \def\sf{\fam\sffam\xiisf}%
    \textfont\sffam=\xiisf \scriptfont\sffam=\ixsf 
       \scriptscriptfont\sffam=\viisf 
    \def\mib{\fam\mibfam\xiimib}%
    \textfont\mibfam=\xiimib \scriptfont\mibfam=\ixmib 
       \scriptscriptfont\mibfam=\viimib 
    \def\bsy{\fam\bsyfam\xiibsy}%
    \textfont\bsyfam=\xiibsy \scriptfont\bsyfam=\ixbsy 
       \scriptscriptfont\bsyfam=\viibsy 
    \def\fm{\fam\fmfam\xiieufm}%
    \textfont\fmfam=\xiieufm \scriptfont\fmfam=\ixeufm 
       \scriptscriptfont\fmfam=\viieufm 
    \def\fb{\fam\fbfam\xiieufb}%
    \textfont\fbfam=\xiieufb \scriptfont\fbfam=\ixeufb 
       \scriptscriptfont\fbfam=\viieufb 
    \def\sam{\fam\samfam\xiimsam}%
    \textfont\samfam=\xiimsam \scriptfont\samfam=\ixmsam 
       \scriptscriptfont\samfam=\viimsam 
    \def\sbm{\fam\sbmfam\xiimsbm}%
    \textfont\sbmfam=\xiimsbm \scriptfont\sbmfam=\ixmsbm 
       \scriptscriptfont\sbmfam=\viimsbm 
    \def\sfb{\fam\sfbfam\xiisfb}%
    \textfont\sfbfam=\xiisfb \scriptfont\sfbfam=\ixsfb 
       \scriptscriptfont\sfbfam=\viisfb 
    \tt \ttglue=.5em plus.25em minus.15em 
  \normalbaselineskip=14,3pt 
  \setbox\strutbox=\hbox{\vrule height8.5pt depth3.5pt width\z@}%
  \normalbaselines\rm}

%
%
%
\def\mittel{\def\rm{\fam\z@\xivrm}%
  \textfont0=\xivrm \scriptfont0=\xirm \scriptscriptfont0=\viiirm 
  \textfont1=\xivmi \scriptfont1=\ximi \scriptscriptfont1=\viiimi 
	\def\mit{\fam\@ne}
  \textfont2=\xivsy \scriptfont2=\xisy \scriptscriptfont2=\viiisy 
	\def\cal{\fam\tw@}
  \textfont3=\xivex \scriptfont3=\xiex \scriptscriptfont3=\viiiex 
    \def\it{\fam\itfam\xivit}%
    \textfont\itfam=\xivit 
    \def\sl{\fam\slfam\xivsl}%
    \textfont\slfam=\xivsl 
    \def\bf{\fam\bffam\xivbf}%
    \textfont\bffam=\xivbf \scriptfont\bffam=\xibf 
       \scriptscriptfont\bffam=\viiibf 
    \def\tt{\fam\ttfam\xivtt}%
    \textfont\ttfam=\xivtt 
    \def\sf{\fam\sffam\xivsf}%
    \textfont\sffam=\xivsf \scriptfont\sffam=\xisf 
       \scriptscriptfont\sffam=\viiisf 
    \def\mib{\fam\mibfam\xivmib}%
    \textfont\mibfam=\xivmib \scriptfont\mibfam=\ximib 
       \scriptscriptfont\mibfam=\viiimib 
    \def\bsy{\fam\bsyfam\xivbsy}%
    \textfont\bsyfam=\xivbsy \scriptfont\bsyfam=\xibsy 
       \scriptscriptfont\bsyfam=\viiibsy 
    \def\fm{\fam\fmfam\xiveufm}%
    \textfont\fmfam=\xiveufm \scriptfont\fmfam=\xieufm 
       \scriptscriptfont\fmfam=\viiieufm 
    \def\fb{\fam\fbfam\xiveufb}%
    \textfont\fbfam=\xiveufb \scriptfont\fbfam=\xieufb 
       \scriptscriptfont\fbfam=\viiieufb 
    \def\sam{\fam\samfam\xivmsam}%
    \textfont\samfam=\xivmsam \scriptfont\samfam=\ximsam 
       \scriptscriptfont\samfam=\viiimsam 
    \def\sbm{\fam\sbmfam\xivmsbm}%
    \textfont\sbmfam=\xivmsbm \scriptfont\sbmfam=\ximsbm 
       \scriptscriptfont\sbmfam=\viiimsbm 
    \def\sfb{\fam\sfbfam\xivsfb}%
    \textfont\sfbfam=\xivsfb \scriptfont\sfbfam=\xisfb 
       \scriptscriptfont\sfbfam=\viiisfb 
    \tt \ttglue=.5em plus.25em minus.15em 
   \normalbaselineskip=17pt 
   \setbox\strutbox=\hbox{\vrule height8.5pt depth3.5pt width\z@}%
   \normalbaselines\rm}

%
%
%
\def\hyppertertia{\def\rm{\fam\z@\xviirm}%
  \textfont0=\xviirm \scriptfont0=\xivrm \scriptscriptfont0=\xiirm 
  \textfont1=\xviimi \scriptfont1=\xivmi \scriptscriptfont1=\xiimi 
	\def\mit{\fam\@ne}
  \textfont2=\xviisy \scriptfont2=\xivsy \scriptscriptfont2=\xiisy 
	\def\cal{\fam\tw@}
  \textfont3=\xviiex \scriptfont3=\xivex \scriptscriptfont3=\xiiex 
    \def\it{\fam\itfam\xviiit}%
    \textfont\itfam=\xviiit 
    \def\sl{\fam\slfam\xviisl}%
    \textfont\slfam=\xviisl 
    \def\bf{\fam\bffam\xviibf}%
    \textfont\bffam=\xviibf \scriptfont\bffam=\xivbf 
       \scriptscriptfont\bffam=\xiibf 
    \def\tt{\fam\ttfam\xviitt}%
    \textfont\ttfam=\xviitt 
    \def\sf{\fam\sffam\xviisf}%
    \textfont\sffam=\xviisf \scriptfont\sffam=\xivsf 
       \scriptscriptfont\sffam=\xiisf 
    \def\mib{\fam\mibfam\xviimib}%
    \textfont\mibfam=\xviimib \scriptfont\mibfam=\xivmib 
       \scriptscriptfont\mibfam=\xiimib 
    \def\bsy{\fam\bsyfam\xviibsy}%
    \textfont\bsyfam=\xviibsy \scriptfont\bsyfam=\xivbsy 
       \scriptscriptfont\bsyfam=\xiibsy 
    \def\fm{\fam\fmfam\xviieufm}%
    \textfont\fmfam=\xviieufm \scriptfont\fmfam=\xiveufm 
       \scriptscriptfont\fmfam=\xiieufm 
    \def\fb{\fam\fbfam\xviieufb}%
    \textfont\fbfam=\xviieufb \scriptfont\fbfam=\xiveufb 
       \scriptscriptfont\fbfam=\xiieufb 
    \def\sam{\fam\samfam\xviimsam}%
    \textfont\samfam=\xviimsam \scriptfont\samfam=\xivmsam 
       \scriptscriptfont\samfam=\xiimsam 
    \def\sbm{\fam\sbmfam\xviimsbm}%
    \textfont\sbmfam=\xviimsbm \scriptfont\sbmfam=\xivmsbm 
       \scriptscriptfont\sbmfam=\xiimsbm 
    \def\sfb{\fam\sfbfam\xviisfb}%
    \textfont\sfbfam=\xviisfb \scriptfont\sfbfam=\xivsfb 
       \scriptscriptfont\sfbfam=\xiisfb 
    \tt \ttglue=.5em plus.25em minus.15em 
  \normalbaselineskip=21pt 
  \setbox\strutbox=\hbox{\vrule height8.5pt depth3.5pt width\z@}%
  \normalbaselines\rm}

%
%
%
\def\canon{\def\rm{\fam\z@\xxrm}%
  \textfont0=\xxrm \scriptfont0=\xviirm \scriptscriptfont0=\xivrm 
  \textfont1=\xxmi \scriptfont1=\xviimi \scriptscriptfont1=\xivmi 
	\def\mit{\fam\@ne}
  \textfont2=\xxsy \scriptfont2=\xviisy \scriptscriptfont2=\xivsy 
	\def\cal{\fam\tw@}
  \textfont3=\xivex \scriptfont3=\xviiex \scriptscriptfont3=\xivex 
    \def\it{\fam\itfam\xxit}%
    \textfont\itfam=\xxit 
    \def\sl{\fam\slfam\xxsl}%
    \textfont\slfam=\xxsl 
    \def\bf{\fam\bffam\xxbf}%
    \textfont\bffam=\xxbf \scriptfont\bffam=\xviibf 
       \scriptscriptfont\bffam=\xivbf 
    \def\tt{\fam\ttfam\xxtt}%
    \textfont\ttfam=\xxtt 
    \def\sf{\fam\sffam\xxsf}%
    \textfont\sffam=\xxsf \scriptfont\sffam=\xviisf 
       \scriptscriptfont\sffam=\xivsf 
    \def\mib{\fam\mibfam\xxmib}%
    \textfont\mibfam=\xxmib \scriptfont\mibfam=\xviimib 
       \scriptscriptfont\mibfam=\xivmib 
    \def\bsy{\fam\bsyfam\xxbsy}%
    \textfont\bsyfam=\xxbsy \scriptfont\bsyfam=\xviibsy 
       \scriptscriptfont\bsyfam=\xivbsy 
    \def\sfb{\fam\sfbfam\xxsfb}%
    \textfont\sfbfam=\xxsfb \scriptfont\sfbfam=\xviisfb 
       \scriptscriptfont\sfbfam=\xivsfb 
    \tt \ttglue=.5em plus.25em minus.15em 
  \normalbaselineskip=27pt 
  \setbox\strutbox=\hbox{\vrule height8.5pt depth3.5pt width\z@}%
  \normalbaselines\rm}

%
%

\def\biscicero{\def\rm{\fam\z@\xxvrm}%
  \textfont0=\xxvrm \scriptfont0=\xxrm \scriptscriptfont0=\xviirm 
  \textfont1=\xxvmi \scriptfont1=\xxmi \scriptscriptfont1=\xviimi 
	\def\mit{\fam\@ne}
  \textfont2=\xxvsy \scriptfont2=\xxsy \scriptscriptfont2=\xviisy 
	\def\cal{\fam\tw@}
%
    \def\it{\fam\itfam\xxvit}%
    \textfont\itfam=\xxvit 
    \def\sl{\fam\slfam\xxvsl}%
    \textfont\slfam=\xxvsl 
    \def\bf{\fam\bffam\xxvbf}%
    \textfont\bffam=\xxvbf \scriptfont\bffam=\xxbf 
       \scriptscriptfont\bffam=\xviibf 
    \def\tt{\fam\ttfam\xxvtt}%
    \textfont\ttfam=\xxvtt 
    \def\sf{\fam\sffam\xxvsf}%
    \textfont\sffam=\xxvsf \scriptfont\sffam=\xxsf 
       \scriptscriptfont\sffam=\xviisf 
    \tt \ttglue=.5em plus.25em minus.15em 
  \normalbaselineskip=32pt 
  \setbox\strutbox=\hbox{\vrule height8.5pt depth3.5pt width\z@}%
  \normalbaselines\rm}

\newfam\itfam
\newfam\slfam 
\newfam\bffam
\newfam\ttfam
\newfam\sffam
\newfam\mibfam
\newfam\bsyfam 
\newfam\fmfam
\newfam\fbfam
\newfam\samfam
\newfam\sbmfam
\newfam\sfbfam

%
%

\message{macros,}

\def\frenchspacing{\sfcode`\.\@m \sfcode`\?\@m \sfcode`\!\@m
  \sfcode`\:\@m \sfcode`\;\@m \sfcode`\,\@m}
\def\nonfrenchspacing{\sfcode`\.3000\sfcode`\?3000\sfcode`\!3000%
  \sfcode`\:2000\sfcode`\;1500\sfcode`\,1250 }

\def\normalbaselines{\lineskip\normallineskip
  \baselineskip\normalbaselineskip \lineskiplimit\normallineskiplimit}

\def\^^M{\ } 
\def\^^I{\ } 

 \let\endline=\cr

\def\space{ }
\def\empty{}
\def\null{\hbox{}}

\let\bgroup={ \let\egroup=}

{\catcode`\^^M=\active 
  \gdef\obeylines{\catcode`\^^M\active \let^^M\par}%
  \global\let^^M\par} 
\def\obeyspaces{\catcode`\ \active}
{\obeyspaces\global\let =\space}

\def\loop#1\repeat{\def\body{#1}\iterate}
\def\iterate{\body \let\next\iterate \else\let\next\relax\fi \next}
\let\repeat=\fi 

\def\thinspace{\kern .16667em }
\def\negthinspace{\kern-.16667em }
\def\enspace{\kern.5em }

\def\quad{\hskip1em\relax}
\def\qquad{\hskip2em\relax}

\def\smallskip{\vskip\smallskipamount}
\def\medskip{\vskip\medskipamount}
\def\bigskip{\vskip\bigskipamount}

\def\nointerlineskip{\prevdepth-1000\p@}
\def\offinterlineskip{\baselineskip-1000\p@
  \lineskip\z@ \lineskiplimit\maxdimen}

\def\vglue{\afterassignment\vgl@\skip@=}
\def\vgl@{\par \dimen@\prevdepth \hrule height\z@
  \nobreak\vskip\skip@ \prevdepth\dimen@}
\def\hglue{\afterassignment\hgl@\skip@=}
\def\hgl@{\leavevmode \count@\spacefactor \vrule width\z@
  \nobreak\hskip\skip@ \spacefactor\count@}

\def~{\penalty\@M \ } 

\def\break{\penalty-\@M}
\def\nobreak{\penalty \@M}
\def\allowbreak{\penalty \z@}

\def\eject{\par\break}
\def\supereject{\par\penalty-\@MM}

\def\removelastskip{\ifdim\lastskip=\z@\else\vskip-\lastskip\fi}
\def\smallbreak{\par\ifdim\lastskip<\smallskipamount
  \removelastskip\penalty-50\smallskip\fi}
\def\medbreak{\par\ifdim\lastskip<\medskipamount
  \removelastskip\penalty-100\medskip\fi}
\def\bigbreak{\par\ifdim\lastskip<\bigskipamount
  \removelastskip\penalty-200\bigskip\fi}

\def\line{\hbox to\hsize}
\def\leftline#1{\line{#1\hss}}
\def\rightline#1{\line{\hss#1}}
\def\centerline#1{\line{\hss#1\hss}}

\def\rlap#1{\hbox to\z@{#1\hss}}
\def\llap#1{\hbox to\z@{\hss#1}}

\def\m@th{\mathsurround=\z@}
\def\underbar#1{$\setbox\z@\hbox{#1}\dp\z@\z@
  \m@th \underline{\box\z@}$}

\newbox\strutbox
\setbox\strutbox=\hbox{\vrule height8.5pt depth3.5pt width\z@}
\def\strut{\relax\ifmmode\copy\strutbox\else\unhcopy\strutbox\fi}

\def\hidewidth{\hskip\hideskip} 
\def\ialign{\everycr{}\tabskip\z@skip\halign} 
\newcount\mscount
\def\multispan#1{\omit \mscount#1
  \loop\ifnum\mscount>\@ne \sp@n\repeat}
\def\sp@n{\span\omit\advance\mscount\m@ne}

\newif\ifus@ \newif\if@cr
\newbox\tabs \newbox\tabsyet \newbox\tabsdone

\def\cleartabs{\global\setbox\tabsyet\null \setbox\tabs\null}
\def\settabs{\setbox\tabs\null \futurelet\next\sett@b}
\let\+=\relax 
\def\sett@b{\ifx\next\+\let\next\relax
    \def\next{\afterassignment\s@tt@b\let\next}%
  \else\let\next\s@tcols\fi\next}
\def\s@tt@b{\let\next\relax\us@false\m@ketabbox}
\def\tabalign{\us@true\m@ketabbox} 
\outer\def\+{\tabalign}
\def\s@tcols#1\columns{\count@#1 \dimen@\hsize
  \loop\ifnum\count@>\z@ \@nother \repeat}
\def\@nother{\dimen@ii\dimen@ \divide\dimen@ii\count@
  \setbox\tabs\hbox{\hbox to\dimen@ii{}\unhbox\tabs}%
  \advance\dimen@-\dimen@ii \advance\count@\m@ne}

\def\m@ketabbox{\begingroup
  \global\setbox\tabsyet\copy\tabs
  \global\setbox\tabsdone\null
  \def\cr{\@crtrue\crcr\egroup\egroup
    \ifus@\unvbox\z@\lastbox\fi\endgroup
    \setbox\tabs\hbox{\unhbox\tabsyet\unhbox\tabsdone}}%
  \setbox\z@\vbox\bgroup\@crfalse
    \ialign\bgroup&\t@bbox##\t@bb@x\crcr}

\def\t@bbox{\setbox\z@\hbox\bgroup}
\def\t@bb@x{\if@cr\egroup 
  \else\hss\egroup \global\setbox\tabsyet\hbox{\unhbox\tabsyet
      \global\setbox\@ne\lastbox}
    \ifvoid\@ne\global\setbox\@ne\hbox to\wd\z@{}%
    \else\setbox\z@\hbox to\wd\@ne{\unhbox\z@}\fi
    \global\setbox\tabsdone\hbox{\box\@ne\unhbox\tabsdone}\fi
  \box\z@}

\def\hang{\hangindent\parindent}
\def\textindent#1{\indent\llap{#1\enspace}\ignorespaces}
\def\item{\par\hang\textindent}

\def\narrower{\advance\leftskip\parindent
  \advance\rightskip\parindent}

\outer\def\beginsection#1\par{\vskip\z@ plus.3\vsize\penalty-250
  \vskip\z@ plus-.3\vsize\bigskip\vskip\parskip
  \message{#1}\leftline{\bf#1}\nobreak\smallskip\noindent}
\outer\def\proclaim #1. #2\par{\medbreak
  \noindent{\bf#1.\enspace}{\sl#2}\par
  \ifdim\lastskip<\medskipamount \removelastskip\penalty55\medskip\fi}

\def\raggedright{\rightskip\z@ plus2em \spaceskip.3333em \xspaceskip.5em\relax}
\def\ttraggedright{\tt\rightskip\z@ plus2em\relax} 

\chardef\%=`\%
\chardef\&=`\&
\chardef\#=`\#
\chardef\$=`\$
\chardef\ss="19
\chardef\ae="1A
\chardef\oe="1B
\chardef\o="1C
\chardef\AE="1D
\chardef\OE="1E
\chardef\O="1F
\chardef\i="10 \chardef\j="11 

\def\L{\leavevmode\setbox0\hbox{L}\hbox to\wd0{\hss\char32L}}

\def\leavevmode{\unhbox\voidb@x} 
\def\_{\leavevmode \kern.06em \vbox{\hrule width.3em}}
\def\AA{\leavevmode\setbox0\hbox{h}\dimen@\ht0\advance\dimen@-1ex%
  \rlap{\raise.67\dimen@\hbox{\char'27}}A}

\def\mathhexbox#1#2#3{\leavevmode
  \hbox{$\m@th \mathchar"#1#2#3$}}

\def\oalign#1{\leavevmode\vtop{\baselineskip\z@skip \lineskip.25ex%
  \ialign{##\crcr#1\crcr}}} 
\def\ooalign{\lineskiplimit-\maxdimen \oalign}

\def\c#1{\setbox\z@\hbox{#1}\ifdim\ht\z@=1ex\accent24 #1%
  \else{\ooalign{\hidewidth\char24\hidewidth\crcr\unhbox\z@}}\fi}
\def\copyright{{\ooalign{\hfil\raise.07ex\hbox{c}\hfil\crcr\mathhexbox20D}}}

\def\dots{\relax\ifmmode\ldots\else$\m@th\ldots\,$\fi}
\def\TeX{T\kern-.1667em\lower.5ex\hbox{E}\kern-.125emX}

\def\`#1{{\accent18 #1}}
\def\'#1{{\accent19 #1}}
\def\v#1{{\accent20 #1}} \let\^^_=\v
\def\u#1{{\accent21 #1}} \let\^^S=\u
\def\=#1{{\accent22 #1}}
\def\^#1{{\accent94 #1}} \let\^^D=\^
\def\.#1{{\accent95 #1}}
\def\H#1{{\accent"7D #1}}
\def\~#1{{\accent"7E #1}}
\def\"#1{{\accent"7F #1}}
\def\t#1{{\edef\next{\the\font}\the\textfont1\accent"7F\next#1}}

\def\dotfill{\cleaders\hbox{$\m@th \mkern1.5mu.\mkern1.5mu$}\hfill}
\def\rightarrowfill{$\m@th\mathord-\mkern-6mu%
  \cleaders\hbox{$\mkern-2mu\mathord-\mkern-2mu$}\hfill
  \mkern-6mu\mathord\rightarrow$}
\def\leftarrowfill{$\m@th\mathord\leftarrow\mkern-6mu%
  \cleaders\hbox{$\mkern-2mu\mathord-\mkern-2mu$}\hfill
  \mkern-6mu\mathord-$}
\mathchardef\braceld="37A \mathchardef\bracerd="37B
\mathchardef\bracelu="37C \mathchardef\braceru="37D
\def\downbracefill{$\m@th\braceld\leaders\vrule\hfill\braceru
  \bracelu\leaders\vrule\hfill\bracerd$}
\def\upbracefill{$\m@th\bracelu\leaders\vrule\hfill\bracerd
  \braceld\leaders\vrule\hfill\braceru$}

\outer\def\bye{\par\vfill\supereject\end}
\message{math definitions,}

\let\sp=^ \let\sb=_
\def\,{\mskip\thinmuskip}
\def\>{\mskip\medmuskip}
\def\;{\mskip\thickmuskip}
\def\!{\mskip-\thinmuskip}
\def\*{\discretionary{\thinspace\the\textfont2\char2}{}{}}
{\catcode`\'=\active \gdef'{^\bgroup\prim@s}}
\def\prim@s{\prime\futurelet\next\pr@m@s}
\def\pr@m@s{\ifx'\next\let\nxt\pr@@@s \else\ifx^\next\let\nxt\pr@@@t
  \else\let\nxt\egroup\fi\fi \nxt}
\def\pr@@@s#1{\prim@s} \def\pr@@@t#1#2{#2\egroup}
{\catcode`\^^Z=\active \gdef^^Z{\not=}} 

{\catcode`\_=\active \global\let_=\_} 

\mathchardef\alpha="710B
\mathchardef\beta="710C
\mathchardef\gamma="710D
\mathchardef\delta="710E
\mathchardef\epsilon="710F
\mathchardef\zeta="7110
\mathchardef\eta="7111
\mathchardef\theta="7112
\mathchardef\iota="7113
\mathchardef\kappa="7114
\mathchardef\lambda="7115
\mathchardef\mu="7116
\mathchardef\nu="7117
\mathchardef\xi="7118
\mathchardef\pi="7119
\mathchardef\rho="711A
\mathchardef\sigma="711B
\mathchardef\tau="711C
\mathchardef\upsilon="711D
\mathchardef\phi="711E
\mathchardef\chi="711F
\mathchardef\psi="7120
\mathchardef\omega="7121
\mathchardef\varepsilon="7122
\mathchardef\vartheta="7123
\mathchardef\varpi="7124
\mathchardef\varrho="7125
\mathchardef\varsigma="7126
\mathchardef\varphi="7127
\mathchardef\Gamma="7000
\mathchardef\Delta="7001
\mathchardef\Theta="7002
\mathchardef\Lambda="7003
\mathchardef\Xi="7004
\mathchardef\Pi="7005
\mathchardef\Sigma="7006
\mathchardef\Upsilon="7007
\mathchardef\Phi="7008
\mathchardef\Psi="7009
\mathchardef\Omega="700A

\mathchardef\aleph="0240
\def\hbar{{\mathchar'26\mkern-9muh}}
\mathchardef\imath="017B
\mathchardef\jmath="017C
\mathchardef\ell="0160
\mathchardef\wp="017D
\mathchardef\Re="023C
\mathchardef\Im="023D
\mathchardef\partial="0140
\mathchardef\infty="0231
\mathchardef\prime="0230
\mathchardef\emptyset="023B
\mathchardef\nabla="0272

\mathchardef\top="023E
\mathchardef\bot="023F
\def\angle{{\vbox{\ialign{$\m@th\scriptstyle##$\crcr
      \not\mathrel{\mkern14mu}\crcr
      \noalign{\nointerlineskip}
      \mkern2.5mu\leaders\hrule height.34pt\hfill\mkern2.5mu\crcr}}}}
\mathchardef\triangle="0234
\mathchardef\forall="0238
\mathchardef\exists="0239
\mathchardef\neg="023A 
\mathchardef\flat="015B
\mathchardef\natural="015C
\mathchardef\sharp="015D
\mathchardef\clubsuit="027C
\mathchardef\diamondsuit="027D
\mathchardef\heartsuit="027E
\mathchardef\spadesuit="027F

\mathchardef\coprod="1360
\mathchardef\bigvee="1357
\mathchardef\bigwedge="1356
\mathchardef\biguplus="1355
\mathchardef\bigcap="1354
\mathchardef\bigcup="1353
\mathchardef\intop="1352 \def\int{\intop\nolimits}
\mathchardef\prod="1351
\mathchardef\sum="1350
\mathchardef\bigotimes="134E
\mathchardef\bigoplus="134C
\mathchardef\bigodot="134A
\mathchardef\ointop="1348 
\mathchardef\bigsqcup="1346
\mathchardef\smallint="1273

\mathchardef\triangleleft="212F
\mathchardef\triangleright="212E
\mathchardef\bigtriangleup="2234
\mathchardef\bigtriangledown="2235
\mathchardef\wedge="225E 
\mathchardef\vee="225F 
\mathchardef\cap="225C
\mathchardef\cup="225B
\mathchardef\ddagger="227A
\mathchardef\dagger="2279
\mathchardef\sqcap="2275
\mathchardef\sqcup="2274
\mathchardef\uplus="225D
\mathchardef\amalg="2271
\mathchardef\diamond="2205
\mathchardef\bullet="220F
\mathchardef\wr="226F
\mathchardef\div="2204
\mathchardef\odot="220C
\mathchardef\oslash="220B
\mathchardef\otimes="220A
\mathchardef\ominus="2209
\mathchardef\oplus="2208
\mathchardef\mp="2207
\mathchardef\pm="2206
\mathchardef\circ="220E
\mathchardef\bigcirc="220D
\mathchardef\setminus="226E 
\mathchardef\cdot="2201
\mathchardef\ast="2203
\mathchardef\times="2202
\mathchardef\star="213F

\mathchardef\propto="322F
\mathchardef\sqsubseteq="3276
\mathchardef\sqsupseteq="3277
\mathchardef\parallel="326B
\mathchardef\mid="326A
\mathchardef\dashv="3261
\mathchardef\vdash="3260
\mathchardef\nearrow="3225
\mathchardef\searrow="3226
\mathchardef\nwarrow="322D
\mathchardef\swarrow="322E
\mathchardef\Leftrightarrow="322C
\mathchardef\Leftarrow="3228
\mathchardef\Rightarrow="3229
 \let\ne=\neq
\mathchardef\leq="3214 
\mathchardef\geq="3215 
\mathchardef\succ="321F
\mathchardef\prec="321E
\mathchardef\approx="3219
\mathchardef\succeq="3217
\mathchardef\preceq="3216
\mathchardef\supset="321B
\mathchardef\subset="321A
\mathchardef\supseteq="3213
\mathchardef\subseteq="3212
\mathchardef\in="3232
\mathchardef\ni="3233 
\mathchardef\gg="321D
\mathchardef\ll="321C
\mathchardef\not="3236
\mathchardef\leftrightarrow="3224
\mathchardef\leftarrow="3220 
\mathchardef\rightarrow="3221 \let\to=\rightarrow
\mathchardef\mapstochar="3237 \def\mapsto{\mapstochar\rightarrow}
\mathchardef\sim="3218
\mathchardef\simeq="3227
\mathchardef\perp="323F
\mathchardef\equiv="3211
\mathchardef\asymp="3210
\mathchardef\smile="315E
\mathchardef\frown="315F
\mathchardef\leftharpoonup="3128
\mathchardef\leftharpoondown="3129
\mathchardef\rightharpoonup="312A
\mathchardef\rightharpoondown="312B

\def\joinrel{\mathrel{\mkern-3mu}}
\def\relbar{\mathrel{\smash-}} 
\def\Relbar{\mathrel=}
\mathchardef\lhook="312C \def\hookrightarrow{\lhook\joinrel\rightarrow}
\mathchardef\rhook="312D

\def\Longrightarrow{\Relbar\joinrel\Rightarrow}
\def\longrightarrow{\relbar\joinrel\rightarrow}

\def\Longleftrightarrow{\Leftarrow\joinrel\Rightarrow}

\mathchardef\ldotp="613A 
\mathchardef\cdotp="6201 
\mathchardef\colon="603A 
\def\ldots{\mathinner{\ldotp\ldotp\ldotp}}

\def\vdots{\vbox{\baselineskip4\p@ \lineskiplimit\z@
    \kern6\p@\hbox{.}\hbox{.}\hbox{.}}}
\def\ddots{\mathinner{\mkern1mu\raise7\p@\vbox{\kern7\p@\hbox{.}}\mkern2mu
    \raise4\p@\hbox{.}\mkern2mu\raise\p@\hbox{.}\mkern1mu}}

\def\tilde{\mathaccent"707E }
\def\bar{\mathaccent"7016 }

\def\hat{\mathaccent"705E }

\def\widehat{\mathaccent"0362 }
\def\overrightarrow#1{\vbox{\ialign{##\crcr
      \rightarrowfill\crcr\noalign{\kern-\p@\nointerlineskip}
      $\hfil\displaystyle{#1}\hfil$\crcr}}}
\def\overleftarrow#1{\vbox{\ialign{##\crcr
      \leftarrowfill\crcr\noalign{\kern-\p@\nointerlineskip}
      $\hfil\displaystyle{#1}\hfil$\crcr}}}
\def\overbrace#1{\mathop{\vbox{\ialign{##\crcr\noalign{\kern3\p@}
      \downbracefill\crcr\noalign{\kern3\p@\nointerlineskip}
      $\hfil\displaystyle{#1}\hfil$\crcr}}}\limits}
\def\underbrace#1{\mathop{\vtop{\ialign{##\crcr
      $\hfil\displaystyle{#1}\hfil$\crcr\noalign{\kern3\p@\nointerlineskip}
      \upbracefill\crcr\noalign{\kern3\p@}}}}\limits}

\def\Vert{\delimiter"26B30D } \let\|=\Vert

\def\rangle{\delimiter"526930B }
\def\langle{\delimiter"426830A }
\def\rbrace{\delimiter"5267309 } \let\}=\rbrace
\def\lbrace{\delimiter"4266308 } \let\{=\lbrace

\def\bigl{\mathopen\big}

\def\bigr{\mathclose\big}
\def\Bigl{\mathopen\Big}

\def\Bigr{\mathclose\Big}
\def\biggl{\mathopen\bigg}

\def\biggr{\mathclose\bigg}
\def\Biggl{\mathopen\Bigg}

\def\Biggr{\mathclose\Bigg}
\def\big#1{{\hbox{$\left#1\vbox to8.5\p@{}\right.\n@space$}}}
\def\Big#1{{\hbox{$\left#1\vbox to11.5\p@{}\right.\n@space$}}}
\def\bigg#1{{\hbox{$\left#1\vbox to14.5\p@{}\right.\n@space$}}}
\def\Bigg#1{{\hbox{$\left#1\vbox to17.5\p@{}\right.\n@space$}}}
\def\n@space{\nulldelimiterspace\z@ \m@th}

\def\choose{\atopwithdelims()}

\def\sqrt{\radical"270370 }

\def\mathpalette#1#2{\mathchoice{#1\displaystyle{#2}}%
  {#1\textstyle{#2}}{#1\scriptstyle{#2}}{#1\scriptscriptstyle{#2}}}
\newbox\rootbox
\def\root#1\of{\setbox\rootbox\hbox{$\m@th\scriptscriptstyle{#1}$}
  \mathpalette\r@@t}
\def\r@@t#1#2{\setbox\z@\hbox{$\m@th#1\sqrt{#2}$}
  \dimen@\ht\z@ \advance\dimen@-\dp\z@
  \mkern5mu\raise.6\dimen@\copy\rootbox \mkern-10mu \box\z@}
\newif\ifv@ \newif\ifh@
\def\vphantom{\v@true\h@false\ph@nt}
\def\hphantom{\v@false\h@true\ph@nt}
\def\phantom{\v@true\h@true\ph@nt}
\def\ph@nt{\ifmmode\def\next{\mathpalette\mathph@nt}%
  \else\let\next\makeph@nt\fi\next}
\def\makeph@nt#1{\setbox\z@\hbox{#1}\finph@nt}
\def\mathph@nt#1#2{\setbox\z@\hbox{$\m@th#1{#2}$}\finph@nt}
\def\finph@nt{\setbox\tw@\null
  \ifv@ \ht\tw@\ht\z@ \dp\tw@\dp\z@\fi
  \ifh@ \wd\tw@\wd\z@\fi \box\tw@}
\def\mathstrut{\vphantom(}
\def\smash{\relax 
  \ifmmode\def\next{\mathpalette\mathsm@sh}\else\let\next\makesm@sh
  \fi\next}
\def\makesm@sh#1{\setbox\z@\hbox{#1}\finsm@sh}
\def\mathsm@sh#1#2{\setbox\z@\hbox{$\m@th#1{#2}$}\finsm@sh}
\def\finsm@sh{\ht\z@\z@ \dp\z@\z@ \box\z@}

\def\cong{\mathrel{\mathpalette\@vereq\sim}} 
\def\@vereq#1#2{\lower.5\p@\vbox{\baselineskip\z@skip\lineskip-.5\p@
    \ialign{$\m@th#1\hfil##\hfil$\crcr#2\crcr=\crcr}}}
\def\notin{\mathrel{\mathpalette\c@ncel\in}}
\def\c@ncel#1#2{\ooalign{$\hfil#1\mkern1mu/\hfil$\crcr$#1#2$}}
\def\rightleftharpoons{\mathrel{\mathpalette\rlh@{}}}
\def\rlh@#1{\vcenter{\hbox{\ooalign{\raise2pt
          \hbox{$#1\rightharpoonup$}\crcr
        $#1\leftharpoondown$}}}}
\def\buildrel#1\over#2{\mathrel{\mathop{\kern\z@#2}\limits^{#1}}}

\def\lim{\mathop{\rm lim}}
\def\limsup{\mathop{\rm lim\,sup}}

\def\max{\mathop{\rm max}}
\def\min{\mathop{\rm min}}
\def\sup{\mathop{\rm sup}}

\def\Pr{\mathop{\rm Pr}}

\def\cases#1{\left\{\,\vcenter{\normalbaselines\m@th
    \ialign{$##\hfil$&\quad##\hfil\crcr#1\crcr}}\right.}
\def\matrix#1{\null\,\vcenter{\normalbaselines\m@th
    \ialign{\hfil$##$\hfil&&\quad\hfil$##$\hfil\crcr
      \mathstrut\crcr\noalign{\kern-\baselineskip}
      #1\crcr\mathstrut\crcr\noalign{\kern-\baselineskip}}}\,}

\newdimen\p@renwd
\setbox0=\hbox{\xex B} \p@renwd=\wd0 
\def\bordermatrix#1{\begingroup \m@th
  \setbox\z@\vbox{\def\cr{\crcr\noalign{\kern2\p@\global\let\cr\endline}}%
    \ialign{$##$\hfil\kern2\p@\kern\p@renwd&\thinspace\hfil$##$\hfil
      &&\quad\hfil$##$\hfil\crcr
      \omit\strut\hfil\crcr\noalign{\kern-\baselineskip}%
      #1\crcr\omit\strut\cr}}%
  \setbox\tw@\vbox{\unvcopy\z@\global\setbox\@ne\lastbox}%
  \setbox\tw@\hbox{\unhbox\@ne\unskip\global\setbox\@ne\lastbox}%
  \setbox\tw@\hbox{$\kern\wd\@ne\kern-\p@renwd\left(\kern-\wd\@ne
    \global\setbox\@ne\vbox{\box\@ne\kern2\p@}%
    \vcenter{\kern-\ht\@ne\unvbox\z@\kern-\baselineskip}\,\right)$}%
  \null\;\vbox{\kern\ht\@ne\box\tw@}\endgroup}

\def\openup{\afterassignment\@penup\dimen@=}
\def\@penup{\advance\lineskip\dimen@
  \advance\baselineskip\dimen@
  \advance\lineskiplimit\dimen@}
\def\eqalign#1{\null\,\vcenter{\openup\jot\m@th
  \ialign{\strut\hfil$\displaystyle{##}$&$\displaystyle{{}##}$\hfil
      \crcr#1\crcr}}\,}
\newif\ifdt@p
\def\displ@y{\global\dt@ptrue\openup\jot\m@th
  \everycr{\noalign{\ifdt@p \global\dt@pfalse
      \vskip-\lineskiplimit \vskip\normallineskiplimit
      \else \penalty\interdisplaylinepenalty \fi}}}
\def\@lign{\tabskip\z@skip\everycr{}} 
\def\displaylines#1{\displ@y
  \halign{\hbox to\displaywidth{$\@lign\hfil\displaystyle##\hfil$}\crcr
    #1\crcr}}
\def\eqalignno#1{\displ@y \tabskip\centering
  \halign to\displaywidth{\hfil$\@lign\displaystyle{##}$\tabskip\z@skip
    &$\@lign\displaystyle{{}##}$\hfil\tabskip\centering
    &\llap{$\@lign##$}\tabskip\z@skip\crcr
    #1\crcr}}
\def\leqalignno#1{\displ@y \tabskip\centering
  \halign to\displaywidth{\hfil$\@lign\displaystyle{##}$\tabskip\z@skip
    &$\@lign\displaystyle{{}##}$\hfil\tabskip\centering
    &\kern-\displaywidth\rlap{$\@lign##$}\tabskip\displaywidth\crcr
    #1\crcr}}

\message{output routines,}

\countdef\pageno=0 \pageno=1 
\newtoks\headline \headline={\hfil} 
\newtoks\footline \footline={\hss\xrm\folio\hss}
\newif\ifr@ggedbottom
\def\raggedbottom{\topskip 10\p@ plus60\p@ \r@ggedbottomtrue}
\def\normalbottom{\topskip 10\p@ \r@ggedbottomfalse} 
\def\folio{\ifnum\pageno<\z@ \romannumeral-\pageno \else\number\pageno \fi}
\def\nopagenumbers{\footline{\hfil}} 
\def\advancepageno{\ifnum\pageno<\z@ \global\advance\pageno\m@ne
  \else\global\advance\pageno\@ne \fi} 

\newinsert\footins
\def\footnote#1{\let\@sf\empty 
  \ifhmode\edef\@sf{\spacefactor\the\spacefactor}\/\fi
  #1\@sf\vfootnote{#1}}
\def\vfootnote#1{\insert\footins\bgroup
  \interlinepenalty\interfootnotelinepenalty
  \splittopskip\ht\strutbox 
  \splitmaxdepth\dp\strutbox \floatingpenalty\@MM
  \leftskip\z@skip \rightskip\z@skip \spaceskip\z@skip \xspaceskip\z@skip
  \textindent{#1}\footstrut\futurelet\next\fo@t}
\def\fo@t{\ifcat\bgroup\noexpand\next \let\next\f@@t
  \else\let\next\f@t\fi \next}
\def\f@@t{\bgroup\aftergroup\@foot\let\next}
\def\f@t#1{#1\@foot}
\def\@foot{\strut\egroup}
\def\footstrut{\vbox to\splittopskip{}}
\skip\footins=\bigskipamount 
\count\footins=1000 
\dimen\footins=8in 

\newinsert\topins
\newif\ifp@ge \newif\if@mid
\def\topinsert{\@midfalse\p@gefalse\@ins}
\def\midinsert{\@midtrue\@ins}
\def\pageinsert{\@midfalse\p@getrue\@ins}
\skip\topins=\z@skip 
\count\topins=1000 
\dimen\topins=\maxdimen 
\def\@ins{\par\begingroup\setbox\z@\vbox\bgroup} 
\def\endinsert{\egroup 
  \if@mid \dimen@\ht\z@ \advance\dimen@\dp\z@
    \advance\dimen@12\p@ \advance\dimen@\pagetotal
    \ifdim\dimen@>\pagegoal\@midfalse\p@gefalse\fi\fi
  \if@mid \bigskip\box\z@\bigbreak
  \else\insert\topins{\penalty100 
    \splittopskip\z@skip
    \splitmaxdepth\maxdimen \floatingpenalty\z@
    \ifp@ge \dimen@\dp\z@
    \vbox to\vsize{\unvbox\z@\kern-\dimen@}
    \else \box\z@\nobreak\bigskip\fi}\fi\endgroup}

\output{\plainoutput}
\def\plainoutput{\shipout\vbox{\makeheadline\pagebody\makefootline}%
  \advancepageno
  \ifnum\outputpenalty>-\@MM \else\dosupereject\fi}
\def\pagebody{\vbox to\vsize{\boxmaxdepth\maxdepth \pagecontents}}
\def\makeheadline{\vbox to\z@{\vskip-22.5\p@
  \line{\vbox to8.5\p@{}\the\headline}\vss}\nointerlineskip}
\def\makefootline{\baselineskip24\p@\line{\the\footline}}
\def\dosupereject{\ifnum\insertpenalties>\z@ 
  \line{}\kern-\topskip\nobreak\vfill\supereject\fi}

\def\pagecontents{\ifvoid\topins\else\unvbox\topins\fi
  \dimen@=\dp\@cclv \unvbox\@cclv 
  \ifvoid\footins\else 
    \vskip\skip\footins
    \footnoterule
    \unvbox\footins\fi
  \ifr@ggedbottom \kern-\dimen@ \vfil \fi}
\def\footnoterule{\kern-3\p@
  \hrule width 2truein \kern 2.6\p@} 
\message{hyphenation}

\def\magnification{\afterassignment\m@g\count@}
\def\m@g{\mag\count@
  \hsize6.5truein\vsize8.9truein\dimen\footins8truein}

\def\tracingall{\tracingonline\@ne\tracingcommands\tw@\tracingstats\tw@
  \tracingpages\@ne\tracingoutput\@ne\tracinglostchars\@ne
  \tracingmacros\tw@\tracingparagraphs\@ne\tracingrestores\@ne
  \showboxbreadth\maxdimen\showboxdepth\maxdimen\errorstopmode}

\def\showhyphens#1{\setbox0\vbox{\parfillskip\z@skip\hsize\maxdimen\xrm
  \pretolerance\m@ne\tolerance\m@ne\hbadness0\showboxdepth0\ #1}}

\def\rm{\fam\z@\xrm}
\normalbaselines\rm 
\nonfrenchspacing 
\catcode`@=12 






\let\centreline=\centerline 

\font\sc=cmcsc11

 
\def\today{\ifcase\day\or 
 1$^{\rm st}$\or 2$^{\rm nd}$\or 3$^{\rm rd}$\or 4$^{\rm th}$\or 5$^{\rm th}$\or 
 6$^{\rm th}$\or 7$^{\rm th}$\or 8$^{\rm th}$\or 9$^{\rm th}$\or 10$^{\rm th}$\or 
 11$^{\rm th}$\or 12$^{\rm th}$\or 13$^{\rm th}$\or 14$^{\rm th}$\or 15$^{\rm th}$\or 
 16$^{\rm th}$\or 17$^{\rm th}$\or 18$^{\rm th}$\or 19$^{\rm th}$\or 20$^{\rm th}$\or 
 21$^{\rm st}$\or 22$^{\rm nd}$\or 23$^{\rm rd}$\or 24$^{\rm th}$\or 25$^{\rm th}$\or 
 26$^{\rm th}$\or 27$^{\rm th}$\or 28$^{\rm th}$\or 29$^{\rm th}$\or 30$^{\rm th}$\or 
 31$^{\rm st}$\fi 
 ~\ifcase\month\or 
 January\or February\or March\or April\or May\or June\or 
 July\or August\or September\or October\or November\or December\fi 
 \space \number\year}

\ifx\pisac\undefined \def\pisac{Nenad Antoni\'c}\fi
\ifx\predmet\undefined \def\predmet{No subject}\fi

\ifx\naslova\undefined \def\naslova{-}\fi
\ifx\naslovb\undefined \def\naslovb{}\fi

\ifx\naslovr\undefined \def\naslovr{\naslova}\fi

\nopagenumbers

\headline={\ifnum \pageno>0  { \xsf \draft }	\hfil
				{\xrm  \hss\naslovr } \else \hss\fi}

\footline={\ifnum \pageno>0  \xsf  \pisac \hfil{\xbf \hss \folio}
				 \fi}

\def\podnaslov#1{{\bf \par\bigskip
{#1}\nobreak
\smallskip\nobreak}
}

\ifx\dnum\undefined
   \def\odjeljak#1{\bigskip{\cicero \bf\par\bigskip\centreline{ #1}\nobreak
	\medskip\nobreak} }
   \def\for{\global\advance\bfor by1 \number\bfor}  
   \def\pbfor#1{\btmp=\bfor \advance\btmp by-#1 \number\btmp}	
\else
   \def\odjeljak#1{\bigskip{\cicero \bf\par\bigskip\centreline{ \global\advance\bodj by1 \number\bodj. #1}\nobreak \medskip\nobreak \global\bfor=0} }

   \def\for{\global\advance\bfor by1 \number\bodj.\number\bfor} 
   \def\pbfor#1{\btmp=\bfor \advance\btmp by-#1 \number\bodj.\number\btmp}	
\fi

\def\appendix#1{{\cicero \bf\par\bigskip\centreline{Appendix. #1}\nobreak \medskip\nobreak \global\bfor=0} 
}
%


\def\plainoutput{\shipout\vbox{\makeheadline\pagebody\makefootline}%
  \advancepageno
  \ifnum\outputpenalty>-2000 \else\dosupereject\fi}

%
 
\def\lA{{\cal A}} 
\def\lB{{\cal B}} 
 
\def\lD{{\cal D}} 
 
\def\lF{{\cal F}}

\def\lL{{\cal L}} 
\def\lM{{\cal M}}

\def\lS{{\cal S}}

\def\mA{{\bf A}} 
\def\mB{{\bf B}} 
\def\mC{{\bf C}}

\def\mL{{\bf L}} 
\def\mM{{\bf M}} 
\def\mN{{\bf N}}

\def\mQ{{\bf Q}} 
\def\mR{{\bf R}}

\def\mZ{{\bf Z}}

\def\mj{{\bf j}}

\def\msnop{{\bf p}}


\def\mx{{\bf x}} 
\def\my{{\bf y}}

\def\mnul{{\bf 0}} 



\def\vF{{\sf F}}

\def\va{{\sf a}} 
\def\vb{{\sf b}}

\def\ve{{\sf e}} 
\def\vf{{\sf f}} 
\def\vg{{\sf g}}

\def\vk{{\sf k}}

\def\vs{{\sf s}} 
 
\def\vu{{\sf u}} 
\def\vv{{\sf v}} 
\def\vw{{\sf w}}

\def\vnul{{\sf 0}}

%

\def\malpha{{\mib \alpha}} 
\def\mbeta{{\mib \beta}}

\def\mzeta{{\mib \zeta}} 
\def\meta{{\mib \eta}}

\def\mlambda{{\mib \lambda}} 
\def\mmu{{\mib \mu}} 
\def\mnu{{\mib \nu}} 
\def\mxi{{\mib \xi}} 
\def\mpi{{\mib \pi}}

%

\def\mpL{{\bsy L}}

%
%
%
\mathchardef\kvtpr="2D02
\mathchardef\kvadratic="0D03
\mathchardef\punikv="0D04
\mathchardef\rnsmjer="0D08
\mathchardef\rpsmjer="0D09
\mathchardef\jkon="3D13
\mathchardef\preth="3D34
\mathchardef\mj="3D36
\mathchardef\mave="3D37
\mathchardef\sljed="3D3C
\mathchardef\vj="3D3E
\mathchardef\vema="3D3F
\mathchardef\trokutic="0D4D
\mathchardef\punitrok="0D4E
\mathchardef\cul="3D62
\mathchardef\ksd="3D63
\mathchardef\mmmanje="3D6E
\mathchardef\vvvece="3D6F
\mathchardef\lgkut="0D70
\mathchardef\dgkut="0D71
\mathchardef\dirsuma="2D75
\mathchardef\ldkut="0D78
\mathchardef\ddkut="0D79
\mathchardef\ominus="0D7F

%
%
\mathchardef\hkprec="0E7D
\mathchardef\hrprec="0E7E






\let\Fou=\lF
\let\N=\mN
\let\R=\mR

\def\rC{{\rm C}}
\def\rL{{\rm L}}


\def\Rd{{{\bf R}^{d}}}

\def\Rpl{{{\bf R}^{+}}}


\let\ptrok=\punitrok
\let\cul=\kul


 
\def\Im{{\sf Im\thinspace}} 
\def\tr{{\sf tr}}

\def\Int{{\sf Int\thinspace}}

\def\Re{{\sf Re}\,}


\def\pS#1{{{\cal S}(#1)}}

\def\L#1{{{\rm L}(#1)}}
\def\LL#1{{{\rm L}}}

\def\Lb#1{{{\rm L}^\infty(#1)}}


\def\Ld#1{{{\rm L}^{2}(#1)}}

\def\Ll#1#2{{{\rm L}^{#1}_{{\rm loc}}(#2)}}

\def\Ldl#1{{{\rm L}^{2}_{{\rm loc}}(#1)}}




\def\M#1#2{{\rm M}_{{#1}\times{#2}}}


\def\H#1#2{{{\rm H}^{#1}(#2)}} 
\def\Hl#1#2{{{\rm H}^{#1}_{{\rm loc}}(#2)}}

\def\Lj#1{{{\rm L}^{1}(#1)}}

\def\Lb#1{{{\rm L}^{\infty}(#1)}} 

\def\Ljl#1{{{\rm L}^{1}_{{\rm loc}}(#1)}}

\def\Cbc#1{{{\rm C}^{\infty}_{c}(#1)}}

\def\Cb#1{{{\rm C}^{\infty}(#1)}} 
\def\Cp#1{{{\rm C}(#1)}}

\def\Cc#1{{{\rm C}_c(#1)}} 
 
\def\pCd#1#2{{{\rm C}_{#1}(#2)}}

\def\Cnl#1{{{\rm C}_0(#1)}}

\def\Kz#1#2{{{\rm K}[{#1},{#2}]}}
\def\Ko#1#2{{{\rm K}({#1},{#2})}}

\def\Sdmj{{\rm S}^{d-1}}



\def\supp{{\rm supp\,}}

\def\sign{{\rm sign\,}}


\def\rest#1{_{\displaystyle |_{#1}}}

\def\nor#1#2{{\| #1 \|}_{#2}}	

\def\Dup#1#2{\langle#1,#2\rangle}


\def\oi#1#2{\langle#1,#2\rangle}		
\def\ozi#1#2{\langle#1,#2]}
\def\zoi#1#2{[#1,#2\rangle}
\def\zi#1#2{[#1,#2]}				


\def\str{\longrightarrow} 
 
\def\Dstr{\relbar\joinrel\longrightarrow}

\def\dscon{\relbar\joinrel\rightharpoonup} 
\def\Dscon{\relbar\joinrel\dscon}

\def\povlaci{\quad\Longrightarrow\quad}
\def\akko{\qquad\Longleftrightarrow\qquad}
 
 
\def\eps{\varepsilon} 
\def\ph{\varphi}

\def\svaki#1{(\forall\,#1)}
\def\postoji#1{(\exists\,#1)}

\chardef\sS="19

\def\ss#1{\;(\hbox{\rm ss } #1)}




\def\kraj{\vrule height0dd width0,3em depth1ex}


\def\povrhsk#1{\smash{
        \mathop{\;\Dscon\;}\limits^{#1}}}

\def\putover#1{\mathop{\vbox{\ialign{##\crcr\noalign{\kern0pt}
             $\hfil\displaystyle{#1}\hfil$\crcr}}}\limits}    



%

\def\ss#1{\;(\hbox{\rm a.e. } #1)}

\newcount\btmp	
\newcount\btm	\btm=0 
\newcount\btma	\btma=0 
\newcount\blm	\blm=0
\newcount\blma	\blma=0
\newcount\bkr	\bkr=0 
\newcount\btv		\btv=0 
\newcount\brz		
\newcount\bfor	\bfor=0		
\newcount\bodj	\bodj=0		

\outer\long\def\tm#1{\medbreak	\global\advance\btm by1 
\noindent {\bf Theorem \number\btm.\enspace}{\sl #1 \hfill\kraj   \par\medbreak}} 
 
\outer\long\def\tma#1{\medbreak	\global\advance\btma by1 
\noindent {\bf Theorem \uppercase\expandafter{\romannumeral\btma}.\enspace}{\sl #1 \hfill\kraj   \par\medbreak}} 

\outer\long\def\lm#1{\medbreak	\global\advance\blm by1 
\noindent {\bf Lemma \number\blm.\enspace}{\sl #1 \hfill\kraj   \par\medbreak}} 
 
\outer\long\def\lma#1{\medbreak	\global\advance\blma by1 
\noindent {\bf Lemma \uppercase\expandafter{\romannumeral\blma}.\enspace}{\sl #1 \hfill\kraj   \par\medbreak}} 
 
\outer\long\def\kr#1{\medbreak	\global\advance\bkr by1 
\noindent {\bf Corollary \number\bkr.\enspace}{\sl #1 \hfill\kraj   \par\medbreak}} 

\outer\long\def\tmd#1#2{\medbreak	\global\advance\btm by1 
\noindent {\bf Theorem \number\btm.\enspace}{\sl #1  \par\smallbreak} 
\noindent {\sf Dem.}\enspace #2 \par \hfill \hbox{\bf Q.E.D.}\par\medbreak}

\outer\long\def\lmd#1#2{\medbreak	\global\advance\blm by1 
\noindent {\bf Lemma \number\blm.\enspace}{\sl #1  \par\smallbreak} 
\noindent {\sf Dem.}\enspace #2 \par \hfill \hbox{\bf Q.E.D.}\par\medbreak}  
 
\outer\long\def\krd#1#2{\medbreak	\global\advance\bkr by1 
\noindent {\bf Corollary \number\bkr.\enspace}{\sl #1  \par\smallbreak} 
\noindent {\sf Dem.}\enspace #2 \par \hfill \hbox{\bf Q.E.D.}\par\medbreak}  
 
\long\def\tvd#1#2{\medbreak	\global\advance\btv by1 
\noindent {\bf Claim \number\btv.\enspace}{\sl #1  \par\smallbreak} 
\noindent {\sf Dem.}\enspace #2 \par \hfill \hbox{\ptrok}\par\medbreak}

\long\def\razn#1{\smallbreak \global\advance\brz by1
\noindent{\bf \number\brz.} #1\par\medbreak}

\long\def\sazetak#1{{\garmond\narrower\narrower\smallskip
#1 \hfill\par\medbreak }}

\pageno=0

\ 
\vskip2mm
\par

{\cicero \bf \centerline{\pisac}}
\vskip8mm

{\mittel \bf \centreline{\naslova}} 
\bigskip

{\mittel \bf \centreline{\naslovb}} 
\vfill
\ifx\draft\undefined \else
{\hyppertertia \sf \draft} \vfill	\fi

{\cicero \bf  \centerline{Abstract}}
\vskip6mm

\hypocicero

\ifx\draft\undefined \def\draft{}\fi

\input epsf

\def\kobb{{\rm K_{0,\infty}}}
\def\kob{{{\rm K_{0,\infty}}}(\Rd)} 
\def\Ldl#1{{{\rm L}^{2}_{{\rm loc}}(#1)}} 

\def\domena{\Omega}

\def\binom#1#2{{#1\choose #2}}
\def\fp#1{\widehat{#1}}

\def\hmv{\mmu_H}
\def\hm{\mu_H}
\def\pkmv{\mmu_{sc}}
\def\pkm{\mu_{sc}}
\def\jhm{\mu_\kobb}
\def\jhmv{\mmu_\kobb}
\def\Sdmj{{\rm S}^{d-1}}
\def\Cr{\mC^r}
\def\M#1#2{{\rm M}_{#1\times #2}}

\def\Mx#1#2{{\rm M_{#1}}(#2)}
\def\mmp{{\bf p}}

\newcount\brpr \brpr=0
\long\def\Pr#1{\medbreak	\global\advance\brpr by1
\noindent {\bf Example \number\brpr.\enspace}{ #1 \hfill\kraj   
\par\medbreak} }

\def\Rdz{\R^d_\ast}
\def\buc{\rC_{ub}(\Rd)} 
\def\Dupp#1#2{\left\langle#1,#2\right\rangle}

\def\Cq{\mC^q}
\newcount\bnap \bnap=0
\long\def\Napbr#1{{\medbreak \global\advance\bnap by1
\noindent {\bf Remark \number\bnap.\enspace}{ #1 \hfill\kraj   \par\medbreak} }}

\def\LLd{{{\rm L}^2}}

\def\pkob#1{{{\rm K_{0,\infty}}}(#1)} 

\font\sc=cmcsc10

\newread\auX
\immediate\openin\auX=\jobname.auX
\ifeof\auX \message{! No file \jobname.auX;}
\else \input \jobname.auX \immediate\closein\auX \fi
\newwrite\auX
\immediate\openout\auX=\jobname.auX

\def\label#1#2{\immediate\write\auX%
{\noexpand\def\expandafter\noexpand\csname#2\endcsname{#1}}}
\def\ref#1{%
\ifundefined{#1}\message{! No ref. to #1;}%
\else\csname#1\endcsname\fi}

\def\ifundefined#1{\expandafter\ifx\csname#1\endcsname\relax}

\newcount\bite \bite=0

\let\itemm=\item

\def\ite#1{\itemm{[\global\advance\bite by1 \number\bite]\def\brite{\number\bite}\label\brite{#1}}}


\def\rfor#1{\global\advance\bfor by1 \number\bodj.\number\bfor \def\brfor{\number\bodj.\number\bfor}\label\brfor{#1}}

\def\cfor#1{%
\ifundefined{#1}\message{! No formula number to #1;}%
\else\csname#1\endcsname\fi}


\newcount\bodj  \bodj=0         

\def\odjeljak#1{\bigskip{\cicero \bf\par\bigskip\centreline{ \global\advance\bodj by1 \number\bodj. #1}\nobreak \medskip\nobreak \global\bfor=0} }

\sazetak{
Microlocal defect functionals (H-measures, H-distributions,
semiclassical measures, etc.) are objects which determine, in some 
sense, the lack of strong compactness for weakly convergent ${\rm L}^p$ 
sequences. Recently, Luc Tartar introduced one-scale H-measures, a 
generalisation of H-measures with a characteristic length, which also
comprehend the notion of semiclassical measures. 
 
We present a self-contained introduction to one-scale H-measures, carrying out
some alternative proofs, and strengthening some results, comparing these
objects to known microlocal defect functionals. Furthermore,
we improve and generalise Tartar's localisation principle for 
these objects from which we are able to derive the known 
localisation principles for both H-measures and semiclassical measures.
Moreover, we develop a variant of compactness by compensation suitable 
for equations with a characteristic length.
}

\vskip1cm
\bigskip

\vfill

\noindent {\bf Keywords:} semiclassical measure, H-measure, one-scale 
H-measure, localisation principle, compactness by compensation
\smallskip

\noindent {\bf Mathematics subject classification: 35B27, 35Q40, 35S05,
46G10}

\vfill

\settabs 8 \columns

\+&Department of Mathematics  	&	&&&   \cr
\+&Faculty of Science           &	&&&   \cr
\+&University of Zagreb       	&&	&&& University of Dubrovnik	 \cr
\+&Bijeni\v{c}ka cesta 30     	&&	&&& \'Cira Cari\'ca 4	 \cr
\+&Zagreb, Croatia		&&	&&& Dubrovnik, Croatia		\cr
\vskip1mm
\+{\tt nenad@math.hr} &&{\tt maerceg@math.hr} 	&&&&{\tt mlazar@unidu.hr}  \cr

\vskip2cm

\vfill
{\borgis \sl  This work has been supported in part by Croatian Science
Foundation under the project 9780 WeConMApp, by University of Zagreb
trough grant 4.1.2.14/2013, by the bilateral Croatian--Montenegro project {\it
Transport in highly heterogeneous media}, as well as by the DAAD
project {\it Centre of Excellence for Applications of Mathematics}.
}
\vskip5mm
\rightline{\borgis\today}

\smallskip

\eject

\pageno=0
\hbox{\ }

\vfill
\centreline{{\mittel\bf Contents}}
\bigskip

\hbox to\hsize{{\bf 1. Introduction}		\dotfill	1}
\item{} Overview
\item{} Notation
\medskip

\hbox to\hsize{{\bf 2. Overview of H-measures and semiclassical measures}	\dotfill	2}
\item{} H-measures
\item{} Semiclassical measures
\item{} Relation between H-measures and semiclassical measures
\item{} First examples
\medskip

\hbox to\hsize{{\bf 3. Tartar's one-scale H-measures}	\dotfill	6}
\item{} Test functions
\item{} Preliminary lemmata
\item{} The existence
\item{} First properties
\medskip

\hbox to\hsize{{\bf 4. Localisation principle}	\dotfill	17}
\item{} A notable condition
\item{} The case $\omega_n=\eps_n$
\item{} The case $c:=\lim_n{\eps_n\over\omega_n}\in\zi0\infty$
\medskip

\hbox to\hsize{{\bf 5. Some applications}	\dotfill	26}
\item{} Localisation principles revisited
\item{} Compactness by compensation with a characteristic length
\item{} Concluding remarks
\medskip

\hbox to\hsize{{\bf References}			\dotfill	30}

\vfill
\vfill
\eject


\odjeljak{Introduction}

\podnaslov{Overview}

In various situations concerning partial differential equations one
often encounters ${\rm L}^2$ weakly converging sequences, which do not
converge strongly. For such a sequence $(u_n)$ it is natural to consider
a ${\rm L}^1$ bounded sequence $|u_n|^2$ which, in general, does not
converge weakly in ${\rm L}^1$, but only weakly $\ast$ in the space of
bounded Radon measures ($\lM_b={\rm C}'_0$), to a defect measure $\nu$.

Essentially, there are two distinctive types of non-compact sequences,
with typical defect measures $\nu$ (for these two examples we fix 
$\ph\in\Cbc\Rd$ such that $\nor\ph{\Ld\Rd}=1$, $\eps_n\to 0^+$, and note 
that in both cases below one has $u_n\dscon0$):

\item{a)} {\sl concentration}\/: $u_n(\mx):=\eps_n^{-d/2}\ph\bigl({\mx-\mx_0 
\over\eps_n}\bigr)$, where $\nu=\delta_{\mx_0}$, and

\item{b)} {\sl oscillation}\/: $u_n(\mx):=\ph(\mx)e^{{2\pi i\mx\cdot
\mxi\over\eps_n}}$, where $\nu$ is actually equal to $|\ph|^2\lambda$ (i.e.~to 
the measure having density $|\ph|^2$ with respect to the 
Lebesgue measure $\lambda$ on $\Rd$).

\noindent However, defect measures are not enough to completely determine 
the difference between various non-compact sequences (e.g.~the direction 
and frequency of oscillations). This lack of information could be (partially) overcome
by introducing objects in full phase space such as H-measures and 
semiclassical measures.  
\smallskip

H-measures, as originally introduced a quarter of century ago by Luc Tartar [\ref{Tprse}] 
and (independently) Patrick G\' erard [\ref{Gmd}], are Radon measures 
on the cospherical bundle $\Omega\times\Sdmj$ over a domain $\Omega\subseteq\Rd$. 
Since their introduction they saw quite a number of successful applications,
many of which depend on the so called {\sl localisation principle}, which is
closely related to the generalisation of compactness by compensation method
to variable coefficients. Let us mention the applications to small-amplitude 
homogenisation [\ref{Tprse}, \ref {ALjmaa}],  to the  control theory [\ref{Burq}, \ref{LZ}], 
to explicit formul\ae\ and bounds in homogenisation  [\ref{Tprse}, \ref {ALrwa}],
to the existence of entropy solutions in the theory of conservation laws [\ref{EJParma}],
and the velocity averaging results [\ref{Gmd}, \ref{LMdp}].
This property is also crucial in applications of the {\sl propagation principle}, 
as it allows certain control over the support of H-measure ([\ref{Tprse}, \ref{Asym}]) .
\smallskip

Related to H-measures, but certainly different objects, are 
{\sl semiclassical measures}\/ first introduced by G\'erard
[\ref{Gmsc}], and later renamed by Pierre-Luis Lions and Thierry 
Paul [\ref{LP}] as {\sl Wigner measures}. 
They are Radon measures on the cotangential bundle $\Omega\times\Rd$
over a domain $\Omega\subseteq\Rd$.

In contrast to the H-measures, semiclassical measures depend upon 
characteristic length $(\omega_n)$, $\omega_n\to 0^+$, which makes 
them more suitable objects for the problems where such a characteristic 
length naturally appears, often related to highly oscillating problems
for partial differential equations (such as the homogenisation limit [\ref{GMMP}, \ref{LP}],
to microlocal energy density for (semi)linear wave equation [\ref{Goce}, \ref{FrInt}],
the semiclassical limit of Schr\" odinger equations [\ref{CFKMS}, \ref{Zh}],
or some other problems related to the quantum theory). 

While the {\sl localisation principle}\/ for semiclassical measures was
mentioned already in the first papers [\ref{Gmsc}, \ref{LP}]
(sometimes called also the {\sl elliptic regularity}\/ [\ref{Burq}]),
it has been playing a less prominent role in applications.
\medskip

In the next section we recall the definitions of both H-measures and
semiclassical measures, state the corresponding localisation principles,
and discuss their differences in the light of some examples.
The third section is devoted to the objects recently introduced by
Luc Tartar [\ref{Tgth}], which we call  {\sl one-scale H-measures}, in line
with the term {\sl multi-scale H-measures}\/ used recently in  [\ref{TmsHm}].
As opposed to the historical approach in [\ref{Tgth}], we have chosen some
variant proofs, and for the sake of completeness provided all the necessary
details.
The following section is dedicated to the localisation principle for one-scale 
H-measures, where we significantly extend Tartar's result from [\ref{Tgth}].
In the last section we show how the localisation principles, both
for H-measures and semiclassical measures, can be derived from the
obtained corresponding principle for the one-scale H-measures. 
We also apply the localisation principle for one-scale H-measures in developing 
a variant of compactness by compensation suitable for equations containing
a characteristic length.
\medskip
\break

\podnaslov{Notation}

Throughout the paper we denote by $\otimes$ the tensor product of vectors 
on $\Cr$, defined by $(\va\otimes\vb) \vv=(\vv\cdot\vb)\va$, where $\cdot$
stands for the (complex) scalar product ($\va\cdot\vb:=\sum_{i=1}^r a_i\bar  b_i$),
resulting in $[\va\otimes \vb]_{ij}=a_i\bar b_j$,
while $\kvtpr$ denotes the tensor product of functions in different variables. 
By $\Dup \cdot\cdot$ we denote any sesquilinear dual product, which we
take to be antilinear in the first variable, and linear in the second.
When matrix functions appear as both arguments in a dual product, 
we interpret it as 
$\Dup\mA\mB = \int\mB\cdot\mA = \int\tr(\mB\mA^\ast)$.
In order to have all formul{\ae} written in matrix notation, we define
dual product also for the case when matrix function appear in the first 
and scalar function in the second argument as 
$\Dup\mA\ph = [\Dup{A^{ij}}\ph]_{ij}$.

By $\mpi(\mxi) := {\mxi\over|\mxi|}$ we denote the projection
on $\Rdz:=\Rd\setminus\{\vnul\}$ along rays to the unit sphere $\Sdmj$. 
The open (closed) ball in $\Rd$ or $\Cr$ centred at point $\mx$ 
with radius $r>0$ we will denote by $\Ko\mx r$ ($\Kz\mx r$).

The Fourier transform, defined as 
$\hat \vu(\mxi):=\Fou \vu(\mxi) :=\int_\Rd e^{-2\pi i \mxi\cdot \mx}
\vu(\mx)\,d\mx$, and its inverse as $(\vu)^\vee(\mxi):=\bar\Fou \vu(\mxi) 
:=\int_\Rd e^{2\pi i \mxi\cdot \mx}\vu(\mx)\,d\mx$, allows for an elegant 
definition of Sobolev spaces on the whole space $\Rd$, in particular those 
modelled on $\Ld\Rd$, which is sufficient for our present purpose. Indeed, 
for $s\in\mR$ we take
$$
\H s\Rd := \Bigl\{u\in\lS' : \Bigl(1+|\mxi|^{|s|}\Bigr)^{\sign s}\hat u \in \Ld\Rd \Bigr\} \;,
$$
with the norm $\nor u{\H s\Rd} := \nor{ \bigl(1+|\mxi|^{|s|}\bigr)^{\sign s}\hat u}{\Ld\Rd}$.

However, in the sequel we shall also deal with functions defined on an open 
$\Omega\subseteq\Rd$. In order to apply the Fourier transform,
as well as to simplify other arguments, we shall identify such a function
by its extension by zero to the whole $\Rd$.

The local Sobolev spaces we are going to need are
$$
\Hl s\Omega := \Bigl\{u\in\lD'(\Omega) : \svaki{\ph\in\Cbc\Omega}\; \ph u \in \H s\Rd \Bigr\} \;,
$$
with the above mentioned identification of a distribution and its extension by zero assumed.
Naturally, $\Hl s\Omega$ will be endowed with the weakest topology in which every
mapping $u\mapsto \ph u$ is continuous (cf.~[\ref{AL}, p.~1207] and the references mentioned there).
In particular, the space $\Ldl{\Omega;\Cr}$ is equipped with the standard Fr\'echet locally convex topology
(cf.~[\ref{ABu}] and the references therein). Hence, a subset of 
$\Ldl{\Omega;\Cr}$ is bounded if and only if it is bounded in the sense of 
seminorms which generate the corresponding locally convex topology (which 
is a stronger notion then metric boundedness).

\odjeljak{Overview of H-measures and semiclassical measures}


\podnaslov{H-measures}

For the sake of generality we shall deal with local spaces. In this context,
the following theorem can be stated (cf.~[\ref{Gmd}, Theorem 1] or [\ref{Tprse}, Theorem 1.1]):

\tm{
{\bf (existence of H-measures)} \  
For a weakly converging sequence $\vu_n\dscon\vnul$ in $\Ldl{\Omega;\Cr}$, 
there exists a subsequence\/ $(\vu_{n'})$ and an $r \times r$
hermitian non-negative matrix Radon measure\/ $\hmv$ on\/ 
$\Omega\times \Sdmj$ such that for any $\ph_1,\ph_2 \in \Cc\Omega$ 
and $\psi \in \Cp{\Sdmj}$ one has:
$$
\eqalign{
\lim_{n'} \int\limits_\Rd \Bigl(\widehat{\ph_1 \vu_{n'}}(\mxi) \otimes
	\widehat{\ph_2 \vu_{n'}}(\mxi)\Bigr) \psi\bigg({\mxi \over |\mxi|}
	\bigg) \,d\mxi 
	& = \Dup{\hmv}{(\ph_1\bar\ph_2) \kvtpr\psi} \cr
& = \int\limits_{\Omega\times\Sdmj} \ph_1(\mx)\bar\ph_2(\mx)\psi(\mxi)
		 \,d\bar\hmv(\mx,\mxi) \,. \cr
}
$$
The above measure $\hmv$ is called\/ {\rm the H-measure} associated to 
the (sub)sequence $(\vu_{n'})$.
}

We shall often abuse the notation and terminology, assuming that 
we have already passed to a subsequence determining an H-measure.
When the {\sl whole}\/ sequence admits the H-measure (i.e.~the definition 
is valid without passing to a subsequence), we say that the 
sequence is {\sl pure}.

The corresponding defect measure of the weakly converging (towards $\vnul$)
sequence $(\vu_n)$ can be obtained simply by integrating the H-measure with 
respect to the Fourier space variable $\mxi$.
As an immediate consequence we have that 
any H-measure associated to a strongly convergent sequence in 
$\Ldl{\Omega;\Cr}$ is necessarily zero, and vice versa, if the
H-measure is trivial, then the corresponding (sub)sequence 
converges strongly in $\Ldl{\Omega;\Cr}$.

The majority of successful applications of H-measures use the localisation
principle (cf.~[\ref{Tprse}], [\ref{Asym}]).

\tm{
{\bf (localisation principle for H-measures)} \ 
Let\/ $\vu_n\dscon\vnul$ in $\Ldl{\Omega;\Cr}$, and let for a given 
$m\in\N$
$$
\sum_{|\malpha|\mj m} \partial_\malpha(\mA^{\malpha}\vu_{n}) 
	\str \vnul \quad \hbox{strongly in} \quad \Hl{-m}{\Omega;\Cq} \,,
\leqno(\rfor{locprEqHmA})
$$
where $\mA^\malpha\in\Cp{\Omega;\M qr(\mC)}$ and 
$\partial_\malpha = {\partial^{\alpha_1}\over \partial x^1} \dots 
{\partial^{\alpha_d}\over \partial x^d}$ denotes partial derivatives 
in variable $\mx$ in the physical space. 

Then for the associated H-measure $\hmv$ we have
$$
\mmp_{pr}\hmv^\top = \mnul \,,
$$
where
$$
\mmp_{pr}(\mx,\mxi) := \sum_{|\malpha|=m}(2\pi i)^m\Bigl({\mxi\over |\mxi|}\Bigr)^\malpha
	\mA^\malpha(\mx)
\leqno(\rfor{Symbppr})
$$
is the principal symbol of the differential operator in (\cfor{locprEqHmA}). 
}

This result implies that the support of $\hmv$ is contained in the set
$$
\Sigma_{\mmp_{pr}} := \Big\{(\mx,\mxi) \in \Omega\times \Sdmj : 
	\hbox{rank}\thinspace\mmp_{pr}(\mx,\mxi) < r \Big\}
$$
of points where $\mmp_{pr}(\mx,\mxi)$ is not left invertible.

\Napbr{%
The main difference between Tartar's and G\'erard's approach in the construction of 
H-measures is in the choice of the space of test functions. Tartar based the proof on  a 
commutation lemma (known under the name {\sl First commutation lemma}) for which 
only the continuity of test functions is sufficient. On the other hand, G\'erard 
used the theory of pseudodifferential operators in which the commutation lemma is 
trivially obtained (the difference of pseudodifferential operators of the same order 
and with the same principal symbol is compact), but requiring smooth test functions. 
}

Recently, H-measures were extended in two directions: to {\sl parabolic 
H-measures}\/ [\ref{ALjmaa}, \ref{AL}], adapted to a different scaling and projection $\mpi$, 
and to {\sl H-distributions}\/ [\ref{AMaaa}], allowing the treatment of $\rL^p$
weakly converging sequences for $p\in\oi1\infty$ (see also [\ref{FRin}]).


\podnaslov{Semiclassical measures}

We present the existence result for semiclassical measures in a 
simpler, but equivalent form to the original G\'erard's definition [\ref{Gmsc}], 
without introducing the notion of (semiclassical) pseudodifferential operators 
(cf.~[\ref{Tgth}, Chapter 32]).

\tm{
{\bf (existence of semiclassical measures)} \
If\/ $(\vu_n)$ is a bounded sequence in\/ $\Ldl{\Omega;\Cr}$, 
and $(\omega_n)$ a sequence of positive numbers such that 
$\omega_n\to 0^+$, then there exists a subsequence\/ $(\vu_{n'})$ and 
an $r \times r$ hermitian non-negative matrix Radon measure\/ 
$\pkmv^{(\omega_{n'})}$ on\/ $\Omega\times\Rd$ such that for any
$\ph_1,\ph_2\in\Cbc\Omega$ and $\psi \in \lS(\Rd)$ one has:
$$
\eqalign{
\lim_{n'} \int\limits_\Rd \Bigl(\widehat{\ph_1 \vu_{n'}}(\mxi) \otimes
	\widehat{\ph_2 \vu_{n'}}(\mxi)\Bigr) \psi(\omega_{n'}\mxi) \,d\mxi 
	& = \Dup{\pkmv^{(\omega_{n'})}}{(\ph_1\bar\ph_2) \kvtpr\psi} \cr
& = \int\limits_{\Omega\times\Rd} \ph_1(\mx)\bar\ph_2(\mx)\psi(\mxi)
		 \,d\bar\pkmv^{(\omega_{n'})}(\mx,\mxi) \,. \cr
}
$$
The above measure $\pkmv^{(\omega_{n'})}$ is called\/ {\rm the semiclassical
measure (with characteristic length $(\omega_{n'})$)}\/ associated to the
(sub)sequence $(\vu_{n'})$.
}

\Napbr{%
For a bounded sequence in $\Ld{\Omega;\Cr}$, both the associated 
H-measure and the semiclassical measure are bounded Radon measures 
and formul{\ae} in Theorem 1 and Theorem 3 make sense for 
$\ph_1, \ph_2 \in \Cnl{\Omega}$.
}

As it is done with H-measures, we shall often abuse the notation and 
terminology, assuming that we have already passed to a subsequence 
determining a semiclassical measure, and reduce the notation to 
$\pkmv=\pkmv^{(\omega_n)}$, except in situations where it will 
not be clear which characteristic length do we use. When the 
{\sl whole}\/ sequence admits a semiclassical measure with characteristic 
length $(\omega_n)$ (i.e.~the definition is valid without passing 
to a subsequence), we say that the sequence is {\sl $(\omega_n)$-pure}
(in terms of semiclassical measures).

In the definition of semiclassical measures it is not necessary that 
the observed sequence tends (weakly) to zero, contrary to the definition
of H-measures. This is the case because here, due to the characteristic length 
$(\omega_n)$, we have a stronger variant of the commutation lemma (see Lemma 2
in the next section). Consequently, in some situations semiclassical measures 
can be used not only to prove strong compactness of a sequence, but also to establish 
that the limit is zero. 
Moreover, we can decompose a semiclassical measure into two parts, one 
depending on the value of the limit, and the other associated to the 
sequence converging to zero. Namely, if $\vu_n\dscon\vu$ in $\Ldl{\Omega;\Cr}$, 
then 
$$
\pkmv = (\vu\otimes\vu)\lambda\kvtpr\delta_\vnul + \mnu_{sc}	\,,
\leqno(\rfor{decompPkm})
$$
where $\mnu_{sc}$ is the semiclassical measure, with the same characteristic 
length as $\pkmv$, associated to the sequence $(\vu_n-\vu)$, while $\lambda$ is the 
Lebesgue measure in $\mx$, and $\delta_\vnul$ the Dirac mass in $\mxi=\vnul$.

In order to obtain the defect measure from a semiclassical measure,
an additional assumption is required. First, recall that a sequence 
$(\vu_n)$ in $\Ldl{\Omega;\Cr}$ is called {\sl $(\omega_n)$-oscillatory}, 
where $\omega_n\to 0^+$, if for any $\ph\in\Cbc\Omega$
$$
\lim_{r\to\infty}\limsup_n \int\limits_{|\mxi|\vj{r\over\omega_n}}
	|\widehat{\ph\vu_n}|^2\,d\mxi = 0 \,.
$$
An interpretation of this property is that information of the observed 
sequence is on the scale not greater than ${1\over\omega_n}$ in the 
Fourier space, so the semiclassical measure with characteristic length 
$(\omega_n)$ is able of capturing all the information and nothing is lost 
at infinity (cf.~[\ref{FrInt}, Def.~3.3] and [\ref{GMMP}, Def.~1.6]).
For example, sequences of concentration and oscillations from the introductory 
section are both $(\omega_n)$-oscillatory if and only if 
$\lim_n{\eps_n\over\omega_n}\in\ozi0\infty$. One can notice here certain ambiguity  
in terminology, suggesting that a concentration effect provides oscillations in some sense.     

If a sequence $(\vu_n)$ is bounded and $(\omega_n)$-oscillatory, then the
corresponding defect measure can be obtained by
integrating $\pkmv^{(\omega_n)}$ with respect to the Fourier space
variable $\mxi$.
As an immediate consequence we have that, if for a bounded sequence 
$(\vu_n)$ in $\Ldl{\Omega;\Cr}$ there exists a sequence $(\omega_n)$,
$\omega_n\to0^+$, such that $(\vu_n)$ is $(\omega_n)$-oscillatory and 
the semiclassical measure with characteristic length $(\omega_n)$ 
associated to $(\vu_n)$ is trivial, then $(\vu_n)$ converges strongly 
in $\Ldl{\Omega;\Cr}$.

As it was the case with H-measures, it is of great interest to localise
the support of semiclassical measures. The statement of the following
theorem is slightly different from those in [\ref{Burq}, \ref{Gmsc}, \ref{LP}],
as we have a system of differential relations, so the coefficients are matrix valued.
Since the proof follows by minor modifications of the usual proof, while 
in the sequel a more general result will be presented, we shall omit it.

\tm{
{\bf (localisation principle for semiclassical measures)} \ 
Let\/ $(\vu_n)$ be bounded in $\Ldl{\Omega;\Cr}$, $\omega_n\to 0^+$, and 
let for $m\in\N$
$$
\sum_{|\malpha|\mj m} \omega_n^{|\malpha|}\partial_\malpha(
	\mA^{\malpha}\vu_{n}) \str \vnul \quad \hbox{strongly in} 
	\quad \Ldl{\Omega;\Cq} \,,
$$
where $\mA^\malpha\in\Cb{\Omega;\M qr(\mC)}$.

Then for the associated semiclassical measure $\pkmv=\pkmv^{(\omega_{n'})}$
we have
$$
\mmp_{sc} \pkmv^\top = \mnul \,,
$$
where
$$
\mmp_{sc}(\mx,\mxi) := \sum_{|\malpha|\mj m}(2\pi i\mxi)^\malpha
	\mA^\malpha(\mx) \,.
$$ 
}

Analogously to the result for H-measures, this theorem constrains the support of
$\pkmv$ within the set of points where the symbol $\mmp_{sc}(\mx,\mxi)$ has no left 
inverse.

\medskip


\podnaslov{Relation between H-measures and semiclassical measures}

In order to compare the two types of objects, let us consider an 
elementary example where both the H-measure and the semiclassical measure
can be explicitly calculated (cf.~[\ref{Tgth}, Lemma 32.2]). 

\Pr{%
{\bf (oscillation: one characteristic length)} \
Let $\alpha>0$, $\vk\in\mZ^d\setminus\{\vnul\}$, and consider
$$
u_n(\mx):= e^{2\pi i n^\alpha\vk\cdot\mx}\dscon 0 \quad \hbox{in} 
	\quad \Ldl\Rd \,.
$$

The sequence $(u_n)$ is pure and the associated H-measure is 
$$
\hm = \lambda\kvtpr\delta_{{\vk \over |\vk|}} \,.
$$
Thus the H-measure contains important information on the direction of
oscillation $\vk$.

Furthermore, the sequence $(u_n)$ is $(\omega_n)$-pure for any
characteristic length $(\omega_n)$, and the corresponding semiclassical
measure depends on the characteristic length:
$$
\pkm^{(\omega_n)} = \lambda\kvtpr 
	\left\{
	\matrix{ 
	\delta_\vnul & , & \lim_n n^\alpha\omega_n = 0 \cr
	\delta_{c\vk} & , & \lim_n n^\alpha\omega_n = c\in\oi0\infty \cr
	0 & , & \lim_n n^\alpha\omega_n = \infty
	} \right. \,.
$$
Only in the case $\lim_n n^\alpha\omega_n\in\oi0\infty$, where the 
characteristic length of semiclassical measure is of the same order as the
wavelength of the observed sequence, the corresponding semiclassical
measure contains information on the direction of oscillation. However, 
contrary to the H-measure, which cannot comprehend the frequency, for
every $(\omega_n)$ the semiclassical measure contains information on
the frequency of oscillation. 
}

In the previous example one can notice that, in the case 
$\lim_n n^\alpha\omega_n\in\oi0\infty$, the H-measure can be derived from 
the corresponding semiclassical measure by taking the projection to the
unit sphere in the Fourier space. This can be justified by the well known result on the
relation of the above measures [\ref{Burq}, Proposition 4],
[\ref{FrInt}, Lemma 3.4].

\tm{
Let\/ $\vu_n\dscon\vu$ in\/ $\Ldl{\Omega;\Cr}$, $\omega_n\to 0^+$, and let 
$(\vu_n)$ be $(\omega_n)$-pure with the 
corresponding semiclassical measure $\pkmv=\pkmv^{(\omega_n)}$.

If\/ $\tr\pkmv(\Omega\times\{\vnul\}) = 0$, then $\vu=\vnul$.
Furthermore, if we additionally assume that $(\vu_n)$ is 
$(\omega_n)$-oscillatory, then the sequence $(\vu_n)$ is pure
(in the sense of H-measures), and for any $\ph\in\Cbc\Omega$ and 
$\psi\in\Cb{\Sdmj}$ one has:
$$
\Dup{\hmv}{\ph\kvtpr\psi} = \Dup{\pkmv}{\ph\kvtpr(\psi\circ\mpi)} \,,
$$
where $\hmv$ is the H-measure associated to $(\vu_n)$.
}

Going back to Example 1, it can be checked (either by straightforward calculation, 
or by using the interpretation of the  $(\omega_n)$-oscillatory property)
that the sequence $(u_n)$ is $(\omega_n)$-oscillatory in the case 
$\lim_n n^\alpha\omega_n\in\zoi0\infty$. Furthermore, with the exception of the 
case when this limit is zero, the set $\Rd\times\{\vnul\}$ is of zero measure with respect 
to the corresponding semiclassical measure as well. Hence, by the previous 
theorem, the comment after Example 1 is justified. 
Moreover, we can see that the previous theorem is sharp, since in other 
two cases the H-measure cannot be derived from the semiclassical
measure.   

Although one can think, after the above arguments, that semiclassical 
measures are more general objects then H-measures (taking into account all 
possible characteristic lengths), this is not the case. Indeed, it 
is not hard to construct an example in which (for any
characteristic length) one of the assumptions of the previous theorem is
not satisfied.

\Pr{%
{\bf (oscillation: two characteristic lengths)} \
Let $0<\alpha<\beta$, while $\vk, \vs\in\mZ^d\setminus\{\vnul\}$, and 
define two weakly convergent sequences in $\Ldl\Rd$
$$
\eqalign{
u_n(\mx) &:= e^{2\pi i n^\alpha\vk\cdot\mx}\dscon 0 \,, \cr
v_n(\mx) &:= e^{2\pi i n^\beta\vs\cdot\mx} \dscon 0 \,. \cr
}
$$
We are interested in the measures associated to the sequence $(u_n+v_n)$, 
which is pure and $(\omega_n)$-pure for any characteristic length $(\omega_n)$.

The corresponding measures $\hm$ and $\pkm^{(\omega_n)}$ are given by:
$$
\eqalign{
\hm &= \lambda\kvtpr\Bigl(\delta_{{\vk \over |\vk|}}+\delta_{{\vs
	\over |\vs|}}\Bigr) \,, \cr
\pkm^{(\omega_n)} &= \lambda\kvtpr 
	\left\{ \matrix{ 
	2\delta_\vnul & , & \lim_n n^\beta\omega_n = 0 \cr
	(\delta_{c\vs}+\delta_\vnul) & , & \lim_n n^\beta\omega_n = 
		c\in\oi{0}{\infty} \cr
	\delta_\vnul & , & \lim_n n^\beta\omega_n = \infty \ \&	 
		\ \lim_n n^\alpha\omega_n = 0 \cr
	\delta_{c\vk} & , & \lim_n n^\alpha\omega_n = c\in\oi{0}{\infty} \cr
	0 & , & \lim_n n^\alpha\omega_n = \infty \cr
	} \right. \,. \cr
}
$$
In neither of these cases the assumptions of the last theorem are 
satisfied, and one can not derive the H-measure from the semiclassical measure. 
}

To conclude, we can say that semiclassical measures lose a part of information
for some (or even all) characteristic lengths, while H-measures cannot catch
frequency of an observed sequence. This was the motivation for Tartar to introduce
a new object, the one-scale H-measure, that overcomes the shortcomings of the 
preceding tools, and which we study in the next sections.


\odjeljak{Tartar's one-scale H-measures}

In the recent book [\ref{Tgth}, Chapter 32] Luc Tartar introduced improved
objects, having the advantages of both H-measures and semiclassical
measures, which he called {\sl H-measure variants with a characteristic
length}\/ there. However, in line with [\ref{TmsHm}], we call them 
{\sl one-scale H-measures}\/ in this paper.

\vfill
\eject

\podnaslov{Test functions}

While the theorems on existence of both H-measures and semiclassical
measures have the same form, the crucial difference is in the choice 
of test functions in the Fourier space.
In the former case functions from $\Cp\Sdmj$ are chosen (in fact, the continuous 
functions on $\Rdz:=\Rd\setminus\{\vnul\}$ being homogeneous of order zero), 
while in the latter functions from the Schwartz space $\pS{\Rd}$ are being used.

For the new variant, test functions are taken to be continuous functions on a
suitable compactification of $\Rdz$ in the Fourier space. More precisely, 
following Tartar, we denote by $\kob$ the compactification of $\Rdz$ obtained by placing 
two unit spheres, $\Sigma_0$ at the origin, and $\Sigma_\infty$ at infinity, 
equipped with the topology making it homeomorphic to the 
{\sl $d$-dimensional spherical shell}, as depicted on Figure 1.
The points on these two spheres we denote by $0^\ve$
and $\infty^\ve$, respectively, where $\ve\in \Sdmj$ is the direction 
of the point (i.e.~$0^\ve$ and $\infty^\ve$ are the {\sl endpoints}\/ 
of a ray passing through $\ve$). 
In the plane, ${\rm K_{0,\infty}}(\R^2)$ is homeomorphic to the closed annulus. 
$$
\epsfbox{slika1.mps}
$$
$$
\hbox{{\borgis{\bf Figure 1.} $\kob$, the compactification of $\Rdz$.}}
$$

As we have already seen in two previous examples, the semiclassical measures 
have a disadvantage of losing information at infinity and mixing 
various pieces of information at the origin, which motivated Tartar to 
introduce this compactification. The idea is to capture lost data on 
$\Sigma_0$ and $\Sigma_\infty$ (cf.~Example 4 below). 

For the sake of simplicity, we shall use the same notation for functions $\psi$ on
$\kob$ and their corresponding pullbacks $\sigma^\ast \psi = 
\psi\circ\sigma$, where $\sigma$ is the embedding of $\Rdz$ into $\kob$, while 
from the context it will be clear which object is used.

According to the previous identification, the space $\Cp\kob$ can be
interpreted as consisting of all continuous functions $\psi$ on $\Rdz$ 
such that there are functions $\psi_0, \psi_\infty\in\Cp\Sdmj$ satisfying 
$$
\eqalign{
& \psi(\mxi)-\psi_0\Bigl( {\mxi\over|\mxi|} \Bigr)\to 0, \quad 
	|\mxi|\to 0 \,, \cr
& \psi(\mxi)-\psi_\infty\Bigl( {\mxi\over|\mxi|} \Bigr)\to 0, \quad 
	|\mxi|\to\infty \,. \cr
}
\leqno(\rfor{psinib})
$$
Furthermore, for an arbitrary $\ve\in\Sdmj$ we have $\psi_0(\ve)=\psi(0^\ve)$
and $\psi_\infty(\ve)=\psi(\infty^\ve)$.
This notation will be kept throughout the paper: for an arbitrary
function $\psi$ from $\Cp\kob$, with indices 0 $(\psi_0)$ and $\infty$ $(\psi_\infty)$, 
we shall denote functions with the above properties .
Equipped by the $\rL^\infty$ norm, $\Cp\kob$ becomes a separable Banach space,
even a Banach algebra with usual multiplication of functions.

For $\psi\in\Cp\kob$, $\psi-\psi_0\circ\mpi$ is constant (with zero
value) on $\Sigma_0$, so it can be extended by continuity to the origin. 
Hence, we can say that $\psi-\psi_0\circ\mpi$ is contained in $\buc$ 
(the space of bounded uniformly continuous functions on $\Rd$), as both $\psi$ 
and $\psi_0\circ\mpi$ are bounded, while outside a sufficiently large 
compact set the difference can be approximated by 
$(\psi_\infty-\psi_0)\circ\mpi$. 
Moreover, if $\psi_0$ is constant than we immediately have that $\psi\in\buc$.

The space $\Cnl\Rd$ is continuously embedded in the space $\Cp\kob$. 
Indeed, a function $\psi\in\Cnl\Rd$ is continuous at the origin and tends 
to zero at infinity, so for constants $\psi_0\equiv \psi(\vnul)$ and 
$\psi_\infty\equiv 0$, $\psi$ satisfies the above conditions. It is even simpler 
to prove that $\mpi^\ast(\Cp\Sdmj):=\{\psi\circ\mpi \ : \ \psi\in\Cp\Sdmj\}\hookrightarrow\Cp\kob$, 
by taking $\psi_0:=\psi$ and $\psi_\infty:=\psi$. This discussion we summarise in the form of a lemma. 

\lm{
\item{$i)$} $\Cnl\Rd\hookrightarrow\Cp\kob$, and

\item{$ii)$} $\{\psi\circ\mpi \ : \ \psi\in\Cp\Sdmj\}\hookrightarrow\Cp\kob$.
}

As $\lS(\Rd) \hookrightarrow \Cnl\Rd$, it follows  that the space $\Cp\kob$ 
subsumes the spaces of test functions for both H-measures and semiclassical 
measures.

An example of a function from $\Cp\kob$ which is not entirely contained
in Lemma 1 and which will be important in the next section when discussing 
the localisation principle for one-scale H-measures, is given in the next example.

\Pr{%
Let $l,m\in\N$ and define 
$\psi^{\malpha}(\mxi):={\mxi^\malpha\over |\mxi|^l+|\mxi|^m}$ 
for $\malpha\in\N^d$, $l\mj |\malpha| \mj m$.
The cases $|\malpha|=l$ and $|\malpha|=m$ are of particular interest,
as they are not necessarily contained in the above lemma. 
Let us show that $\psi^{\malpha}$ is (a pullback of a function) in 
$\Cp\kob$.

It is obvious that $\psi^{\malpha}$ is continuous on $\Rdz$. For $l=m$,
$\psi^{\malpha}$ is homogeneous of order zero, so Lemma 1 gives the 
statement. On the other hand, for $l\mj m-1$ we need to make an analysis in
dependence on $\malpha$. For $|\malpha|=l$ we have 
$\psi^{\malpha}_0(\mxi)={\mxi^\malpha \over |\mxi|^l}$, while in the
case $|\malpha|\vj l+1$, $\psi^{\malpha}$ tends to zero at the origin,
which implies $\psi^{\malpha}_0\equiv 0$.
At $\Sigma_\infty$ we have an opposite situation. For 
$|\malpha|\mj m-1$, $\psi^{\malpha}$ tends to zero at infinity, which
implies $\psi^{\malpha}_\infty\equiv 0$, but 
$\psi^{\malpha}_\infty(\mxi)={\mxi^\malpha \over |\mxi|^m}$ for 
$|\malpha|=m$. 
}

\podnaslov{Preliminary lemmata}

The existence result for one-scale H-measures, resembling the one for H-measures,
relies on a variant of the First commutation lemma 
for which we  introduce two basic types of operators. For $\psi\in\buc$, 
as well as for $\psi\in\Cp\kob$, one can define {\sl the Fourier multiplier 
operator}\/:
$$
\lA_\psi : \Ld\Rd \str \Ld\Rd \,, \quad \lA_\psi\vu := (\psi\hat\vu)^\vee \,.
$$    
The definition is justified, as in both cases $\psi$ is in $\Lb\Rd$ (in the latter
case, we assume that $\psi$ is first restricted to $\Rdz$, and then interpreted as
an $\rL^\infty$ function). Also, for $\ph\in\Lb\Rd$ by 
$$
B_\ph : \Ld\Rd \str \Ld\Rd \,, \quad B_\ph\vu := \ph\vu \,,
$$ 
we denote the operator of multiplication by $\ph$. In both cases, the above operators 
are bounded on $\Ld\Rd$, with the norm equal to the $\rL^\infty$ norm of $\psi$, 
respectively $\ph$.

In [\ref{Tgth}, Lemma 32.4] the following result is given:

\lm{
Let\/ $\psi\in\buc$, $\ph\in\Cnl\Rd$, $\omega_n\to 0^+$, and 
denote\/ $\psi_n(\mxi) := \psi(\omega_n\mxi)$.
Then the commutator\/ $C_n:=[B_\ph, \lA_{\psi_n}] = B_\ph\lA_{\psi_n} - \lA_{\psi_n}B_\ph$
tends to zero in the operator norm on $\lL(\Ld\Rd)$.
}

\Napbr{%
Lemma 2 can be used to prove the known result for semiclassical
measures given in Theorem 3 not only for (infinitely) smooth, but 
merely for continuous test functions both in the physical and Fourier 
space [cf. \ref{TmsHm}].
}

In comparison to the classical First commutation lemma, where the
commutator is only a compact operator, here we have the convergence 
in the uniform (operator) norm to zero. This would allow us to weaken the
requirement that $(\vu_n)$ weakly converges to zero, allowing any limit 
if  test functions were taken from $\buc$ in the Fourier space. 
Unfortunately, functions from $\Cp\kob$ are not continuous at the 
origin, therefore not in $\buc$, so we cannot avoid the 
additional assumption that the weak limit of the observed sequence is zero.
For that reason the following variant of the previous lemma will be used in 
the definition of one-scale H-measures.

\lmd{
Let\/ $\psi\in\Cp\kob$, $\ph\in\Cnl\Rd$, $\omega_n\to 0^+$, and 
denote\/ $\psi_n(\mxi) := \psi(\omega_n\mxi)$.
Then the commutator can be expressed as a sum
$$
C_n:=[B_\ph, \lA_{\psi_n}] = \tilde C_n + K \,,
$$
where $K$ is a compact operator on $\Ld\Rd$, while $\tilde C_n \str 0$ in the operator norm
on $\lL(\Ld\Rd)$.
}
{
Since
$$
\lA_{\psi_n} = \lA_{\psi_n-\psi_0\circ\mpi} + \lA_{\psi_0\circ\mpi} \,,
$$
we have that $C_n = [B_\ph, \lA_{\psi_n-\psi_0\circ\mpi}] + 
[B_\ph, \lA_{\psi_0\circ\mpi}]=: \tilde C_n+K$. As $\psi - \psi_0\circ\mpi 
\in\buc$, we can apply Lemma 2 to this term in order to show $\tilde C_n\str 0$, while
$K=[B_\ph, \lA_{\psi_0\circ\mpi}]$ is compact by the classical 
First commutation lemma [\ref{Tprse}, Lemma 1.7].
}

In [\ref{Tgth}, Chapter 32] and [\ref{TmsHm}] one can find the proof of existence of one-scale H-measures
(albeit only for the weakly converging sequences in ${\rm L}^2$, and not in the local space), starting
from the corresponding H-measure in one dimension more. 
In the same references one can find the idea of an alternative proof based 
on Tartar's original proof [\ref{Tprse}], which uses the construction with Hilbert-Schmidt operators
[\ref{Tprse}, Lemma 1.10]. 
As we need several modifications, we shall refrain from stating that the construction
can be done along the same lines, and try to outline all the necessary steps. 

Another venue for the proof could be by using the Schwartz kernel theorem, but this would require
introduction of a differential structure on $\kob$, and the space of distributions on it,
which we choose to avoid here.
However, it could be said that we are providing a version of the kernel theorem in
this special case.

\lmd{Let $X$ and $Y$ be open unit balls in Euclidean spaces (of dimension $d$ and $d'$),
and let $B$ be a non-negative continuous bilinear form on $\Cc X\times \Cc Y$.
Then there exists a Radon measure $\mu\in\lM(X\times Y)$ such that
$$
\svaki{f\in\Cc X} \svaki{g\in\Cc Y}\quad  B(f,g)=\langle \mu, f\kvtpr g\rangle \;.
$$
Furthermore, the above remains valid if we replace ${\rm C}_c$ by ${\rm C}_0$, and
$\lM$ by $\lM_b$ (the space of bounded Radon measures, i.e.~the dual of Banach space ${\rm C}_0$).
}
{{\bf I. } {\sl Friedrichs mollifiers and nested balls}
\smallskip

Following Tartar's ideas, we take standard Friedrichs mollifiers
$$
(M_n f)(\mx):= \int_\Rd m_n(\mx,\mx')f(\mx')\,d\mx' \;,
$$
where the kernel $m_n(\mx,\mx'):=\rho_n(\mx-\mx')$, while $\rho_n(\mx):=n^d\rho(n\mx)$ and
$$
\rho(\mx):= C \chi_\Ko\vnul1(\mx) e^{-{1\over 1-|\mx|^2}}\;,
$$
with $C$ chosen such that $\int\rho=\int\rho_n = 1$ ($C$ depends on dimension $d$).

These kernels $m_n$ are non-negative, continuous, supported in $\{(\mx,\mx')\in X\times X:|\mx-\mx'|\mj1/n\}$,
and smoothing; in particular $M_n$ maps Lebesgue functions into smooth ones
(recall that we identify functions defined on subsets of $\Rd$ to their continuation by zero
to $\Rd$).

For an open unit ball $X=\Ko\vnul1$ in $\Rd$ we have an increasing sequence of compacts $K_m$
contained in it, such that their union is the whole ball, while $K_m \subseteq \Int K_{m+1}$.
For definiteness, let us take closed balls $K_m:=\Kz\vnul{1-1/m}$.
Recall that $\Cc X$ is a strict inductive limit of Banach spaces
${\rm C}_{K_m}(X):= \Bigl\{\ph\in\Cc X : \supp\ph\subseteq K_m \Bigr\}$.

Of course, we can repeat the same construction on $Y$, obtaining operators $N_n$ and kernels $n_n$, 
while the corresponding sequence of compacts we denote by $L_m$ (we also write $\rho'_n, C'$ 
instead of $\rho_n, C$).

The continuity of $B$, as spaces ${\rm C}_c$ are strict inductive limits, can be expressed by the
continuity of restrictions, i.e.~by
$$
\svaki{m\in\mN}\postoji{C_m>0}\svaki{f\in\pCd{K_m}X}\svaki{g\in\pCd{L_m}Y}\quad
	|B(f,g)| \mj C_m \nor f{\Lb{K_m}}\nor g{\Lb{L_m}} \;.
\leqno(\rfor{Bcont})
$$

\noindent {\bf II. } {\sl Construction of approximative Hilbert-Schmidt kernels}
\smallskip

Let us define a bilinear functional $B^m_n: \Ld{K_m}\times \Ld{L_m} \str\mC$ by the relation
$ B^m_n(f,g):=B(M_nf, N_ng)$. Note that for $n> m(m+1)$ the function $M_n f$ belongs to $\pCd{K_m}X$, 
and similarly for $N_n g$.
Furthermore, we have the following estimate (by $\omega_d$ we denote the volume of $d$-dimensional
unit ball):
$$
\eqalign{	|B^m_n(f,g)| 	& = |B(M_nf, N_ng)|	\cr
	& \mj 	C_{m+1} \nor{M_n f}{\Lb{K_{m+1}}} \nor{N_n g}{\Lb{L_{m+1}}}\cr
	& \mj	C_{m+1} \nor{\rho_n}{\Ld\Rd}\nor{f}{\Ld{K_{m}}} \nor{\rho'_n}{\Ld{\mR^{d'}}}\nor{g}{\Ld{L_{m}}}\cr
	& \mj	C_{m+1} CC' e^{-2}\sqrt{\omega_d\omega_{d'}} n^{d+d'\over2} \nor{f}{\Ld{K_{m}}}\nor{g}{\Ld{L_{m}}} \;,\cr
}
$$
so for any (large enough) $n$ we have that $B^m_n$ is a continuous bilinear form on $\Ld{K_{m}}\times\Ld{L_{m}}$.
Therefore, it can be represented by an operator $T^m_n\in\lL(\Ld{K_{m}};\Ld{L_{m}})$ in the sense that
$$
\Dup{\overline{T^m_n f}}g = B^m_n(f,g) = B(M_nf, N_ng)\;.
$$
Of course, $T^m_n$ is not uniformly bounded with respect to $n$.

The operators $M_n$ and $N_n$ are Hilbert-Schmidt operators, with corresponding kernels
$m_n$ and $n_n$ (for standard material on Hilbert-Schmidt operators the reader might consult 
[\cfor{JPAub}, Chapter 12]). If we restrict operator $M_n$ to $\Ld{K_m}$, in the sense that we also 
restrict the resulting function to $K_m$, we again obtain a Hilbert-Schmidt operator $M^m_n$,
with kernel $m^m_n$ equal to $m_n$ restricted to $K_m\times K_m$. As $T_n^m$ is a bounded operator, the
composition $T_n^m M^m_n$ is again a Hilbert-Schmidt operator, but from $\Ld{K_m}$ to $\Ld{L_m}$, with a
kernel which we denote by $k^m_n\in\Ld{K_m\times L_m}$.
\medskip

\noindent {\bf III. } {\sl Positivity and boundedness of approximative Hilbert-Schmidt kernels}
\smallskip

Since $M_n$ and $N_n$ are positive operators, for non-negative functions $f$ and $g$ we have that
$$
\Dup{\overline{T^m_n f}}g = B(M_n f, N_n g) \vj0 \;,
$$
so $T^m_n$ is a positive operator, thus $T^m_n M^m_n$ as well, which implies that its kernel
$k^m_n \vj0$.

On the other hand, for constants 1 on compacts $K_m$ and $L_m$ we have:
$$
\eqalign{	0\mj \Dup{T^m_n  M^m_n \chi_{K_m}}{\chi_{L_m}}	& = B(M_n M^m_n \chi_{K_m}, N_n\chi_{L_m}) \cr
	& = |B(M_n M^m_n \chi_{K_m}, N_n \chi_{L_m})| 	\cr
	& \mj C_{m+1} \nor{M_n M^m_n \chi_{K_m}}{\Lb X}
	\nor{N_n\chi_{L_m}}{\Lb Y} \mj C_{m+1} \;. \cr
}
$$
Note that $M^m_n \chi_{K_m}$ is identical to 1 on $K_{m-1}$, while supported in $K_m$.
Further,  $M_n M^m_n \chi_{K_m}$ is identical to 1 on $K_{m-2}$ (here we need to take $m\vj 4$), 
and supported in $K_{m+1}$. For $N_n\chi_{L_m}$ we similarly have that it is equal to 1 on 
$L_{m-1}$, while it is supported in $L_{m+1}$. Therefore, the ${\rm L}^\infty$ norms in the 
above inequality are both equal to 1.

In particular, from the above we can conclude that $\nor{k^m_n}{\Lj{K_m\times L_m}} \mj C_{m+1}$.
It might be of interest to notice that positivity was crucial in obtaining the boundedness.
\medskip

\noindent {\bf IV. } {\sl Passage to the limit}
\smallskip

By the continuity of bilinear form, for $m$ large enough we get
$$
B(f,g) = \lim_n B(M_n M^m_n f , N_n g) =  \lim_n \Dup{\overline{T^m_n M^m_n f}}g 
	= \lim_n \Dup{k^m_n}{f\kvtpr g} \;.
\leqno(\rfor{BzaT})
$$
Indeed, take $f\in\Cc X$ and $g\in\Cc Y$. Naturally, for some large enough $m$, $\supp f\subseteq K_{m-1}$
and $\supp g\subseteq L_{m-1}$, so we have $f\in\pCd{K_{m-1}}X$ and $g\in\pCd{L_{m-1}}Y$.
By standard properties of mollifiers, the latter immediately implies that $N_n g\str g$ in the space $\pCd{L_{m}}Y$.
For the former, we first notice that $M^m_nf\str f$ in $\pCd{K_m}X$, thus $h_n:=M^m_n f - f\str0$ in
the same space. Therefore
$$
|(M_n h_n)(\mx)| \mj \int_{K_m}\rho_n(\mx-\my) |h_n(\my)| d\my < \eps\int_{K_m} \rho_n(\mx-\my) d\my \mj \eps \;,
$$
for $n\vj n_0$ such that $\sup_{\mx\in X}|h_n(\mx)| < \eps$, which furnishes the proof of (\cfor{BzaT}).

As the sequence $(k^m_n)_{n\in\mN}$ is bounded in $\Lj{K_{m-1}\times L_{m-1}}$ by $C_{m+1}$
(actually, here we consider the restrictions of $k^m_n$ to $K_{m-1}\times L_{m-1}$, without
unnecessarily reflecting that in the notation; this restriction is needed in order to use (\cfor{BzaT})),
it is weakly $\ast$ precompact in the space of bounded Radon measures.
Since any accumulation point satisfies (\cfor{BzaT}), they all coincide on tensor products, so by density
also on $\Cp{K_{m-1}\times L_{m-1}}$. Therefore the accumulation point is unique, and the whole sequence
$k^m_n$ converges to a non-negative Radon measure; let us denote it by $\mu^m$. 

Clearly, by (\cfor{BzaT}) again, we have that $\mu^{m+1}$ and $\mu^m$ coincide on
$\Cp{K_{m-1}\times L_{m-1}}$, so for any $\Phi\in\Cc{X\times Y}$ we can define $\mu$
by
$$
\Dup\mu\Phi := \Dup{\mu^m}\Phi	\;,
$$
for any $m$ such that $\supp\Phi\subseteq K_{m-1}\times L_{m-1}$, and the definition is good,
thus obtaining a non-negative unbounded Radon measure on $X\times Y$.
This measure $\mu$ is indeed a representation of bilinear functional $B$.
\medskip

Furthermore, if $B$ is a non-negative continuous bilinear form on $\Cnl X\times\Cnl Y$, then 
the sequence of kernels $(k_n^m)$ is uniformly bounded in the space $\Lj{K_m\times L_m}$
both in $n$ and $m$ (as in this case one constant $C_0$ is good for any $m$ in (\cfor{Bcont})), 
thus $(\mu^m)$ is also uniformly bounded in $m$ which shows that $\mu$ is a bounded Radon 
measure.  
}

In fact, we want to apply the above lemma to a bilinear form on $\Cc\Omega\times\Cp\kob$.
To this end, we can consider both $\Omega$ and $\kob$ as (nice) topological manifolds
(for standard material on manifolds the reader might consult [\ref{Lang}, Chapter XXII]),
which can be covered by a countable atlas (even finite for the compact $\kob$), with the
image of charts being the unit balls in the Euclidean space (for the $\kob$, which is a manifold with
boundary, some charts have unit semiballs ${\rm K}^+(\vnul, 1):=\{\mx\in\Ko\vnul1: x^d\vj0\}$ as images).
The coordinate patches on both $\Omega$ and $\kob$ can be chosen in such a way that any compact
set contained in the manifold intersects only a finite number of patches.
In such an open cover we can inscribe a partition of unity, thus reducing the problem
to the local one, considered in the lemma.

However, for the manifold with boundary $\kob$, some coordinate images might be semiballs 
${\rm K}^+(\vnul, 1)$, for which we first introduce an extended bilinear form $\tilde B$ on 
$\Cc{\Ko\vnul 1}\times\Cc{\Ko\vnul 1}$ by $\tilde B(f,g) := B(f,g\rest{{\rm K}^+(\vnul,1)})$, 
and then apply Lemma 4 obtaining measure $\tilde \mu$ on $\Ko\vnul 1\times\Ko\vnul 1$.
Then we define 
$$
\Dup\mu\Phi := \Dup{\tilde\mu}{\tilde\Phi}\;,
$$
where $\tilde\Phi$ is an extension of $\Phi$ by reflection to the lower semiball in the second argument.

The above discussion, which can easily be expressed in terms of general manifolds $X$ and $Y$
instead of $\Omega$ and $\kob$, we summarise as a lemma.

\lm{
Let $X$ and $Y$ be two Hausdorff second countable topological manifolds (with or without a boundary),
and let $B$ be a non-negative continuous bilinear form on $\Cc X\times \Cc Y$.
Then the conclusions of Lemma 4 hold.
}

\podnaslov{The existence}

Similarly as we have done in the proof of Lemma 4 (for open unit ball $X$),
let us denote by $(K_m)$ a sequence of compacts in $\Omega$ which exhaust
$\Omega$; more precisely, such that $K_m \subseteq \Int K_{m+1}$ and $\bigcup_m K_m = \Omega$,
and recall that $\Cbc\Omega$ is a strict inductive limit of spaces ${\rm C}_{K_m}(\Omega)$.
Moreover, note that each ${\rm C}_{K_m}(\Omega)$ (with appropriate restriction to $K_m$) is a 
closed subspace of separable Banach space $\Cp{K_m}$, therefore separable itself.

\tmd{ {\bf (existence of one-scale H-measures)} \ 
If \/ $\vu_n\dscon\vnul$ in $\Ldl{\Omega;\Cr}$ and
$\omega_n\to 0^+$, then there exists a subsequence\/ $(\vu_{n'})$ 
and an $r \times r$ hermitian non-negative matrix Radon measure\/ 
$\jhmv^{(\omega_{n'})}$ on\/ $\Omega\times\kob$ such that for any 
$\ph_1,\ph_2\in\Cc\Omega$ and $\psi \in \Cp\kob$ one has:
$$
\eqalign{
\lim_{n'} \int\limits_\Rd \Bigl(\widehat{\ph_1\vu_{n'}}(\mxi) \otimes 
	\widehat{\ph_2\vu_{n'}}(\mxi)\Bigr) \psi(\omega_{n'}\mxi) \,d\mxi
	&= \Dupp{\jhmv^{(\omega_{n'})}}{(\ph_1 \bar\ph_2)\kvtpr\psi} \cr
& = \int\limits_{\Omega\times\kob} \ph_1(\mx)\bar\ph_2(\mx)\psi(\mxi) 
	\,d\bar\jhmv^{(\omega_{n'})}(\mx,\mxi) \,. \cr
}
$$
The above measure $\jhmv^{(\omega_{n'})}$ is called\/ {\rm the one-scale 
H-measure (with characteristic length $(\omega_{n'})$)}\/ associated to the 
(sub)sequence $(\vu_{n'})$.
}
{%
For simplicity, we shall omit writing explicitly the characteristic length of
the one-scale H-measure in the proof. 

For any test functions $\ph_1$ and $\ph_2$ supported in a $K_m$, and $\psi\in\Cp\kob$,
by the Plancherel formula the following sequence of integrals is bounded:
$$
\eqalign{
\biggl|\, \int\limits_\Rd \Bigl(\widehat{\ph_1\vu_{n}}(\mxi) \otimes 
		\widehat{\ph_2\vu_{n}} & (\mxi)\Bigr)\psi(\omega_{n}\mxi) \,d\mxi \biggr| \cr
&\mj \Bigl(   \limsup_n \nor{\vu_n}{\Ld{K_m;\Cr}}^2   \Bigr) 
		\nor{\ph_1}{\Lb{K_m}} \nor{\ph_2}{\Lb{K_m}} \nor{\psi}{\Lb\Rd} \;, \cr
}
\leqno(\rfor{estimA})
$$
hence we can pass to a subsequence converging to a limit satisfying the same bound.

This passage should be done in such a way that we obtain a subsequence $(\vu_{n'})$ 
good for any choice of test functions. 
The space $\Cp\kob$ is a separable Banach space, so we can consider only
a countable dense set of $\psi\in\Cp\kob$, and use the Cantor diagonal
procedure. However, for $\Cc\Omega$ we cannot apply the same idea, even
though it is separable. We have to consider separable Banach spaces of
continuous functions with ${\rm L}^\infty$ norm, and then use the fact that
the original space is a strict inductive limit of such spaces.

Thus for each $m\in\mN$ we consider countable dense subsets of ${\rm C}_{K_m}(\Omega)$ 
and $\Cp\kob$, and choose a subsequence of integrals in (\cfor{estimA})  that
converge to a limit denoted by $\mL(\ph_1,\ph_2,\psi)\in\M rr$, for any
test functions from the chosen countable sets, and then by continuity for any
$\ph_1,\ph_2\in{\rm C}_{K_m}(\Omega)$ and $\psi\in\Cp\kob$.
The details in a similar setting could be found in the proof of [\ref{ALjmaa}, Theorem 6],
so we omit them here.

After choosing a good subsequence for $m$, we pass to a subsequence again
for $m+1$. Applying the Cantor diagonal procedure once more, finally we arrive
at a subsequence, denoted by $n'$, which is good for any choice of
test functions $\ph_1,\ph_2\in\Cc\Omega$ as well.

Furthermore, $\mL(\ph_1,\ph_2,\psi)$ depends only on the product $\ph_1\bar\ph_2$ 
(not on both $\ph_1$ and $\ph_2$ independently). Indeed, according to
Lemma 3 we have:
$(\lA_{\psi_n}B_{\bar\ph_2} - B_{\bar\ph_2}\lA_{\psi_n})B_{\ph_1}\vu_n \str \vnul$ in 
$\Ld{\Rd;\Cr}$, hence by the Plancherel formula we have
$$
\eqalign{
\lim_{n'} \int\limits_\Rd \Bigl(\widehat{\ph_1\vu_{n'}}(\mxi)\otimes 
	\widehat{\ph_2\vu_{n'}}(\mxi)\Bigr) \psi(\omega_{n'}\mxi) \,d\mxi
	&= \lim_{n'} \int\limits_\Rd B_{\bar\ph_2}\lA_{\psi_{n'}}B_{\ph_1}\vu_{n'}
	\otimes \vu_{n'} \, d\mx \cr
&= \lim_{n'} \int\limits_\Rd \lA_{\psi_{n'}}B_{\bar\ph_2}B_{\ph_1}\vu_{n'}
	\otimes \vu_{n'} \, d\mx \cr
&= \lim_{n'} \int\limits_\Rd \Bigl(\widehat{(\ph_1\bar\ph_2\vu_{n'})}
	(\mxi) \otimes \widehat{\vu_{n'}}(\mxi)\Bigr) \psi(\omega_{n'}\mxi)
	 \,d\mxi \,, \cr
}
$$
so we can define $\mpL(\ph_1\bar\ph_2,\psi) := \mL(\ph_1,\ph_2,\psi)$. 
In particular, if we take $\ph_2$ to be 
equal to 1 on  $\supp\ph_1$, by
$$
\mpL(\ph_1,\psi) =\mL(\ph_1,\ph_2,\psi) \,
$$ 
we define an $r\times r$ matrix, its entries being continuous bilinear forms on 
$\Cc\Omega\times\Cp\kob$.

For real functions $\ph$ and $\psi$ we have that $\lL(\ph,\psi)$ is a
hermitian matrix. Indeed, 
$$
\eqalign{
\mpL(\ph,\psi) &= \mpL(\ph\ph_2,\psi) \cr
&= \lim_{n'} \int\limits_\Rd \Bigl(\widehat{\ph\vu_{n'}}(\mxi)
	\otimes \widehat{\ph_2\vu_{n'}}(\mxi)\Bigr) \psi(\omega_{n'}
	\mxi) \,d\mxi \cr
&= \lim_{n'} \int\limits_\Rd \Bigl(\widehat{\ph_2\vu_{n'}}(\mxi)
	\otimes \widehat{\ph\vu_{n'}}(\mxi)\Bigr)^\ast \psi
	(\omega_{n'}\mxi) \,d\mxi \cr
&= \biggl(\lim_{n'} \int\limits_\Rd \Bigl(\widehat{\ph_2\vu_{n'}}(\mxi)
	\otimes \widehat{\ph\vu_{n'}}(\mxi)\Bigr) \psi(\omega_{n'}\mxi) 
	\,d\mxi\biggr)^\ast = \mpL(\ph,\psi)^\ast \,, \cr
}
$$
where $\ph_2$ is again taken to be equal to $1$ on the support of $\ph$.
By the just proved hermitian property we have for $\mlambda_1, \mlambda_2\in\Cr$
and real $\ph\in\Cc\Omega$ and $\psi\in\Cp\kob$ that
$$
\eqalign{
2\Re(\mpL(\ph,\psi)\mlambda_1\cdot\mlambda_2) &= \mpL(\ph,\psi)(\mlambda_1+\mlambda_2)
	\cdot(\mlambda_1+\mlambda_2) - \mpL(\ph,\psi)\mlambda_1\cdot\mlambda_1 - 
	\mpL(\ph,\psi)\mlambda_2\cdot\mlambda_2 \,, \cr
2\Im(\mpL(\ph,\psi)\mlambda_1\cdot\mlambda_2) &= \mpL(\ph,\psi)(\mlambda_1+i\mlambda_2)
	\cdot(\mlambda_1+i\mlambda_2) - \mpL(\ph,\psi)\mlambda_1\cdot\mlambda_1 - 
	\mpL(\ph,\psi)\mlambda_2\cdot\mlambda_2 \,. \cr
}
\leqno(\for)
$$
Moreover, for $\ph,\psi\vj 0$ and $\mlambda\in\Cr$ we have
$$
\mpL(\ph,\psi)\mlambda\cdot\mlambda = \mpL(\sqrt{\ph}\sqrt{\ph},\psi)\mlambda\cdot\mlambda = \lim_{n'}\int\limits_\Rd 
	|\widehat{\sqrt{\ph}\vu_{n'}}(\mxi)\cdot\mlambda|^2\psi(\omega_{n'}\mxi) \,d\mxi \vj 0 \,. 
$$
Therefore, for any $\mlambda\in\Cr$ $\mpL(\cdot,\cdot)\mlambda\cdot\mlambda$ is a non-negative continuous
bilinear form on $\Cc\Omega\times\Cp\kob$. Hence, by Lemma 5 for any 
$\mlambda\in\Cr$ there exists $\nu_\mlambda\in\lM(\Omega\times\kob)$ such that 
$$
\svaki{\ph\in\Cc\Omega}\svaki{\psi\in\Cp\kob} \quad \mpL(\ph,\psi)\mlambda\cdot\mlambda = 
	\Dup{\nu_\mlambda}{\ph\kvtpr\psi} \,.
$$ 

Finally, let us define a matrix Radon measure $\jhmv=[\jhm^{ij}]$ by
$$
\eqalign{
\Re\jhm^{ij} &:= {1\over 2}(\nu_{\ve_j+\ve_i} - \nu_{\ve_j} - \nu_{\ve_i}) \,,\cr
\Im\jhm^{ij} &:= {1\over 2}(\nu_{\ve_j+i\ve_i} - \nu_{\ve_j} - \nu_{\ve_i}) \,,\cr
}
$$ 
where vectors $\ve_i$, $i\in1..r$, form a canonical basis of $\Cr$. By (\pbfor0) we have that $\mpL$ and 
$\jhmv$ coincide on real functions, and by linearity we have the same for arbitrary $\ph\in\Cc\Omega$
and $\psi\in\Cp\kob$. In addition, $\jhmv$ inherits the hermitian and positivity property of $\mpL$.
}

As with H-measures before, when there is no fear of ambiguity, we assume that 
we have already passed to a subsequence determining an one-scale 
H-measure, and reduce the notation to $\jhmv=\jhmv^{(\omega_n)}$. 
When the {\sl whole} sequence admits an one-scale H-measure
with characteristic length $(\omega_n)$ (i.e.~the definition is valid
without passing to a subsequence), we say that the sequence is 
{\sl $(\omega_n)$-pure} (in terms of one-scale H-measures). 
For the sake of simplicity, here we use the same notion as in the case
with semiclassical measures, while in the situations where there will
not be clear which definition is used, additional explanation will be 
given. Let us just emphasise that from Corollary 6 it will be clear 
that the property of a sequence being $(\omega_n)$-pure {\sl in terms of
one-scale H-measures} implies both that the sequence is pure and 
$(\omega_n)$-pure {\sl in terms of semiclassical measures}.

\medskip
\break

\podnaslov{First properties}

An immediate consequence of the definition is the following simple localisation 
property.
 
\krd{
Let $\jhmv$ be a one-scale H-measure with an arbitrary characteristic
length, determined by the sequence $(\vu_{n})$. 
If all components $u_{n}^i=\vu_{n}\cdot\ve_i$ have their supports in
closed sets $K_i\subseteq\Omega$ (respectively), then the support of 
$\jhm^{ij}=\jhmv\ve_j\cdot\ve_i$ is necessary contained in 
$(K_i\cap K_j)\times\kob$. 
}
{%
For $\ph_1\in\Cc\Omega$ such that $\supp\ph_1\subseteq\Omega
\setminus(K_i\cap K_j)$ we take $\ph_2\in\Cc\Omega$  being equal to 1 on 
the support of $\ph_1$ and vanishing on $K_i\cap K_j$. Since both 
$\ph_1 u_{n}^i$ and $\ph_2 u_{n}^j$ are zero, we have
$$
\eqalign{
\Dup{\jhm^{ij}}{\ph_1\kvtpr\psi} =& \Dup{\jhm^{ij}}{\ph_1\bar\ph_2
	\kvtpr\psi} \cr
=& \lim_{n} \int\limits_\Rd \widehat{\ph_1 u_{n}^i}(\mxi)
	\overline{\widehat{\ph_2 u_{n}^j}}(\mxi) \psi(\omega_{n}\mxi) 
	\,d\mxi
	= 0 \,, \cr
}
$$
where $\psi\in\Cp\kob$ is an arbitrary test function and $(\omega_n)$
is the characteristic length of the measure. Thus we get the
claim.
}

The more precise relation between supports of diagonal and non-diagonal
elements is a consequence of the Cauchy-Schwartz-Bunjakovskij inequality.

\krd{
Let $\jhmv$ be a one-scale H-measure with an arbitrary characteristic
length, determined by the sequence $(\vu_{n})$. 
The support of $\jhm^{ij}$ is contained in the intersection of supports
of the corresponding diagonal elements $\jhm^{ii}$ and $\jhm^{jj}$, i.e.
$$
\supp\jhm^{ij} \subseteq \supp\jhm^{ii} \cap \supp\jhm^{jj} \,.
$$
\vskip-2mm
}
{
Let $\ph_1, \ph_2\in\Cc\Omega$ and $\psi\in\Cp\kob$.
By the definition of one-scale H-measures and the 
Cauchy-Schwartz-Bunjakovskij inequality we have
$$
\eqalign{
|\Dup{\jhm^{ij}}{\ph_1\bar\ph_2\kvtpr\psi}| &\mj \limsup_{n} 
	\int\limits_\Rd \Bigl(|\widehat{\ph_1 u_{n}^i}(\mxi)|
	\sqrt{|\psi|}(\omega_{n}\mxi)\Bigr) \Bigl(|\widehat{\ph_2 u_{n}^j}
	(\mxi)|\sqrt{|\psi|}(\omega_{n}\mxi)\Bigr) \,d\mxi \cr
&\mj  \biggl(\lim_{n}\int\limits_\Rd |\widehat{\ph_1 u_{n}^i}(\mxi)|^2
	|\psi|(\omega_{n}\mxi) \,d\mxi \biggr)^{{1 \over 2}}
	\biggl(\lim_{n}\int\limits_\Rd |\widehat{\ph_1 u_{n}^j}(\mxi)|^2
	|\psi|(\omega_{n}\mxi) \,d\mxi \biggr)^{{1 \over 2}} \cr
&= \sqrt{\Dup{\jhm^{ii}}{|\ph_1|^2\kvtpr|\psi|}} 
	\sqrt{\Dup{\jhm^{jj}}{|\ph_2|^2\kvtpr|\psi|}} \,, \cr
}
$$
so $\Dup{\jhm^{ii}}{|\ph_1|^2\kvtpr|\psi|}=0$ or $\Dup{\jhm^{jj}}{|\ph_2|^2
\kvtpr|\psi|}=0$ implies $\Dup{\jhm^{ij}}{\ph_1\bar\ph_2\kvtpr\psi}=0$.
}

One-scale H-measures are in general associated to complex vector valued 
functions. However, if the observed sequence $(\vu_n)$ is real, 
the associated one-scale H-measure has an additional feature.

\krd{
Let $(\vu_n)$ be an $(\omega_n)$-pure sequence in $\Ldl{\Omega;\Cr}$, 
and $\jhmv$ the corresponding one-scale H-measure. 
Then the sequence $(\bar\vu_n)$ is $(\omega_n)$-pure with
associated one-scale H-measure $\mnu_\kobb$, such that 
$\mnu_\kobb(\mx,\mxi) = \jhmv^\top(\mx,-\mxi)$.

In particular, a one-scale H-measure $\jhm$ associated to a real scalar
sequence is antipodally symmetric, i.e.~$\jhm(\mx,\mxi) = \jhm(\mx,-\mxi)$.
}
{
Let us take test functions $\ph_1, \ph_2\in\Cc\Omega$ and 
$\psi\in\Cp\kob$, and calculate the one-scale H-measure associated to 
the sequence $(\bar\vu_n)$.
$$
\eqalign{
\int\limits_\Rd \Bigl(\widehat{\ph_1\bar\vu_n}(\mxi)\otimes
	\widehat{\ph_2\bar\vu_n}(\mxi)\Bigr) \psi(\omega_n\mxi) \,d\mxi
	=& \overline{\int\limits_\Rd \Bigl((\bar\ph_1\vu_n)^\vee(\mxi)
	\otimes (\bar\ph_2\vu_n)^\vee(\mxi)\Bigr) \bar\psi(\omega_n\mxi)}
	\,d\mxi \cr
=& \overline{\int\limits_\Rd \Bigl((\bar\ph_1\vu_n)^\vee(-\meta)
	\otimes (\bar\ph_2\vu_n)^\vee(-\meta)\Bigr) \bar\psi
	(-\omega_n\meta)} \,d\meta \cr
=& \overline{\int\limits_\Rd \Bigl(\widehat{\bar\ph_1\vu_n}(\meta)
	\otimes \widehat{\bar\ph_2\vu_n}(\meta)\Bigr) \tilde{\bar\psi}
	(\omega_n\meta)} \,d\meta \,, \cr
}
$$
where we made use of a change of variable ($\meta=-\mxi$) in the 
second equality, while $\tilde\psi$ denotes the sign change in argument $\mxi$
($\tilde\psi(\mxi) = \psi(-\mxi)$). Hence, the limit of the term on the
left hand side exists and we have  
$$
\eqalign{
\Dup{\mnu_\kobb}{\ph_1\bar\ph_2\kvtpr\psi}
=& \overline{\Dup{\jhmv}{\bar\ph_1\ph_2\kvtpr\tilde{\bar\psi}}} \cr
=& \overline{\Dup{\tilde\jhmv}{\bar\ph_1\ph_2\kvtpr{\bar\psi}}}
	= \Dup{{\tilde\jhmv}^\top}{\ph_1\bar\ph_2\kvtpr\psi} \,. \cr
}
$$
The last equality is a consequence of the hermitian property of one-scale 
H-measures. 
}

If $(\vu_n)$ is a weakly convergent sequence in $\Ldl{\Omega;\Cr}$, then 
$\vu_n\otimes\vu_n$ is bounded in the space $\Ljl{\Omega;\M rr}$, thus on 
a subsequence converges weakly to a hermitian non-negative Radon measure 
$\mnu$ (the defect measure). The relationship between one-scale H-measure 
and the defect measure is given by the following corollary.

\krd{
If $\vu_n\dscon\vnul$ in $\Ldl{\Omega;\Cr}$ is such that 
$\vu_n\otimes\vu_n \povrhsk\ast \mnu$, the limit being
in the space of Radon measures, then for any $\ph\in\Cc\Omega$
$$
\Dup\mnu\ph = \Dup\jhmv{\ph\kvtpr 1} \,,
$$
where $\jhmv$ is a one-scale H-measure with an arbitrary characteristic 
length associated to a (sub)sequence $(\vu_n)$.
}
{%
By inserting $\psi\equiv 1$ in the definition of one-scale H-measures 
and using the Plancherel formula we get the claim. 
}

From the previous corollary it is straightforward to conclude that a trivial 
one-scale H-measure of an arbitrary characteristic length implies 
strong convergence in $\Ldl{\Omega;\Cr}$ of the observed sequence.  
On the other hand, for a strongly convergent sequence the associated 
one-scale H-measure (of an arbitrary characteristic length) is trivial. 
Let us remember that for this property of semiclassical measures we
need an additional assumption that the corresponding sequence is 
$(\omega_n)$-oscillatory. As semiclassical measures associated to such 
sequences  do not lose information at infinity (in the dual space), a 
simple characterisation of that assumption in terms of one-scale 
H-measures easily follows. 

\krd{
If $\vu_n\dscon\vnul$ in $\Ldl{\Omega;\Cr}$ is $(\omega_n)$-pure (in the sense of one-scale H-measures),
then $\tr\jhmv^{(\omega_n)}(\Omega\times \Sigma_\infty)=0$ if and only if $(\vu_n)$ is 
$(\omega_n)$-oscillatory.
}  
{
By definition we have the following equivalence
$$
\tr\jhmv(\Omega\times\Sigma_\infty)=0 \akko
	\svaki{\ph\in\Cc\Omega} \quad \int\limits_{\Omega\times\kob}\!\!|\ph(\mx)|^2
	\chi_{\Sigma_\infty}(\mxi) \,d\tr\jhmv(\mx,\mxi) =0 \,,
$$
where $\chi_{\Sigma_\infty}$ denotes the characteristic function of the sphere at $\infty$  
in $\kob$.

By $\zeta\in\Cbc\Rd$ let us denote a smooth cutoff function such that 
$0\mj\zeta\mj 1$ and such that it is identically equal to $1$ on $\Kz{\vnul}{1}$, while
$\supp\zeta\subseteq\Ko\vnul 2$. Further, we define $\zeta_m:=\zeta({\cdot\over m})$.

By the Lebesgue dominated convergence theorem, non-negativity of diagonal elements and the 
definition of one-scale H-measures for an arbitrary $i\in 1..d$ we have
$$
\eqalign{
0=\int\limits_{\Omega\times\kob}|\ph(\mx)|^2\chi_{\Sigma_\infty}(\mxi) \,d\jhm^{ii}(\mx,\mxi) =& 
	\lim_m\int\limits_{\Omega\times\kob}|\ph(\mx)|^2(1-\zeta_m(\mxi)) \,d\jhm^{ii}(\mx,\mxi) \cr
=& \lim_m\lim_n \int\limits_\Rd |\widehat{\ph u_n^i}(\mxi)|^2 
	(1-\zeta_m(\omega_n\mxi)) \,d\mxi \cr
\vj& \limsup_m\limsup_n \int\limits_{|\omega_n\mxi|\vj 2m} 
	|\widehat{\ph u_n^i}(\mxi)|^2 \,d\mxi \vj 0 \,, \cr
}   
$$
where in the last step we have used the fact that $\supp\zeta_m\subseteq\Ko{\vnul}{2m}$.
Hence 
$$
\svaki{i\in 1..d} \qquad \lim_m\limsup_n \int\limits_{|\mxi|\vj{2m\over\omega_n}} 
	|\widehat{\ph u_n^i}(\mxi)|^2 \,d\mxi = 0 \,,
$$
and therefore $(\vu_n)$ is $(\omega_n)$-oscillatory. 

The opposite implication follows from the estimate
$$
\svaki{i\in 1..d} \qquad \int\limits_\Rd |\widehat{\ph u_n^i}(\mxi)|^2 
	(1-\zeta_m(\omega_n\mxi))\,d\mxi 
	\mj \int\limits_{|\mxi|\vj {m\over\omega_n}} 
	|\widehat{\ph u_n^i}(\mxi)|^2 \,d\mxi \,.
$$
\vskip-5mm
}

We shall illustrate further advantages of one-scale H-measures by 
reconsidering sequences from  examples 1 and 2. 

\Pr{%
The one-scale H-measure with characteristic length 
$(\omega_n)$  associated to the sequence $(u_n)$ from Example 1 is given by
$$
\jhm^{(\omega_n)} = \lambda\kvtpr 
	\left\{
	\matrix{ 
	\delta_{\vnul^{{\vk\over|\vk|}}} & , & \lim_n n^\alpha\omega_n 
		= 0 \cr
	\delta_{c\vk} & , & \lim_n n^\alpha\omega_n = c\in\oi0\infty \cr
	\delta_{\infty^{{\vk\over|\vk|}}} & , & \lim_n n^\alpha\omega_n
		 = \infty \cr
	} \right. \,.
$$ 
Contrary to the semiclassical measure, one-scale H-measure contains
information on the direction of oscillation for every characteristic
length $(\omega_n)$.

For the sequence $(u_n+v_n)$ observed in Example 2, the associated one-scale 
H-measure with characteristic length $(\omega_n)$ is given by
$$
\jhm^{(\omega_n)} = \lambda\kvtpr 
	\left\{ \matrix{ 
	(\delta_{\vnul^{{\vk\over|\vk|}}}+\delta_{\vnul^{{\vs\over|\vs|}}})
		& , & \lim_n n^\beta\omega_n = 0 \cr
	(\delta_{\vnul^{{\vk\over|\vk|}}}+\delta_{c\vs}) & , & \lim_n 
		n^\beta\omega_n = c\in\oi{0}{\infty} \cr
	(\delta_{\vnul^{{\vk\over|\vk|}}} + \delta_{\infty^{{\vs\over 
		|\vs|}}}) & , & \lim_n n^\beta\omega_n = \infty \ \&	 
		\ \lim_n n^\alpha\omega_n = 0 \cr
	(\delta_{c\vk}+\delta_{\infty^{{\vs\over|\vs|}}}) & , & \lim_n 
		n^\alpha\omega_n = c\in\oi{0}{\infty} \cr
	(\delta_{\infty^{{\vk\over|\vk|}}}+\delta_{\infty^{{\vs\over 
		|\vs|}}}) & , & \lim_n n^\alpha\omega_n = \infty \cr
	} \right. \,.
$$

It is important to notice that in both examples the respective H-measure can be 
derived from the one-scale H-measure with arbitrary characteristic
length, which was not the case with the semiclassical measures.  
}

By Lemma 1 and the known results for semiclassical and one-scale H-measures, we have that the 
semiclassical measure $\pkmv$ and the one-scale H-measure $\jhmv$ associated 
to the same (sub)sequence and with the same characteristic length, coincide 
on $\Cbc\Omega\times\lS(\Rd)$ under necessary assumptions that the weak limit of the observed 
(sub)sequence is zero. Of course, by continuity we can extend this
conclusion to $\Cc\Omega\times\Cnl\Rd$.

Analogously, by using Lemma 1 and the definitions of H-measures and one-scale 
H-measures, we have:
$$
\eqalign{
\Dup{\jhmv}{(\ph_1\bar\ph_2)\kvtpr(\psi\circ\mpi)} &= 
	\lim_{n}\int\limits_\Rd \Bigl(\widehat{\ph_1\vu_{n}}(\mxi)\otimes
	\widehat{\ph_2\vu_{n}}(\mxi)\Bigr)	\psi\Bigl(\mpi(\omega_{n}\mxi)
	\Bigr) \,d\mxi	\cr
&= \lim_{n} \int\limits_\Rd \Bigl(\widehat{\ph_1\vu_{n}}(\mxi)\otimes
	\widehat{\ph_2\vu_{n}}(\mxi)\Bigr)	\psi\biggl({\mxi \over |\mxi|}
	\biggr) \,d\mxi 
	= \Dup{\hmv}{(\ph_1\bar\ph_2)\kvtpr\psi}  \;, \cr
}
$$
where $\ph_1,\ph_2\in\Cc\Omega$ and $\psi\in\Cp\Sdmj$ are arbitrary 
test functions, while the H-measure $\hmv$ and the one-scale H-measure 
$\jhmv$ are associated to the same (sub)sequence.
In particular, this relationship justifies the name {\sl one-scale 
H-measure}, as an H-measure with one characteristic length.
Let us summarise the previous observations.

\kr{
Let $\omega_n\str0^+$ and let $\hmv$, $\pkmv^{(\omega_n)}$ and $\jhmv^{(\omega_n)}$ be the H-measure,
the semiclassical and the one-scale H-measure (respectively), associated to the same
(sub)sequence, weakly converging to zero. 
Then for any choice of test functions $\ph_1,\ph_2\in\Cc\Rd$ one has

\item{$i)$} $\svaki{\psi\in\Cp\Sdmj}$ \ 
$\Dup{\jhmv^{(\omega_n)}}{(\ph_1\bar\ph_2)\kvtpr (\psi\circ\mpi)} = 
	\Dup{\hmv}{(\ph_1\bar\ph_2)\kvtpr \psi}$,

\item{$ii)$} $\svaki{\psi\in\Cnl{\Rd}}$ \ 
$\Dup{\jhmv^{(\omega_n)}}{(\ph_1\bar\ph_2)\kvtpr \psi} = 
	\Dup{\pkmv^{(\omega_n)}}{(\ph_1\bar\ph_2)\kvtpr \psi}$.
}
\smallskip

\Napbr{%
Although in the case when $\vu_n\dscon\vu\not=\vnul$ in 
$\Ldl{\Omega;\Cq}$ we cannot associate one-scale H-measure to 
$(\vu_n)$, the limit of (sub)sequence of integrals in Theorem 6
can be explicitly calculated. Indeed, inserting 
$\vu_n=(\vu_n-\vu)+\vu$, by Theorem 6 and the fact that the product 
of a weakly and a strongly converging sequence converges in $\mC$, we get
$$
\lim_{n'} \int\limits_\Rd \Bigl(\widehat{\ph_1\vu_{n'}}(\mxi) \otimes 
	\widehat{\ph_2\vu_{n'}}(\mxi)\Bigr) \psi(\omega_{n'}\mxi) \,d\mxi
	= \Dup{\mnu_{\kobb}}{\ph_1\bar\ph_2\kvtpr\psi} +
	\int\limits_\Rd\Bigl(\widehat{\ph_1\vu}(\mxi) \otimes 
	\widehat{\ph_2\vu}(\mxi)\Bigr) (\psi_0\circ\mpi)(\mxi) \,d\mxi
	\,,
$$
where $\mnu_\kobb$ is the one-scale H-measure associated to 
$(\vu_n-\vu)$, and $\psi_0\in\Cp{\Sdmj}$ is defined by (\cfor{psinib}). 
Moreover, if we take $\psi\in\Cnl{\Rd}$ (hence 
$\psi_0\equiv \psi(\vnul)$) and apply the Plancherel
formula, by Corollary 6 the above resembles (\cfor{decompPkm}).
Therefore, in the case when the limit of the observed sequence
is known, we can reconstruct the semiclassical measure starting 
from the one-scale H-measure. 
}


\odjeljak{Localisation principle}

Our next goal is to investigate the properties of one-scale H-measures determined by 
solutions to linear differential relations with variable coefficients.
More precisely, we consider a sequence of systems of partial differential
equations of order $m\in\N$ with a characteristic length
$$
\sum_{l\mj|\malpha|\mj m}\eps_{n}^{|\malpha|-l}\partial_\malpha
	(\mA^{\malpha}\vu_{n}) = \vf_n \,,
\leqno(\rfor{locprEqJhmA})
$$
where $(\eps_n)$ is a bounded sequence
of positive numbers, and $l\in0..m$. Our aim is to develop a localisation principle 
for the corresponding one-scale H-measure, and eventually to deduce the conditions 
implying strong convergence of $(\vu_n)$.

\podnaslov{A notable condition}

To this effect we have to impose appropriate conditions on the terms entering (\cfor{locprEqJhmA}). 
The characteristic length of the problem, expressed by the sequence $(\eps_n)$,
naturally appears in the assumption on the convergence of the right hand side.
Namely, for a sequence $(\vf_n)$ in $\H{-m}{\Rd;\Cq}$, the following
convergence turns out to be important
$$
{\widehat\vf_n\over 1+k_n} \str \vnul 
	\quad \hbox{in} \quad  \Ld{\Rd;\Cq} \,,
\leqno(\rfor{ConvRhs}) 
$$ 
where $k_n(\mxi):=\sum_{s=l}^m\eps_n^{s-l}|\mxi|^s$.

Actually, the local version will turn out to be more natural: if the sequence $(\vf_n)$ 
belongs only to $\Hl{-m}{\domena;\Cq}$, we assume the local version
of the above condition
$$
\svaki{\ph\in\Cbc{\domena}} \qquad {\widehat{\varphi \vf_n}\over 1+k_n} \str \vnul  
	\quad \hbox{in} \quad  \Ld{\Rd;\Cq} \,.
\leqno(\rfor{ConvRhsLoc})
$$ 

The characterisation of the above convergence is given by the following lemma.

\lmd{
Let $(\eps_n)$ be a sequence of positive numbers. Then (\cfor{ConvRhs}) is equivalent to 
$$
{\widehat\vf_n\over 1+h_n} 
	\str \vnul \quad \hbox{in} \quad  \Ld{\Rd;\Cq} \,,
$$ 
where $h_n(\mxi):=|\mxi|^l+\eps_n^{m-l}|\mxi|^m$. 

Furthermore, if we additionally assume that there are $\eps_0,\eps_\infty>0$ such that
$\eps_0\mj\eps_n \mj\eps_\infty$, then (\cfor{ConvRhs}) is equivalent to the 
convergence in  $\H{-m}{\Rd;\Cq}$.
}
{%
In order to prove the first statement, it is enough to show that
${1+h_n\over 1+k_n}$ is bounded on $\Rd$ by positive constants both from below and above. 
Clearly, it is bounded above by 1 for $l\ne m$, while for $l=m$ it is bounded by 2. From
$$
\max\{1,\eps_n|\mxi|,\dots,(\eps_n|\mxi|)^{m-l}\} = \max\{1,(\eps_n|\mxi|)^{m-l}\} \,,
$$
it follows that $1+k_n(\mxi)=1+|\mxi|^l\sum_{s=l}^m(\eps_n|\mxi|)^{s-l}\mj 1+(m-l+1) h_n(\mxi)$, 
which implies the lower bound.

In order to deduce the second statement, it is enough to observe the following elementary estimate 
$$
{\min\{1,\eps_\infty^{l-m}\}\over 2}{1\over 1+|\mxi|^m} 
\mj {1\over 1+|\mxi|^l+\eps_n^{m-l}|\mxi|^m} 
	\mj \max\{1,\eps_0^{l-m}\}{1\over 1+|\mxi|^m} \,,
$$
and to apply the sandwich theorem.
}

From the proof of the last lemma we conclude that if the sequence  $(\eps_n)$ is bounded from above,
the condition (\cfor{ConvRhs}) implies strong convergence to $\vnul$ in $\H{-m}{\Rd;\Cq}$.
On the other hand, the strong convergence  in $\H{-l}{\Rd;\Cq}$ is obviously a stronger notion 
than (\cfor{ConvRhs}), implying that the latter is an intermediate convergence between those in
${\rm H}^{-l}$ and ${\rm H}^{-m}$ in this case.

Further relations among the above convergences are provided under additional assumption that $(\vf_n)$
is an $(\eps_n)$-oscillating sequence. More precisely, the following result holds.

\tmd{
Let $(\eps_n)$ be a sequence in $\Rpl$, and $(\vf_n)$ a sequence of vector functions in
$\H{-l}{\Rd;\Cq}$ satisfying (\cfor{ConvRhs}) (with $m\geq l$) and 
$$
\lim_{r\to\infty}\limsup_{n} \int\limits_{|\mxi|\vj{r\over\eps_n}} 
	\Bigl|{\widehat{\vf}_n\over 1+|\mxi|^l}\Bigr|^2 \,d\mxi = 0 \,.
\leqno(\rfor{GenOmegaTit})
$$\vskip-2mm\noindent
Then $\vf_n\str \vnul$ in $\H{-l}{\Rd;\Cq}$. 
}
{
Let us split $\|\vf_n\|_{\H{-l}{\Rd;\Cq}}^2$ into two integrals:
$$
\int\limits_{\Rd} 
	\Bigl|{\widehat{\vf}_n\over 1+|\mxi|^l}\Bigr|^2 \,d\mxi =
	\int\limits_{|\mxi|\vj{r\over\eps_n}}\Bigl|{\widehat{\vf}_n\over 
	1+|\mxi|^l}\Bigr|^2 \,d\mxi
	+ \int\limits_{|\mxi|<{r\over\eps_n}} \Bigl|{\widehat{\vf}_n\over 
	1+|\mxi|^l}\Bigr|^2 \,d\mxi \,.
$$
By (\pbfor0) the first term on the right hand side is arbitrary small for $r$ large enough
(uniformly with respect to large enough $n$).  For the second term we have
$$
\eqalign{
\int\limits_{|\mxi|<{r\over\eps_n}} \Bigl|{\widehat{\vf}_n\over 
	1+|\mxi|^l}\Bigr|^2 \,d\mxi &\mj \int\limits_{|\mxi|<{r\over\eps_n}} 
	\Bigl|{\widehat{\vf}_n\over 1+|\mxi|^l+\eps_n^{m-l}|\mxi|^m}\Bigr|^2 
	\Bigl(1+{|\eps_n\mxi|^m\over\eps_n^l+|\eps_n\mxi|^l}\Bigr)^2\,d\mxi \cr
&\mj (1+r^{m-l})^2 \int\limits_{\Rd} 
	\Bigl|{\widehat{\vf}_n\over 1+|\mxi|^l+\eps_n^{m-l}|\mxi|^m}\Bigr|^2 \,d\mxi \,, \cr
}
$$
which by the assumptions goes to  zero for every $r$ as  $n\to\infty$.  
}

\Napbr{%
Let us note that for $l=0$ relation (\pbfor0) resembles the assumption  of the Riesz-Kolmogorov
compactness theorem [\ref{Bfa}, Section 4.5],  whereas the latter one is
stronger with no $(\eps_n)$ scale included. 
The localised variant of condition (\pbfor0) (with  $\vf_n$ replaced by $\ph\vf_n$ for an arbitrary
$\ph\in\Cbc\domena$) is just a weaker version of the $(\eps_n)$-oscillatory property adapted to 
${\rm H}^{-l}$ setting. 
}

In general, we shall deal with ${\rm H}^{-k}$ precompact sequences, where $k$ is an arbitrary integer
between $l$ and $m$, for which the following result holds.

\lmd{
Let $(\eps_n)$ be a sequence in $\Rpl$ bounded from above, and $(\vf_n)$ a sequence of vector functions
for some $k\in l..m$ strongly converging to $\vf$ in the space $\H{-k}{\Rd;\Cq}$.
Then the sequence $(\eps_n^{k-l}(\vf_n-\vf))$ satisfies condition (\cfor{ConvRhs}).

If we additionally assume that $\eps_n$ converges to zero, then 
$$
{\eps_n^{k-l}\widehat{\vf_n}\over 1+k_n}\quad {\buildrel {\rm L}^2 \over\Dstr }
\quad 	\left\{
	\matrix{
	{\hat\vf \over 1+|\mxi|^l} ,& k=l \cr
	\vnul ,& k\vj l+1 \cr
	}\right. \,.
\leqno(\for)
$$
\vskip-2mm 
}
{
Let $\eps_\infty>0$ be an upper bound for $\eps_n$; from the estimate
$$
{\eps_n^{k-l} |(\widehat\vf_n-\widehat\vf)(\mxi)|\over 1+k_n(\mxi)} \mj 
	{\eps_n^{k-l}|(\widehat\vf_n-\widehat\vf)(\mxi)|\over 1+\eps_n^{k-l}|\mxi|^k} 
	\mj C{|(\widehat\vf_n-\widehat\vf)(\mxi)|\over 1+|\mxi|^k} \,,
$$
where $C=\max\{1,\eps_\infty^{k-l}\}$, we get the claim. 

In order to prove the second statement, one rewrites the left hand side term of (\pbfor0) by adding 
and  subtracting $\vf$ from $\vf_n$. The term containing the difference $\vf_n-\vf$ converges 
strongly to zero by the first part of the lemma, while the limit of the remaining one is 
obtained by the Lebesgue dominated convergence theorem. 
}

If $(\vf_n)$ satisfies (\cfor{ConvRhsLoc}), an application of the above results 
to $\ph\vf_n$, for an arbitrary $\ph\in\Cbc\domena$, gives us that all the conclusions 
obtained for convergence (\cfor{ConvRhs}) are also valid for its local counterpart 
(\cfor{ConvRhsLoc}), the only difference being the spaces $\H{-k}{\Rd;\Cq}$ 
replaced by their local versions $\Hl{-k}{\Omega;\Cq}$.

\podnaslov{The case $\omega_n=\eps_n$}

In the special case $l=1$, Luc Tartar [\ref{Tgth}, Lemma 32.7] proved the 
localisation principle for one-scale H-measures with characteristic length 
equal to the characteristic length of the observed system (\cfor{locprEqJhmA}). 
However, this result does not provide  any information on the structure of one-scale 
H-measures on $\Sigma_0$. In the following theorems we extend that result to more 
general systems, simultaneously  overcoming the mentioned shortcoming.  

First we present results for one-scale H-measures with characteristic length 
$(\omega_n)$ equal to the one of the system under consideration (like in Tartar's 
approach), and later we generalise them to the case of an arbitrary characteristic length.

\tmd{
Let $\vu_n\dscon\vnul$ in $\Ll2{\domena;\Cr}$ and satisfy (\cfor{locprEqJhmA}),
where $\mA^{\malpha}\in\Cp{\domena;\Mx{q\times r}{\mC}}$, $\eps_n\to 0^+$, while $(\vf_n)$ is a 
sequence in  $\Hl{-m}{\domena;\Cq}$ satisfying (\cfor{ConvRhsLoc}).

Then the associated one-scale H-measure $\jhmv$ with characteristic length $(\eps_n)$
satisfies:
$$
\mmp_1\jhmv^\top = \mnul \,,
\leqno(\rfor{locprJhm})
$$
where
$$
\mmp_1(\mx,\mxi):= \sum_{l\mj|\malpha|\mj m}(2\pi i)^{|\malpha|} {\mxi^{\malpha}\over 
	|\mxi|^l+|\mxi|^m}\mA^\malpha(\mx) \,.
\leqno(\rfor{SymbJhmA})
$$
\vskip-2mm 
}
{
{\noindent\bf I.}
In the first step we localise (\cfor{locprEqJhmA}) by
multiplying it by a test function $\ph\in\Cbc\domena$. With the aid of Lemma 8 
below we get
$$
\sum_{l\mj|\malpha|\mj m}\sum_{\vnul\mj\mbeta\mj\malpha}(-1)^{|\mbeta|}
	\binom{\malpha}{\mbeta}\eps_n^{|\malpha|-l}
	\partial_{\malpha-\mbeta}\Bigl((\partial_\mbeta\ph)\mA^\malpha\vu_n \Bigr) 
	= \ph\vf_n \,.
$$
Our goal is to show that the terms with $\mbeta\ne\vnul$ on the left hand side of this
equality satisfy the convergence condition (\cfor{ConvRhs}), i.e.~that we can
subsume them under the right hand side.

To this end, note that for given $\malpha$ and $\mbeta$ the sequence 
$\partial_{\malpha-\mbeta}\bigl((\partial_\mbeta\ph)\mA^\malpha\vu_n \bigr)$ is supported in
a fixed compact ($\supp \ph$), hence by the Rellich compactness theorem for $\mbeta\ne\vnul$
it is strongly precompact in $\H{-|\malpha|}{\Rd;\Cq}$. Thus by Lemma 7  we get that 
$\eps_n^{|\malpha|-l}\partial_{\malpha-\mbeta}\bigl((\partial_\mbeta\ph) \mA^\malpha\vu_n \bigr)$
satisfies (\cfor{ConvRhs}).

Therefore the above equality can be simplified into the form
$$
\sum_{l\mj|\malpha|\mj m}\eps_{n}^{|\malpha|-l} \partial_\malpha(\mA^\malpha\ph\vu_n) 
	= \tilde\vf_n \,,
\leqno(\rfor{locprEqLoc})
$$
where sequence $(\tilde\vf_n)$ satisfies (\cfor{ConvRhs}).
\smallskip

{\noindent\bf II.}
Next we  express the convergence in the terms of convergence in $\Ld{\Rd;\Cq}$:
by applying the Fourier transform to (\cfor{locprEqLoc}), and multiplying the whole equality by 
${1\over 1+|\mxi|^l+\eps_n^{m-l}|\mxi|^m}$ we get
$$
\sum_{l\mj|\malpha|\mj m}\eps_n^{|\malpha|-l}(2\pi i)^{|\malpha|} 
	{\mxi^\malpha\fp{\mA^\malpha\ph\vu_n}\over 1+|\mxi|^l+\eps_n^{m-l}|\mxi|^m}
	= {\fp{\tilde\vf_n}\over 1+|\mxi|^l+\eps_n^{m-l}|\mxi|^m} \,.
$$
By the assumption, after applying Lemma 6, the right hand side tends to $\vnul$ strongly in
$\Ld{\Rd;\Cq}$, so the same holds for the left hand side. Moreover, by Lemma 9 below
we also obtain that
$$
\sum_{l\mj|\malpha|\mj m}(2\pi i)^{|\malpha|}{\eps_n^{|\malpha|-l}\mxi^\malpha\over |\mxi|^l
	+\eps_n^{m-l}|\mxi|^m} \fp{\mA^\malpha\ph\vu_n} \str \vnul \quad \hbox{in} 
	\quad  \Ld{\Rd;\Cq} \,.
$$
Indeed, for $|\mxi|\vj r$ we have $|\mxi|^l+\eps_n^{m-l}|\mxi|^m\vj r^l$, while
$$
\left|{\eps_n^{|\malpha|-l}\mxi^\malpha\over |\mxi|^l + \eps_n^{m-l}|\mxi|^m}\right| 
	\mj {|\eps_n\mxi|^{|\malpha|-l}\over 1+|\eps_n\mxi|^{m-l}} \mj 1 \,,
$$
which justifies the application of Lemma 9.
\smallskip

{\noindent\bf III.} 
In order to introduce one-scale H-measures, we  rewrite the last sum as
$$
\sum_{l\mj|\malpha|\mj m}(2\pi i)^{|\malpha|}{(\eps_n\mxi)^\malpha\over |\eps_n\mxi|^l
	+|\eps_n\mxi|^m} \fp{\mA^\malpha\ph\vu_n} \str \vnul \quad \hbox{in} 
	\quad  \Ld{\Rd;\Cq} \,,
$$
where the function $\mxi\mapsto {\mxi^\malpha \over |\mxi|^l + |\mxi|^m}$ belongs to
$\Cp\kob$ (see Example 3).

After multiplying the above sum by $\psi(\eps_n\cdot)$, where $\psi\in\Cp{\kob}$, and
forming a tensor product with $\fp{\ph_1\vu_n}$, for $\ph_1\in\Cbc\domena$, the
definition of one-scale H-measures implies
$$
\eqalign{
\mnul =& \lim_n \int_\Rd\psi(\eps_n\mxi) \Biggl(\sum_{l\mj|\malpha|\mj m}(2\pi i)^{|\malpha|}
	{(\eps_n\mxi)^\malpha \over |\eps_n\mxi|^l+|\eps_n\mxi|^m} 
	\fp{\mA^{\malpha}\ph\vu_n}\Biggr) \otimes \Bigl(\fp{\ph_1\vu_n}\Bigr) \,d\mxi \;\cr
=& \Dupp{\overline{\sum_{l\mj|\malpha|\mj m}(2\pi i)^{|\malpha|} {\mxi^\malpha 
	\over |\mxi|^l+|\mxi|^m} \mA^{\malpha}}\jhmv}{\ph\bar{\ph_1}\kvtpr\psi} \,. \cr
}
$$
In particular, by choosing $\ph_1$ to be 1 on the support of $\ph$ and using the hermitian 
property of $\jhmv$, we get the claim.
}

\lm{
For any $\malpha\in\mN^d_0$,  $\ph\in\Cb{\domena}$ and a distribution $S\in\lD'(\domena)$ one has:
$$
\ph\partial_\malpha S=\sum_{\vnul\mj\mbeta\mj\malpha}(-1)^{|\mbeta|}\binom{\malpha}{\mbeta}
	\partial_{\malpha-\mbeta}\Bigl((\partial_\mbeta\ph)S\Bigr) \,.
$$
\vskip-2mm
}

The proof of this lemma follows directly from the generalised Leibniz rule.
\smallskip

Next lemma is a generalisation of [\ref{ALrwa}, Lemma 3].

\lmd{
Let $(\vf_n)$ be a sequence of measurable vector functions on $\Rd$ and $(h_n)$ a sequence of
positive scalar functions which are uniformly bounded from below outside an arbitrary neighbourhood
of the origin, i.e.~which satisfy 
$$
\svaki{r>0}\postoji{\tilde C>0}\svaki{n\in\N}\svaki{\mxi\in\Rd\setminus\Ko\vnul r}
	\qquad h_n(\mxi)\vj \tilde C \,.
$$
Furthermore, let $(\vu_n)$ be a bounded sequence in $\Ld{\Rd;\Cr}\cap\Lj{\Rd;\Cr}$, such that
$$
{\vf_n \over 1+h_n}\cdot\hat\vu_n \str 0 \quad \hbox{in} \quad \Ld{\Rd} \,.
$$
If $(h_n^{-2}|\vf_n|^2)$ is equiintegrable, i.e.
$$
\svaki{\eps>0}\postoji{\delta>0}\svaki{n\in\N}\svaki{A\in\lB} \quad \lambda(A)<\delta 
	\povlaci \int_A \biggl({|\vf_n|\over h_n}\biggr)^2 d\mxi <\eps \,, 
$$ 
where $\lB$ denotes the Borel $\sigma$-algebra and $\lambda$ the Lebesgue measure on $\Rd$, then we also have
$$
{\vf_n \over h_n}\cdot\hat\vu_n \str 0 \quad \hbox{in} \quad \Ld{\Rd} \,.
$$
\vskip-2mm
}
{
By the continuity of Fourier transform from $\rL^1$ to $\rL^\infty$, $S:=\sup_n  \|\hat\vu_n\|_{\rL^\infty}$
is finite. From the equiintegrability of $(h_n^{-2}|\vf_n|^2)$ we have that for given $\eps>0$ we can take
a ball $K=\Ko \vnul r$ such that $\|h_n^{-1}\vf_n\|_{\Ld{K;\Cr}}<{\eps\over 2S}$. 

Thus we get
$$
\biggl\|{\vf_n \over h_n}\cdot\hat\vu_n\biggr\|_{\Ld{K}} 
 \mj \biggl\|{\vf_n\over h_n}\biggr\|_{\Ld{K;\Cr}}\|\hat\vu_n\|_{\Lb{\Rd;\Cr}} 
	< {\eps\over 2} \,,
$$
and in order to prove the lemma it remains  to estimate the same term over $\Rd\setminus K$.

By uniform boundedness from below, for $C={1/\tilde C}+1$, outside $K$ we have
$$
{|\vf_n\cdot\hat\vu_n|\over h_n} \mj C{|\vf_n\cdot\hat\vu_n|\over 1+h_n} \,,
$$
implying
$$
\biggl\|{\vf_n \over h_n}\cdot\hat\vu_n\biggr\|_{\Ld{\Rd\setminus K}} \mj C\biggl\|{\vf_n \over 1+h_n}
	\cdot\hat\vu_n\biggr\|_{\Ld{\Rd\setminus K}} 
	\mj C\biggl\|{\vf_n \over 1+h_n}\cdot\hat\vu_n \biggr\|_{\Ld\Rd} \,.
$$
As the last term above tends to 0 by assumption, there is an  $n_0\in\mN$ such that
for any $n\vj n_0$ we have $\|{\vf_n \over 1+h_n}\cdot\hat\vu_n \|_{\Ld\Rd}<{\eps\over 2C}$.
}

The localisation principle has been  obtained under  the assumption on the right hand side (\cfor{ConvRhsLoc}), 
which is adapted to the problem under consideration. However, we shall prove that the assumed convergence is 
also a necessary condition, similarly as it was the case with the localisation principle for H-measures 
[\ref{Tgth}, Lemma 28.18].

\tmd{
Let $\vu_n\dscon\vnul$ in $\Ll2{\domena;\Cr}$ define one-scale H-measure $\jhmv$ with the characteristic
length $(\eps_n)$, and $\mA^\malpha\in\Cp{\domena;\Mx{q\times r}{\mC}}$ for $l\mj|\malpha|\mj m$. 

If the conclusion of previous theorem holds, then (\cfor{ConvRhsLoc}) must also hold, with $\vf_n$ 
being defined by (\cfor{locprEqJhmA}). 
}
{
By part I of the proof of Theorem 8 and Lemma 6 it is equivalent to prove 
that for an arbitrary $\ph\in\Cbc\domena$
$$
{\lF\Bigl(\sum_{l\mj|\malpha|\mj m}\eps_n^{|\malpha|-l}\partial_\malpha(\ph\mA^\malpha\vu_n)
	\Bigr) \over 1+h_n} \str \vnul
$$
in $\Ld{\Rd;\Cq}$.
Moreover, after applying the Fourier transform and taking into account simple inequality
$1+|\mxi|^l+\eps_n^{m-l}|\mxi|^m \vj |\mxi|^l+\eps_n^{m-l}|\mxi|^m$,
it is sufficient to show that
$$
\vw_n := \sum_{l\mj|\malpha|\mj m}(2\pi i)^{|\malpha|} {(\eps_n\mxi)^\malpha \over 
	|\eps_n\mxi|^l + |\eps_n\mxi|^m} \widehat{\mA^\malpha\ph\vu_n} \str \vnul
$$
in $\Ld{\Rd;\Cq}$.
However, by the definition of one-scale H-measure we have
$$
\lim_n \int\limits_\Rd \vw_n\otimes\vw_n \,d\mxi = 
	\Dupp{\overline{\mmp_1\jhmv^\top\mmp_1^\ast}}{|\ph|^2\kvtpr 1} = \mnul \,,
$$
where in the last equality we have used the assumption (\cfor{locprJhm}), with $\mmp_1$  given by 
(\cfor{SymbJhmA}), which completes the proof.  
}

Let us now go back to our main result, 
Theorem 8, which easily extends to the case where the coefficients suitably depend on $n$.

\krd{
Under the same assumptions on $(\vf_n)$ and $(\eps_n)$ as in Theorem 8,
for $\vu_n \dscon \vnul$ in $\Ll2{\domena;\Cr}$ additionally assume that
$$
\sum_{l\mj|\malpha|\mj m}\eps_n^{|\malpha|-l}\partial_\malpha(\mA_n^{\malpha}\vu_{n}) 
	= \vf_n \,,
\leqno(\rfor{locprEqJhmB})
$$
where $\mA^\malpha_n\in\Cp{\domena;\Mx{q\times r}{\mC}}$, such that for any 
$\malpha$ the sequence $\mA^\malpha_n \str \mA^\malpha$ in the space $\Cp{\domena;\Mx{q\times r}{\mC}}$
(in other words, $\mA^\malpha_n$ converges locally uniformly to $\mA^\malpha$).
Then we obtain the same conclusion as in Theorem 8, with the same symbol given by (\cfor{SymbJhmA}).
}
{
It is enough to prove that 
$$
\sum_{l\mj|\malpha|\mj m}\eps_n^{|\malpha|-l}\partial_\malpha\Bigl((\mA^{\malpha}_n-\mA^\malpha)
	\vu_{n}\Bigr)
$$
satisfies the convergence condition for the right hand side (\cfor{ConvRhsLoc}).
As $(\mA^\malpha_n-\mA^\malpha)\vu_n$ strongly converges to $\vnul$ in $\Ll{2}{\domena;\Cq}$,
so $\partial_\malpha\Bigl((\mA^{\malpha}_n-\mA^\malpha) \vu_{n}\Bigr)$ strongly converges to
$\vnul$ in $\Hl{-|\malpha|}{\domena;\Cq}$, hence by Lemma 7 we have the claim. 
}

\Napbr{%
The powers of scale $(\eps_n)$ in (\cfor{locprEqJhmB}) strictly follow the order of derivatives 
(up to constant $l$). However, by a slight modification of the convergence assumption on the 
right side, the above result  easily extends to equations allowing  more general relations between 
these terms.

Namely, let us suppose that we start with the equalities
$$
\sum_{l\mj|\malpha|\mj m}\eps_n^{|\malpha|-l+\delta_\alpha}\partial_\malpha(\mA_n^{\malpha}\vu_{n}) 
	= \vf_n \,,
$$
with $\delta_\alpha\in \R$ . Multiplication of the above relation by $\eps_n^{-\delta}$, 
where $\delta:= \min \delta_\alpha$, implies
$$
\sum_{l\mj|\malpha|\mj m}\eps_n^{|\malpha|-l}\partial_\malpha\left((\eps_n^{\delta_\alpha-\delta}
	\mA_n^{\malpha})\vu_{n}\right) =\eps_n^{-\delta}  \vf_n \,.
$$
For sequence $(\vf_n)$ we assume that $(\eps_n^{-\delta} \vf_n)$ satisfies (\cfor{ConvRhsLoc}), 
which is a stronger assumption for $\delta>0$, and weaker for $\delta<0$, compared to the assumptions 
of the last corollary. Its application finally provides the localisation result (with symbol 
$\mmp_1$ not containing terms for which $\delta_\alpha>\delta$). 
}
\smallskip

\podnaslov{The case $c:=\lim_n{\eps_n\over\omega_n}\in\zi0\infty$}


\tmd{
Let $\vu_n\dscon\vnul$ in $\Ll2{\domena;\Cr}$ and satisfy (\cfor{locprEqJhmB}),
where $\mA^\malpha_n \str \mA^\malpha$ in the space $\Cp{\domena;\Mx{q\times r}{\mC}}$,
while $(\vf_n)$ is a sequence in\/  $\Hl{-m}{\domena;\Cq}$ satisfying (\cfor{ConvRhsLoc}).

If $(\eps_n)$ and  $(\omega_n)$ are sequences of positive numbers such that 
$c=\lim_n{\eps_n\over\omega_n}$ exists (in $\zi0\infty$),
then the associated one-scale H-measure $\jhmv^{(\omega_n)}$ satisfies:
$$
\mmp_c\Bigl(\jhmv^{(\omega_n)}\Bigr)^\top = \mnul \,,
\leqno(\rfor{locprJhmc})
$$
where with respect to the value of $c$ we have
\smallskip

\item{a)} $c=0$: 
$$
\mmp_0(\mx,\mxi):= \sum_{|\malpha|=l}(2\pi i)^l {\mxi^{\malpha}\over |\mxi|^l+|\mxi|^m} 
	\mA^\malpha(\mx)	\;,
\leqno(\rfor{SymbJhmB})
$$
 
\item{b)} $c\in\oi0\infty$:
$$
\mmp_c(\mx,\mxi):= \sum_{l\mj|\malpha|\mj m}(2\pi i c)^{|\malpha|}	
	{\mxi^{\malpha}\over |\mxi|^l+|\mxi|^m} \mA^\malpha(\mx)	\;,
\leqno(\rfor{SymbJhmC})
$$   

\item{c)} $c=\infty$:
$$
\mmp_\infty(\mx,\mxi):= \sum_{|\malpha|=m}(2\pi i)^m {\mxi^{\malpha}\over |\mxi|^l+|\mxi|^m} 
	\mA^\malpha(\mx) \,.
\leqno(\rfor{SymbJhmD})
$$ 
\vskip-2mm
}
{
In both (a) and (b) cases, by rewriting (\cfor{locprEqJhmB}) we infer
$$
\sum_{l\mj|\malpha|\mj m}\omega_n^{|\malpha|-l}\partial_\malpha(\mB_n^{\malpha}\vu_{n}) 
	= \vf_n \,,
$$
where $\mB_n^\malpha:=\Bigl({\eps_n\over\omega_n}\Bigr)^{|\malpha|-l}\mA_n^\malpha$.
As for $n$ large enough there is $\kappa>0$ such that $\omega_n>\kappa
\eps_n$, it can been easily seen  that the right hand side also satisfies 
$$
\svaki{\ph\in\Cbc{\domena}} \qquad {\fp{\ph\vf_n}\over 1+|\mxi|^l+\omega_n^{m-l}|\mxi|^m} \str \vnul  
	\quad \hbox{in} \quad  \Ld{\Rd;\Cq} \,.
$$ 
Taking into account Corollary 7,  and the fact that $\mB_n^\malpha$ converges to 
$c^{|\malpha|-l}\mA^\malpha$ (in particular to $\mA^\malpha$ for $c=0$ and $|\malpha|= l$), we get the claim. 
\smallskip

\noindent {\bf c)}\ 
After multiplying (\cfor{locprEqJhmB}) by  $\Bigl({\omega_n\over\eps_n}\Bigr)^{m-l}$, we obtain
$$
\sum_{l\mj|\malpha|\mj m}\omega_n^{|\malpha|-l}\partial_\malpha(\tilde\mB_n^{\malpha}\vu_{n}) 
	= \vg_n \,,
$$
where $\tilde\mB_n^\malpha:=\Bigl({\omega_n\over\eps_n}\Bigr)^{m-|\malpha|}\mA_n^\malpha$, while
$\vg_n:=\Bigl({\omega_n\over\eps_n}\Bigr)^{m-l}\vf_n$. In order to apply Corollary 7 it remains 
to check the condition satisfied by the right hand side $\vg_n$.

For $\ph\in\Cbc\domena$ we have
$$
\eqalign{
{|\fp{\ph \vg_n|}\over 1+|\mxi|^l+\omega_n^{m-l}|\mxi|^m} 
&= {|\fp{\ph \vf_n|}\over \Bigl({\eps_n\over\omega_n}\Bigr)^{m-l}(1+|\mxi|^l)+\eps_n^{m-l}|\mxi|^m} \cr
&< {|\fp{\ph \vf_n|}\over 1+|\mxi|^l+\eps_n^{m-l}|\mxi|^m} \str \vnul \quad \hbox{in} \quad 
	\Ld{\Rd} \,, \cr
}
$$
where the last inequality is valid for sufficiently large $n$. Thus Corollary 7,
together with the fact that $\tilde\mB_n^\malpha$ converges to zero for $|\malpha|\mj m-1$,
provides the claim.
}

\Napbr{%
Note that for $c=0$ and $c=\infty$ (with $l<m$) the previous theorem does not provide any information 
about the structure of one-scale H-measure on $\Sigma_\infty$ and $\Sigma_0$, respectively. 
This is caused by the fact that (\cfor{SymbJhmB}) is zero on $\Sigma_\infty$, while
(\cfor{SymbJhmD}) vanishes on $\Sigma_0$. 
}
\smallskip


Part (c) of the previous theorem, in the special case when $(\eps_n)$ is bounded by positive constants 
both from below and above, can be improved as a consequence of Lemma 6.

\tmd{
Under the assumptions of Theorem 10 for $c=\infty$, additionally assume that $0<\eps_0\mj\eps_n\mj\eps_\infty$.
Then the associated one-scale H-measure $\jhmv^{(\omega_n)}$ satisfies:
$$
\mmp_{pr}\Bigl(\jhmv^{(\omega_n)}\Bigr)^\top = \mnul \,,
$$
where
$$
\mmp_{pr}(\mx,\mxi)= \sum_{|\malpha|=m}(2\pi i)^m \Bigl({\mxi\over |\mxi|}\Bigr)^\malpha \mA^\malpha(\mx) \,.
\leqno(\cfor{Symbppr})
$$
\vskip-2mm 
}
{
Without loosing generality, we may assume that the coefficients do not depend on $n$, 
i.e.~$\mA^\malpha_n =\mA^\malpha$. The general case can be treated as in Corollary 7.

By repeating part I of the proof of Theorem 8, we get (\cfor{locprEqLoc}). We aim to show 
that the lower order terms converge to zero in the sense of (the right hand side) 
convergence (\cfor{ConvRhs}). In fact, Lemma 6 additionally gives us that $\tilde\vf_n$ 
converges to $\vnul$ in $\H{-m}{\Rd;\Cq}$, same as terms 
$\eps_n^{|\malpha|-l}\partial_\malpha(\ph\mA^\malpha\vu_n)$ for $|\malpha|<m$. This implies
$$
\sum_{|\malpha|=m} \eps_n^{m-l}\partial_\malpha(\mA^\malpha\ph\vu_n) \str \vnul \quad \hbox{in} 
	\quad \H{-m}{\Rd;\Cq} \,,
$$
and further
$$
\sum_{|\malpha|=m}(2\pi i)^m {\mxi^\malpha\fp{\mA^\malpha\ph\vu_n}\over 1+|\mxi|^m} \str 
	\vnul \quad \hbox{in} \quad \Ld{\Rd;\Cq} \,,
$$
as $\eps_n>\eps_0$. Furthermore, Lemma 9 gives us
$$
\sum_{|\malpha|=m}(2\pi i)^m {\mxi^\malpha\over |\mxi|^m} \fp{\mA^\malpha\ph\vu_n} \str \vnul
	\quad \hbox{in} \quad \Ld{\Rd;\Cq} \,.
$$
After forming the tensor product with $\psi(\omega_n\cdot)\fp{\ph_1\vu_n}$, where  
$\psi\in\Cp{\kob}$, $\ph_1\in\Cbc\domena$, and integrating, we obtain 
$$
\eqalign{
\mnul =& \lim_n \int_\Rd\psi(\omega_n\mxi) \Biggl(\sum_{|\malpha|= m}(2\pi i)^m
	{(\omega_n\mxi)^\malpha \over |\omega_n\mxi|^m} 
	\fp{\mA^{\malpha}\ph\vu_n}\Biggr) \otimes \Bigl(\fp{\ph_1\vu_n}\Bigr) \,d\mxi \cr
=& \Dupp{\overline{\sum_{|\malpha|= m}(2\pi i)^m {\mxi^\malpha 
	\over |\mxi|^m} \mA^{\malpha}}\jhmv}{\ph\bar{\ph_1}\kvtpr\psi} \,. \cr
}
$$
As $\ph_1$ is arbitrary, we can choose it  equal to $1$ on $\supp\ph$,
which gives the claim.
}

Let us note that in the previous proof we have not assumed that $(\eps_n)$
has a unique accumulation point. 

Theorem 11 is an improvement over Theorem 10(c) indeed, as the symbol (\cfor{Symbppr}) 
does not vanish on the sphere $\Sigma_0$ in the compactified domain $\kob$, as it 
was the case with symbol (\cfor{SymbJhmD}). In such a way it provides information 
on the structure of the considered one-scale H-measure on its entire domain.

\Napbr{%
If $\vu_n\dscon\vu\not= \vnul$ in $\Ldl{\Omega;\Cr}$, we cannot associate a one-scale 
H-measure to $(\vu_n)$,  but, like in [\ref{Tgth}, Theorem 28.7], we can still derive the
localisation principle for corresponding null sequence $(\vu_n-\vu)$. Namely, let 
$\vu_n$ satisfy (\cfor{locprEqJhmB}), where $\mA^\malpha_n$ converges to $\mA^\malpha$ 
in $\Cp{\Omega;\Mx{q\times r}\mC}$, let $\eps_n\to 0^+$,
and let $\vf_n\in\Hl{-m}{\Omega;\Cq}$ be such that for any $\ph\in\Cbc\Omega$, 
${\widehat{\ph\vf_n} \over 1+k_n}$ is precompact in the strong topology of $\Ld{\Rd;\Cq}$.

As in the first step of the proof of Theorem 8 we localise the equation by
multiplying it by $\ph\in\Cbc\Omega$ in order to get
$$
\sum_{l\mj|\malpha|\mj m} \eps_n^{|\malpha|-l}\partial_{\malpha}
	(\mA_n^\malpha\ph\vu_n) = \vf_{1,n} \,,
$$
where 
$$
\vf_{1,n} := \ph\vf_n - \sum_{l\mj|\malpha|\mj m}\sum_{\vnul\ne\mbeta\mj\malpha}
	(-1)^{|\mbeta|}\binom{\malpha}{\mbeta}\eps_n^{|\malpha|-l}
	\partial_{\malpha-\mbeta}\Bigl( (\partial_\mbeta\ph)\mA_n^\malpha\vu_n \Bigr) \,.
$$
Since $\partial_{\malpha-\mbeta}\Bigl((\partial_\mbeta\ph)\mA_n^\malpha\vu_n\Bigr)\str
\partial_{\malpha-\mbeta}\Bigl((\partial_\mbeta\ph)\mA^\malpha\vu\Bigr)$ in 
$\H{-|\malpha|}{\Rd;\Cq}$, by Lemma 7 and the assumption on $(\vf_n)$, 
we have that ${\widehat{\vf_{1,n}}\over 1+k_n}$ is precompact in the strong topology of 
$\Ld{\Rd;\Cq}$.

Next, after subtracting $\sum_{l\mj|\malpha|\mj m}\eps_n^{|\malpha|-l}
\partial_{\malpha}(\mA_n^\malpha\ph\vu)$ in the localised equation, by the 
linearity of differential operator on the left hand side, $\vu_n$ is 
replaced by $\vu_n-\vu$, so the whole left hand side converges to zero 
in the sense of distributions. On the other hand,  on the right hand side 
we obtained a sequence $(\vf_{2,n})$ which still satisfies that 
${\widehat{\vf_{2,n}}\over 1+k_n}$ is precompact, while
$\vf_{2,{n}} \dscon \vnul$ in the sense of distributions. By 
Lemma 10 below we get that $(\vf_{2,n})$ satisfies (\cfor{ConvRhs}). Hence all the 
assumptions of Theorem 10 are satisfied and the one-scale H-measure 
 associated to $(\vu_n - \vu)$ satisfies (\cfor{locprJhmc}).
}
\smallskip

\lmd{
Let $\eps_n\to0^+$ and let $(\vf_n)$ be a sequence of vector functions from 
$\H{-m}{\Rd;\Cq}$ such that ${\widehat\vf_n\over 1+k_n}$ is precompact 
in the strong topology of\/ $\Ld{\Rd;\Cq}$. Then there exists a subsequence 
$(\vf_{n'})$ and an $\vf\in\H{-l}{\Rd;\Cq}$ such that
$$
{\widehat\vf_{n'} - \widehat\vf \over 1+k_{n'}} \str \vnul
$$
in $\Ld{\Rd;\Cq}$. In particular, if $\vf_n\dscon \vnul$ in the sense of distributions, 
then $\vf\equiv\vnul$ and the whole sequence $(\vf_n)$ satisfies (\cfor{ConvRhs}).
}
{
By the assumptions, there exists a subsequence $(\vf_{n'})$ and 
$\vF\in\Ld{\Rd;\Cr}$ such that 
${\widehat\vf_{n'}\over 1+k_{n'}} \str \vF$ in $\Ld{\Rd;\Cq}$.
By using the Lebesgue dominated convergence theorem it follows that
$$
{\widehat\vf_{n'}-\widehat\vf\over 1+k_{n'}} \str \vnul
$$ 
in $\Ld{\Rd;\Cq}$, with $\vf:=\bar\lF\Bigl((1+|\mxi|^l)\vF\Bigr)
\in\H{-l}{\Rd;\Cq}$.
By the remark after Lemma 6, the above implies $\vf_{n'}\str\vf$ in $\H{-m}{\Rd;\Cq}$.

If in addition we have that $\vf_n\dscon\vnul$ in the sense of distribution, then 
by the uniqueness of the limit, $\vf\equiv\vnul$ and hence $(\vf_{n'})$
satisfies (\cfor{ConvRhs}). Since this holds for every subsequence the claim follows. 
}


\odjeljak{Some applications}

\podnaslov{Localisation principles revisited}

As we have already noticed, from the one-scale H-measure one can get both the corresponding
semiclassical measure, as well as the corresponding H-measure. In the following corollaries we
show that the localisation principles for these two objects are also consequences of the
localisation principle for one-scale H-measures.

\krd{ {\bf (localisation principle for H-measures)} \
Under the assumptions of Theorem 11 one has the conclusion of Theorem 2.
}
{
Since the symbol (\cfor{Symbppr}) is homogeneous of order zero in $\mxi$, by Theorem
11 and Corollary 6 we get the claim.
}

The statement of the previous corollary is more general then the 
one of Theorem 2 since it allows for coefficients depending on  $n$ in the system under consideration.
However, such a generalisation  can also be obtained directly, by using the 
standard techniques for H-measures and the approach presented in Corollary 7.

\krd{ {\bf (localisation principle for semiclassical measures)} \
Under the assumptions of Theorem 10 we have 
$$
\mmp\pkmv^\top = \mnul \,,
$$
where $\pkmv$ is a semiclassical measure with characteristic length $(\omega_n)$ associated 
to the sequence $(\vu_n)$, while
$$
\mmp(\mx,\mxi):= 
\left\{
\matrix{
\sum_{|\malpha|=l}(2\pi i)^l\mxi^\malpha\mA^\malpha(\mx) 
	& , & \lim_n{\eps_n\over\omega_n}= 0 \cr
\sum_{l\mj|\malpha|\mj m} (2\pi i c)^{|\malpha|}\mxi^\malpha\mA^\malpha(\mx)
	& , & \lim_n{\eps_n\over\omega_n}=c\in\oi 0\infty \cr
\sum_{|\malpha|=m}(2\pi i)^m\mxi^\malpha\mA^\malpha(\mx) 
	& , & \lim_n{\eps_n\over\omega_n}= \infty \cr
}
\right. \,.
$$
\vskip-2mm  
}
{
Let $\lim_n{\eps_n\over\omega_n}=c\in\oi0\infty$. For 
$\psi\in\lS(\Rd)$ we have that $\mxi\mapsto (|\mxi|^l+|\mxi|^m)\psi(\mxi)$
is in $\Cp{\kob}$, so after applying the localisation principle, i.e.~Theorem 10(b),
to $\ph\kvtpr(|\mxi|^l+|\mxi|^m)\psi$, for $\ph\in\Cc\domena$
one obtains 
$$
\eqalign{
\mnul &= \Dupp{\sum_{l\mj|\malpha|\mj m}(2\pi i c)^{|\malpha|}
	{\mxi^\malpha \over |\mxi|^l+|\mxi|^m} \mA^{\malpha}\jhmv^\top }
	{\ph\kvtpr (|\mxi|^l+|\mxi|^m)\psi} \cr
&= \sum_{l\mj|\malpha|\mj m}\Dupp{\mA^\malpha\jhmv^\top }{\overline{(2\pi i c)^{|\malpha|}}
	\ph\kvtpr \mxi^\malpha\psi} \cr
&= \sum_{l\mj|\malpha|\mj m}\Dupp{\mA^\malpha\pkmv^\top }{\overline{(2\pi i c)^{|\malpha|}} 
	\ph\kvtpr \mxi^\malpha\psi} 
	= \Dupp{\sum_{l\mj|\malpha|\mj m} (2\pi i c)^{|\malpha|}
	\mxi^\malpha\mA^\malpha \pkmv^\top }{\ph\kvtpr\psi} \,, \cr
}
$$
where in the third equality the fact that $\mxi^\malpha\psi\in\lS(\Rd)$ was used,
as well as Corollary 6.

The cases $\lim_n{\eps_n\over\omega_n}=0$ and $\lim_n{\eps_n\over\omega_n}=\infty$ follow
in the same manner. 
}

The theory of semiclassical measures was developed by means of
pseudodifferential calculus which requires smooth test functions. 
That was the main obstruction to weaker assumptions on the 
coefficients of the observed system. 
Following Tartar's approach, the theory of one scale H-measures is constructed without 
the smoothness assumption, thus enabling us  to generalise existing results for 
semiclassical measures to merely continuous functions without any technical
difficulties. Furthermore, another improvement is based on the
convergence of the right hand side (\cfor{ConvRhsLoc}), which is weaker 
than $\rL^2_{\rm loc}$  convergence required in Theorem 4. 

In addition, Corollary 9 provides  a complete analysis of dependence on the relation 
between characteristic lengths of semiclassical measures and systems. This is useful in 
situations where the characteristic length of the observed problem cannot be easily identified, 
or where the same measure has to be applied to two (or more) systems with different characteristic 
scales.    

Let us stress once more that in Theorem 4 the statement was valid for 
an arbitrary bounded sequence in $\rL^2_{\rm loc}$, which has not been the 
case in the previous theorem. However, according to Remark 8 at the 
end of the previous section, we can still obtain information about the
semiclassical measure associated to $(\vu_n-\vu)$, where $\vu_n\rightharpoonup\vu$. 
Furthermore, this result is not weaker then the result of Theorem 4 since 
$(\vu\otimes\vu)\lambda\kvtpr\delta_\vnul$ cancels all the terms of symbol
$\mmp_{sc}$ from Theorem 4. Indeed, nontrivial is only the fact that 
$\mA^\vnul(\vu\otimes\vu)\lambda\kvtpr\delta_\vnul=\mnul$, which is a 
consequence of $\mA^\vnul\vu=\vnul$, following from the assumptions of 
Theorem 4.

\podnaslov{Compactness by compensation with a characteristic length}

One of the main purposes of one-scale H-measures, like of all microlocal defect tools, 
is to  investigate conditions resulting in strong convergence of a bounded $\LLd$ sequence 
under consideration. Hereby,  localisation principles presented above play a significant r\^ole
as they, possibly accompanied by some additional information, enable one to deduce that a 
(component of) measure $\jhmv$ associated to $\vu_n$ is trivial. This implies then either strong 
$\LLd$ convergence of $(u_n^i)$  (for a zero diagonal component) or vague convergence of 
$(u_n^i \bar u_n^j)$ (for a vanishing corresponding off-diagonal element). In the 
following example we shall present such a result obtained by the localisation principle 
derived in the previous section.

\Pr{%
Let $\domena\subseteq\R^2$ be open, and let $\vu_n:=(u_n^1, u_n^2) \dscon \vnul$
in $\Ll{2}{\domena;\mC^2}$ satisfies
$$
\left\{
\eqalign{
u_n^1 + \eps_n\partial_{x_1}(a_1u_n^1) &= f_n^1 \cr
u_n^2 + \eps_n\partial_{x_2}(a_2u_n^2) &= f_n^2 \cr
}\right. \,,
$$
where $\eps_n\to0^+$, $\vf_n:=(f_n^1, \ f_n^2)\in\Hl{-1}{\domena;\mC^2}$ satisfies 
(\cfor{ConvRhsLoc}) (with $l=0, m=1$), while $a_1, a_2\in\Cp{\domena;\R}$, $a_1, a_2\ne 0$ 
everywhere.

By the localisation principle for one-scale H-measure $\jhmv$ with characteristic length 
$(\eps_n)$ (i.e.~$c=1$) associated to $(\vu_n)$ we get the relation
$$
\left( {1\over 1+|\mxi|}
	\left[
	\matrix{
	1&0\cr
	0&1\cr
	}\right]
	+{2\pi i\xi_1\over 1+|\mxi|}
	\left[
	\matrix{
	a_1(\mx)&0\cr
	0&0\cr
	}\right]
	+{2\pi i\xi_2\over 1+|\mxi|}
	\left[
	\matrix{
	0&0\cr
	0&a_2(\mx)\cr
	}\right]
	\right) \jhmv^\top = \mnul \,,
$$
whose $(1,1)$ component reads
$$
\left( {1\over 1+|\mxi|} + i{2\pi\xi_1 \over 1+|\mxi|}a_1(\mx) \right)\jhm^{11} = 0 \,.
$$
Since measure $\jhm^{11}$ is non-negative, we can separate its real and imaginary part in order to obtain
$$
{1\over 1+|\mxi|}\jhm^{11} = 0 \,, \quad {\xi_1\over 1+|\mxi|}\jhm^{11} = 0 \,,
\leqno(\rfor{locprAppl})
$$
where we have used the assumption that $a_1\ne 0$. Since ${1\over 1+|\mxi|}$ differs from zero on 
$\pkob\Rd\setminus\Sigma_\infty$, from the first equation we get
$\supp\jhmv^{11}\subseteq\domena\times\Sigma_\infty$. Similarly,  from the second 
one we deduce $\supp\jhm^{11}\subseteq\domena\times(\Sigma_0 \cup \{\xi_1=0\})$, so at the 
end we have 
$$
\supp\jhm^{11}\subseteq \domena\times\{\infty^{(0,-1)}, \infty^{(0,1)}\} \,.
 \leqno(\rfor{suppmu})
$$
Analogously, from the $(2,2)$ component we get 
$$
\supp\jhm^{22} \subseteq \domena\times\{\infty^{(-1,0)}, \infty^{(1,0)}\} \,,
$$
which by Corollary 2 implies that $\jhm^{12} = \jhm^{21} = 0$. The very definition of 
one-scale H-measures gives $u_n^1 \bar{u_n^2} {\buildrel \ast \over\dscon} \, 0$.
}

\Napbr{%
If in the previous example we additionally have that $(u_n^1)$ is $(\eps_n)$-oscillatory, then 
Corollary 5 implies $\jhm^{11}(\domena\times\Sigma_\infty)=0$, hence by (\cfor{suppmu}) $\jhm^{11}$
vanishes everywhere, which gives $u_n^1\str 0$ in $\Ldl{\domena;\mC}$. As the product of a strongly 
and a weakly converging 
sequence, the convergence $u_n^1 \bar{u_n^2} {\buildrel \ast \over\dscon} \, 0$ is now a 
trivial consequence. An example of such a situation is the case when $f_n^1$ is a sequence 
of oscillations, resulting in a solution $u_n^1$  oscillating at the same 
frequency as $f_n^1$. Therefore, with  an appropriate choice of the frequency, one
obtains solutions satisfying the $(\eps_n)$-oscillatory property. 

Moreover, having an  $(\eps_n)$-oscillating sequence of solutions,  we can remove the assumption 
$a_1\ne 0$ in Example 5, since from the first equation in (\cfor{locprAppl}) we already have 
$\supp\jhm^{11}\subseteq \domena\times\Sigma_\infty$, hence still $\jhm^{11}=0$. Thus in the special 
case when $a_1$ is trivial, this provides an alternative proof of Theorem 7 for the case $l=0$.   
}
\smallskip

The previous result is an example of a compactness by compensation theory, investigating conditions 
under which a  product of weakly convergent sequences converges (in some suitable sense) to the 
product of corresponding limits. The theory was introduced by Fran\c cois Murat and Luc Tartar 
(see [\ref{BCDM}] and references therein), first for constant coefficients and then generalised 
to continuous coefficients in [\ref{Tprse}] by means of H-measures  (cf.~[\ref{Tgth}, Cor.~28.11] 
and [\ref{Gmd}, Theorem 2]; some recent extensions can be found in [\ref{MMip}, \ref{EJPahp}]).
Of course, the obtained results do not treat equations containing a
characteristic length, which is now possible by  the obtained localisation principle for 
one-scale H-measures.

\tmd{
{\bf (a variant of compactness by compensation)} \
Let $\vu_n\dscon\vu$ in $\Ll{2}{\domena;\Cr}$ satisfy (\cfor{locprEqJhmB}), where 
$\mA^\malpha_n\str\mA^\malpha$ in $\Cp{\domena;\Mx{q\times r}\mC}$, let $\eps_n\to0^+$, and 
$\vf_n\in\Hl{-m}{\domena;\Cq}$ be such that for any $\ph\in\Cbc\domena$
$$
{\widehat{\ph\vf_n}\over 1+k_n} 
$$
is precompact in $\Ld{\Rd;\Cq}$. Furthermore, let $Q(\mx;\mlambda):=
\mQ(\mx)\mlambda\cdot\mlambda$, where $\mQ\in\Cp{\domena;\Mx{r}\mC}$, $\mQ^\ast=\mQ$, 
is such that $Q(\cdot;\vu_n){\buildrel \ast \over\dscon}\,\nu$ in $\lM(\domena)$.

Then we have
\item{a)} $\postoji{c\in\zi0\infty}\svaki{(\mx,\mxi)\in\domena\times\kob}
	\svaki{\mlambda\in\Lambda_{c;\mx,\mxi}}\ Q(\mx;\mlambda)\vj 0 \povlaci \nu\vj Q(\cdot,\vu)$,
\item{b)} $\postoji{c\in\zi0\infty}\svaki{(\mx,\mxi)\in\domena\times\kob}
	\svaki{\mlambda\in\Lambda_{c;\mx,\mxi}}\ Q(\mx;\mlambda)=0 \povlaci \nu=Q(\cdot,\vu)$,

\noindent where
$$
\Lambda_{c;\mx,\mxi} := \{\mlambda\in\Cr : \mmp_c(\mx,\mxi)\mlambda=\vnul\} \,,
$$ 
and $\mmp_c$ is given in Theorem 10. 
}
{
Let $\jhmv$ be the one-scale H-measure associated to $(\vu_n-\vu)$, with a characteristic length $(\omega_n)$
such that $\lim_n{\eps_n\over\omega_n}=c$.

The limit $\nu$ we can express as a sum of $Q(\cdot;\vu)$ and the corresponding linear combinations of 
components of matrix defect measure $\mnu_D$ associated to $(\vu_n-\vu)$. Namely, for an arbitrary 
$\ph\in\Cc\domena$ we have
$$
\eqalign{
\Dupp\nu\ph =& \Dupp{Q(\cdot;\vu)}{\ph} + \Dupp{\mQ\cdot\mnu_D^\top}{\ph} \cr
=& \Dupp{Q(\cdot;\vu)}{\ph} + \Dupp{\mQ\cdot\jhmv^\top}{\ph\kvtpr 1} \,, \cr
}
$$
where  in the second equality we have used Corollary 4. 

We shall demonstrate that $\mQ\cdot\jhmv^\top$ is a non-negative measure (or trivial measure 
in part (b)), which implies the claim. 

Let us introduce a non-negative Radon measure $\rho:=\tr \jhmv=\sum_{i=1}^r\jhm^{ii}$. According to 
Corollary 2, for each $i,j\in1..r$ we have that $\jhm^{ij}\ll\rho$, so by the Lebesgue-Radon-Nikod\'ym theorem 
there exists a measurable function $M^{ij}(\mx,\mxi)$ such that $\jhm^{ij} = M^{ij}\rho$. 
Since $\jhmv$ is hermitian and non-negative, $\mM:=[M^{ij}]$ is also hermitian and non-negative. In 
terms of the matrix function $\mM$ it remains to prove that for $\rho$-a.e.~$(\mx,\mxi)$ we have 
$\mQ\cdot\mM^\top\vj 0$  (or equal to 0 in part (b)). 

By the localisation principle given in Theorem 10 (which remains valid for $\vu\not=\vnul$ based on Remark 8) 
we have $\mmp_c\jhmv^\top=\mnul$, hence $\mmp_c\mM^\top=\mnul$ ($\rho$-a.e.), which implies 
that for $\rho$-a.e.~$(\mx,\mxi)$ the columns of $\mM^\top$ are in the characteristic space $\Lambda_{c;\mx,\mxi}$. 
Let us take $(\mx,\mxi)\in\domena\times\kob$ arbitrary, but fixed. As $\mM^\top(\mx,\mxi)$ is
hermitian, there is an orthonormal basis $\{\mzeta_1,\dots,\mzeta_r\}$ of $\Cr$ composed of
eigenvectors of $\mM^\top(\mx,\mxi)$, which are of course in $\Lambda_{c;\mx,\mxi}$. 
Moreover, the non-negativity of $\mM^\top(\mx,\mxi)$ implies that all eigenvalues are non-negative,
therefore equal to $\|\mM^\top(\mx,\mxi)\mzeta_i\|$. Let us define
$$
\mlambda_i := \sqrt{\|\mM^\top(\mx,\mxi)\mzeta_i\|} \mzeta_i \,, 
	\quad i\in1..r \,. 
$$
We have $\mlambda_i\in\Lambda_{c;\mx,\mxi}$, $i\in1..r$, and it is not hard to check that 
$\mM^\top(\mx,\mxi) = \sum_{i=1}^r \mlambda_i\otimes\mlambda_i$, so by the assumption it follows that
$\mQ(\mx)\cdot\mM^\top(\mx,\mxi) =\sum_i Q(\mx; \mlambda_i)\vj 0$ (or equal to 0 in part (b)). 
The arbitrariness of $(\mx,\mxi)$ gives us the claim.    
}

\Napbr{%
Although the statement of the previous theorem is valid for all $c\in\zi0\infty$, the interesting case 
is only $c\in\oi0\infty$, as for the remaining ones the statement is trivial, at least for $l<m$. 
Indeed, in the case $c=0$ we have $\msnop_0(\mx,\mxi)=\mnul$ for $\mxi\in\Sigma_\infty$, implying 
$\Lambda_{0;\mx,\mxi}=\Cr$ for every $\mx$. Similarly, $\msnop_\infty$ vanishes for $\mxi\in\Sigma_0$, 
again resulting in $\mQ\vj \mnul$ (or equal to $\mnul$ in part (b)).
}

\Napbr{%
By splitting the matrix into hermitian and antihermitian part, we can generalise part (b) 
of the previous theorem to an arbitrary $\mQ\in\Cp{\domena;\Mx{r}{\mC}}$  
(cf.~[\ref{Gmd}, Theorem 2(ii)]).
}

Going back to Example 5, we can obtain the same conclusion by applying the last theorem.
Indeed, a straightforward computation leads to
$$
\bigcup_{\mx,\mxi}\Lambda_{1;\mx,\mxi} = \{(\lambda_1, 0) : \lambda_1\in\mC\}
	\cup\{(0, \lambda_2) : \lambda_2\in\mC\} \,,
$$
and on that set quadratic form $Q(\lambda):=\lambda_1\bar{\lambda_2}$ is zero. 
One needs to be aware that $Q$ is not a hermitian form, but Theorem 12 remains applicable according to the 
last remark.

\podnaslov{Concluding remarks}

We have shown that one-scale H-measures encompass features of both H-measures and semiclassical measures 
(Corollary 6), at the same time providing new properties (Corollary 5) and results.
Namely, by improving and generalising Tartar's localisation principle (theorems 10 and 11, Remark 8),
we have obtained a compactness by compensation result suitable for problems with a characteristic length 
(Theorem 12), opening a new approach for investigating partial differential equations with a characteristic 
length.  
Moreover, a more detailed insight into the condition satisfied by the right hand side (\cfor{ConvRhsLoc})
leads to various useful properties (lemmas 6, 7 and 10) and, more importantly, nice behaviour with respect to
$(\eps_n)$-oscillatory sequences (Theorem 7) and the corresponding partial differential equations with a 
characteristic length (Theorem 9), so it seems to be the right choice for further studies. 
Of course, in applications a need for different scaling $(\eps_n)$ in different variables on the left hand 
side of (\cfor{locprEqJhmB}) might arise, which can be partially overcome by the procedure presented in 
Remark 6.  However, a more systematic approach can be devised by introducing suitable variants [\ref{TmsHm}] 
(see also [\ref{AL}] for problems without a characteristic length and the scaling 1:2).

Further improvement might consist in deriving the {\sl propagation principle}, as it was the case with
classical [\ref{Tprse}, \ref{Asym}] and parabolic [\ref{AL}] H-measures or, on the other hand, with
semiclassical measures [\ref{CFKMS}, \ref{FKGL}, \ref{GMMP}]. According to that principle, the 
measure propagates along bicharacteristics of a differential operator under consideration. By means of 
Corollary 5, its direct consequence in the case of a finite speed of propagation (i.e.~for hyperbolic 
operators) will be $(\eps_n)$-oscillating property of a corresponding sequence of solutions, under 
the assumption of $(\eps_n)$-oscillatory initial conditions. Further results will rely on generalisations 
of existing applications of the propagation principle to problems with a characteristic length, as it 
was done here with the localisation principle.

\bigskip

\garmond

\noindent{\bf Conflict of Interest: }
The authors declare that they have no conflict of interest.

\bigskip

\podnaslov{References}

\tolerance=2000

\ite{Asym}
{\sc Nenad Antoni\'c}:
{\sl H-measures applied to symmetric systems},
{\it Proc.~Roy.~Soc.~Edinburgh}
{\bf 126A} (1996) 1133--1155.

\ite{ABu}
{\sc Nenad Antoni\'c, Kre\v simir Burazin}:
{\sl On certain properties of spaces of locally Sobolev functions, } in
{\it Proceedings of the Conference on applied mathematics and scientific computing},
Zlatko Drma\v c et al.~(eds.), Springer, 2005, pp.~109--120.


\ite{ALjmaa}
{\sc Nenad Antoni\'c, Martin Lazar}:
{\sl H-measures and variants applied to parabolic equations},
{\it J. Math.~Anal.~Appl.}~{\bf 343} (2008) 207--225.

\ite{ALrwa}
{\sc Nenad Antoni\'{c}, Martin Lazar}:
{\sl Parabolic variant of H-measures in homogenisation of a model problem based on Navier-Stokes equation},
{\it Nonlinear Anal.~B: Real World Appl.} {\bf 11} (2010) 4500--4512.

\ite{AL}   
{\sc Nenad Antoni\'c, Martin Lazar}:
{\sl Parabolic H-measures},
{\it J. Functional Analysis}
{\bf 265} (2013) 1190--1239.

\ite{AMaaa}
{\sc Nenad Antoni\'c, Darko Mitrovi\'c}:
{\sl H-distributions: an extension of H-measures to an ${\rm L}^p-{\rm L}^q$ setting},
{\it Abs.~Appl.~Analysis}
{\bf 2011} Article ID 901084 (2011) 12 pp.

\ite{JPAub}
{\sc Jean-Pierre Aubin}:
{\it Applied functional analysis},
Wiley,
2000.

\ite{Bfa}
{\sc Ha\"\i m Brezis}:
{\it Functional analysis, Sobolev spaces and partial differential equations},
Springer, 2011.

\ite{BCDM}
{\sc Marc Briane, Juan Casado-D\'iaz, Fran\c cois Murat}:
{\sl The div-curl lemma\/ ``trente ans apr\`es": an extension and an application to the 
G-convergence of unbounded monotone operators},
{\it J.~Math.~Pures Appl.}
{\bf 91} (2009) 476-494.

\ite{Burq}
{\sc Nicolas Burq}:
{\sl Mesures semi-classiques et mesures de d\' efaut,} 
{\it S\' eminare Bourbaki} (1996-97).

\ite{CFKMS}
{\sc R\'emi Carles, Clotilde Fermanian-Kammerer, Norbert J.~Mauser, 
	Hans Peter Stimming}:
{\sl On the time evolution of Wigner measures for Schr\"odinger 
	equations,}
{\it Comm.~Pure Appl.~Analysis}
{\bf 8} (2009) 559--585.

\ite{FKGL}
{\sc Clotilde Fermanian-Kammerer, Patrick G\'erard, Caroline Lasser}:
{\sl Wigner measure propagation and conical singularity for general initial data},
{\it Arch.~Ration.~Mech.~Anal.}
{\bf 209} (2013) 209--236.

\ite{FrInt}
{\sc Gilles A. Francfort}:
{\sl An introduction to H-measures and their applications,} in
{\it Progress in nonlinear differential equations and their applications} {\bf 68},
Birkh\"auser, 2006.

\ite{Gmd}
{\sc Patrick G\'erard}:
{\sl Microlocal defect measures,}
{\it Comm.~Partial Diff.~Eq.}
{\bf 16} (1991) 1761--1794.

\ite{Gmsc}
{\sc Patrick G\'erard}:
{\sl Mesures semi-classiques et ondes de Bloch,}
{\it Sem.~EDP 1990--91 (exp.~n$^\circ$ XVI)}, Ecole Polytechnique, Palaiseau (1991).

\ite{Goce}
{\sc Patrick G\'erard}:
{\sl Oscillations and concentration effects in semilinear dispersive wave equations,}
{\it J. Functional Analysis} {\bf 141} (1996) 60--98.

\ite{GMMP}
{\sc Patrick G\'erard, Peter A.~Markowich, Norbert J.~Mauser, 
	Fr\'ed\'eric Poupaud}:
{\sl Homogenization limits and Wigner transforms,}
{\it Comm.~Pure Appl.~Math.}
{\bf L} (1997) 323--379.

\ite{Lang}
{\sc  Serge Lang}:
{\it Real and functional analysis}, third edition,
Springer,
1993.

\ite{LMdp}
{\sc Martin Lazar, Darko Mitrovi\'c}:
{\sl Velocity averaging --- a general framework},
{\it  Dyn.~Partial Differ.~Equ.}
{\bf 9} (2012) 239--260.

\ite{LZ}
{\sc Martin Lazar, Enrique Zuazua}: 
{\sl Averaged control and observation of parameter-depending wave equations},  
{\it  C. R. Acad.~Sci.~Paris, Ser.~I}
{\bf 352} (2014) 497--502.

\ite{LP}
{\sc Pierre Louis Lions, Thierry Paul}:
{\sl Sur les mesures de Wigner,} 
{\it Revista Mat.~Iberoamericana } {\bf 9} (1993) 553--618.

\ite{MMip}
{\sc Marin Mi\v{s}ur, Darko Mitrovi\'c}:
{\sl On a generalisation of compensated compactness in the $L^p-L^q$ setting},
{\it J. Functional Analysis}
{\bf 268} (2015) 1904-1927.

\ite{EJParma}
{\sc Evgenij Jurjevi\v c Panov}:
{\sl Existence and strong pre-compactness properties for entropy solutions
of a first-order quasilinear equation with discontinuous flux},
{\it Arch.~Ration.~Mech.~Anal.}
{\bf 195} (2010) 643--673.

\ite{EJPahp}
{\sc Evgenij Jurjevi\v c Panov}:
{\sl Ultra-parabolic H-measures and compensated compactness},
{\it Ann.~Inst.~H. Poincar\'e Anal.~Non Lin\'eaire}
{\bf 28} (2011) 47--62.

\ite{FRin}
{\sc Filip Rindler}: 
{\sl Directional oscillations, concentrations, and compensated compactness via microlocal compactness forms},
{\it Arch.~Ration.~Mech.~Anal.}
{\bf 215} (2015) 1--63.

\ite{Tprse}
{\sc Luc Tartar}:
{\sl H-measures, a new approach for studying homogenisation, oscillations and concentration effects in partial differential equations, }
{\it Proc.~Roy.~Soc.~Edinburgh}\/
{\bf 115A} (1990) 193--230.


\ite{Tgth}
{\sc Luc Tartar}:
{\it The general theory of homogenization: A personalized introduction,}
Springer,
2009.

\ite{TmsHm}
{\sc Luc Tartar}:
{\sl Multi-scale H-measures, }
{\it Discrete and Continuous Dynamical Systems, S}\/
{\bf 8} (2015) 77--90.

\ite{Zh}
{\sc Ping Zhang}:
{\it Wigner measure and semiclassical limits of nonlinear Schr\" odinger
	equations},
AMS, 2008.

\bye